%% file: 23wz1l.tex
\documentclass[12pt]{amsproc}

\usepackage{amsmath}
\usepackage{amsthm}
\usepackage{amscd}
\usepackage{amssymb}
\usepackage{url}
\usepackage{nicefrac}
\usepackage[pdftex]{graphicx}%
\usepackage[outercaption]{sidecap}
\usepackage{amsfonts,epigraph}%
\usepackage{booktabs} 
\usepackage[all,cmtip]{xy}
\newtheorem{theorem}{Theorem}[section]
\newtheorem{corollary}[theorem]{Corollary}
\newtheorem{lemma}[theorem]{Lemma}
\newtheorem{proposition}[theorem]{Proposition}

\theoremstyle{definition}
\newtheorem{definition}[theorem]{Definition}

\newtheorem{choices-notations}[theorem]{Choices and Notations}

\theoremstyle{remark}
\newtheorem{data}{Data\!\!}

\newcommand{\exendproof}{\renewcommand{\qed}{\relax}\end{proof}}

\newsavebox{\SmallMathBox}

\DeclareRobustCommand*{\nicefrac}[2]{\ifmmode\mathnicefrac{#1}
{ #2}%
  \else\textnicefrac{#1}{#2}\fi}
\newcommand*{\textnicefrac}[2]{\check@mathfonts%
\mbox{\raisebox{.5ex}{\fontsize\sf@size\z@\selectfont#1}\kern-.
1em%
/\kern-.1em\raisebox{- .25ex}{\fontsize\sf@size\z@\selectfont#2} }}
\newcommand*{\mathnicefrac}[2]{%
  \mathchoice
    {\m@fr@c{\scriptstyle}{#1}{#2}}
    {\m@fr@c{\scriptstyle}{#1}{#2}}
    {\m@fr@c{\scriptscriptstyle}{#1}{#2}}
    {\m@fr@c{\scriptscriptstyle}{#1}{#2}}}

\def\lla{\langle}

\def\rra{\rangle}
\def\sqm1{\sqrt{-1}}

\def\too{\longrightarrow}

\def\wt{\widetilde}
\def\x{\times}

\def\={\cong}
\def\>{\supset}
\def\<{\subset}

\def\12{\frac{1}{2}}
\def\0{^{\circ}}

\def\CC{{\mathbb C}}

\def\KK{{\mathbb K}}

\def\NN{{\mathbb N}}

\def\RR{{\mathbb R}}

\def\ZZ{{\mathbb Z}}

\def\Bb{{\mathcal B}}
\def\Cc{{\mathcal C}}

\def\Ff{{\mathcal F}}

\def\Ii{{\mathcal I}}

\def\Ll{{\mathcal L}}

\def\Pp{{\mathcal P}}

\def\Ss{{\mathcal S}}

\def\CLR{{\mathcal{CLR}}}

\newcommand{\mi}{\mathrm i}

\def\C{\CC}

\def\la{\lambda}

\def\N{\NN}

\def\R{\RR}

\def\w{\omega}

\def\Z{\ZZ}

\newcommand{\defeq}{\mathrel{\mathop:}=}

 \DeclareMathOperator{\dist}{dist}

\DeclareMathOperator{\Gr}{Gr}
\DeclareMathOperator{\dom}{dom}

\def\flg#1{\Ff\Ll_{#1}}

\DeclareMathOperator{\GL}{GL} \DeclareMathOperator{\gl}{gl} \DeclareMathOperator{\Graph}{graph}
 
\DeclareMathOperator{\Hom}{Hom}

 \DeclareMathOperator{\image}{im}
\DeclareMathOperator{\Index}{index}

\DeclareMathOperator{\Mas}{Mas}

 \DeclareMathOperator{\ran}{im} 
  
 \DeclareMathOperator{\romS}{S} 
\DeclareMathOperator{\sign}{sign}

 \DeclareMathOperator{\Sp}{Sp}

\begin{document}

\title{A formula of local Maslov index and applications}
%
%
%
%


\author{Li Wu}
\address{Department of Mathematics, Shandong University, Jinan, Shandong 250100, P. R. China, https://orcid.org/0000-0003-2605-3824} \email{vvvli@sdu.edu.cn}

\author{Chaofeng Zhu}
\address{Chern Institute of Mathematics and Laboratory of Pure Mathematics and Combinatorics (LPMC), Nankai University,
	Tianjin 300071, P. R. China, https://orcid.org/0000-0003-4600-4253} \email{zhucf@nankai.edu.cn}
\thanks{Corresponding author: C. Zhu [\texttt{zhucf@nankai.edu.cn}]\\	
Chaofeng Zhu is supported by National Key R\&D Program of China (2020YFA0713300), NSFC,  China Grants (11971245, 11771331), Nankai Zhide Foundation and Nankai University.}

\date{today}




\begin{abstract}
In this paper, we explicitly express the local Maslov index by a Maslov index in finite dimensional case without symplectic reduction. Then we calculate the Maslov index for the path of pairs of Lagrangian subspaces in triangular form. In particular, we get the Maslov-type index of a given symplectic path in triangle form. As applications, we calculate the splitting numbers of the symplectic matrix in triangle form, dependence of iteration theory on triangular frames and mod 2 Maslov-type index for a real symplectic path. We study the continuity of families of bounded linear relations and families of bounded linear operators acting on closed linear subspaces as technique preparations. 
\end{abstract}

\maketitle

\include{ls1}

\include{s7}

\include{applications}

\include{operator-space}


\bibliography{Hamiltonian}
\bibliographystyle{amsplain-jl}

\end{document}

%% file: ls1.tex
\section{Introduction}\label{s:introduction}

Since the legendary work of V.P. Maslov \cite{Maslov} in the mid 1960s and the supplementary
explanations by V. Arnol'd \cite{Ar67}, there has been a continuing
interest in the Maslov index for curves of pairs of Lagrangian subspaces in symplectic space. As explained by Maslov and
Arnol'd, the interest arises from the study of dynamical systems in
classical mechanics and related problems in Morse theory. In the finite dimensional case, the notion of Maslov index was generalized by J. Robbin and D. Salamon \cite{RoSa93} in 1993, and S. E. Cappell, R. Lee and E. Y. Miller \cite{CaLeMi94} in 1994. In the infinite dimensional case, it was  generalized by R. C. Swanson \cite{Sw78a} in 1978, B. Booss and K. Furutani \cite{BoFu98} in 1998, and B. Booss and the second author \cite{BoZh18} in 2018. Under a nondegenerate condition, the Maslov index was calculated by J. J. Duistermaat \cite[Lemma 2.5]{Du76}. The calculation of the Maslov index was systematically studied in \cite[Section 3.2]{BoZh18}.

Throughout this paper, we denote by $\KK$ the field of real numbers or complex numbers. We denote by $\N$, $\Z$, $\R$ and $\C$ the sets of all natural, integral, real and complex numbers respectively. By $S^1$ we denote the unit circle in the complex plane. We denote by $\mi$ the imaginary unit. We denote by $I_X$ the identity map on a set $X$. If there is no confusion, we will omit the subindex $X$. For two vector spaces $X$ and $Y$ over $\KK$, we denote by $\Hom(X,Y)$ the set of linear maps from $X$ to $Y$ and $\dim X$ the dimension of $X$ respectively. 

For a Banach space $X$, we denote by $\Bb(X)$ the set of bounded linear operators. We equip the the set of closed linear subspaces $\Ss(X)$ with {\em gap distance} $\hat\delta$ (\cite[Section IV.2.1]{Ka95}). We denote the Fredholm index for two linear subspaces $(M,N)$ of $X$ by $\Index(M,N;X)$. We denote by $\Index(M,N)=\Index(M,N;X)$ if there is no confusion. 

Let $V$ be a vector space over $\KK$ and $Q\colon V\times V\to\KK$ be a symmetric form.
Then $m^{\pm}(Q)$ and $m^0(Q)$ denotes the {\em Morse positive (or negative) index} and the {\em nullity} of $Q$ respectively.  

Let $X$ be a complex vector space. A mapping
\[
  \omega\colon X\times X\too \C
\]
is called a {\em symplectic form} on $X$, if it is a
non-degenerate skew symmetric form. Then we call $(X,\w)$ a {\em symplectic vector space}. A {\em symplectic Banach space} $(X,\omega)$ is a Banach space $X$ together with a symplectic form $\w$ such that $|\w(x,y)|\le C\|x\|\|y\|$.
The {\em $\omega$-annihilator} of a subspace $\la$ of $X$ is defined by
\[\la^{\omega}\ :=\{x\in X;\omega(x,y)=0\text{ for all }y\in\la\}.\]
A subspace ${\la}$ of is called {\em symplectic},
{\em isotropic}, {\em coisotropic}, or {\em Lagrangian} if
\[
\la\cap{\la}^{\w}\ =\ \{0\}\,,\quad{\la} \,\<\, {\la}^{\w}\,,\quad {\la}\,\>\, {\la}^{\w}\,,\quad
{\la}\,\ =\ \, {\la}^{\w}\,,
\]
respectively.
The {\em Lagrangian Grassmannian}
$\Ll(X,\w)$ consists of all Lagrangian subspaces of
$(X,\w)$. We write $\Ll(X)\ :=\Ll(X,\w)$ if there is no confusion.

Let $X$ be a vector space with two linear subspaces $\la,\mu$. The pair $(\la,\mu)$ is said to be {\em Fredholm} if
$\dim(\la\cap\mu)<+\infty$ and $\dim X/(\la+\mu)<+\infty$. The {\em Fredholm index} of $(\la,\mu)$ is defined by
\begin{equation}\label{e:fp-index}
\Index(\la,\mu)\ :=\ \dim(\la\cap\mu) - \dim X/(\la+\mu).
\end{equation}

The space of \emph{Fredholm pairs} of Lagrangian subspaces of a symplectic vector space
$(X,\omega)$ is defined by
\begin{equation}\label{e:fp-lag-alg}
\Ff\Ll(X)\ :=\ \{(\lambda,\mu)\in\Ll(X)\x\Ll(X);\; (\la,\mu)\text{ is Fredholm}\}
\end{equation}
For $k\in\Z$ we define
\begin{equation}\label{e:fl-index-k}
\Ff\Ll_k(X)\ :=\ \{(\lambda,\mu)\in\Ff\Ll(X); \Index(\la,\mu)=k\}.
\end{equation}

Let $a<b$ be two real numbers. Let $X$ be a Banach space with continuously varying symplectic structures $\omega(s))$, $s\in[a,b]$. Let $\la(s),\mu(s)$, $s\in[a,b]$ be two paths of closed subspaces of $X$ such that $(\la(s),\mu(s))\in\Ff\Ll_0(X)$. The Maslov indices (cf. \cite[Definition 3.1.4]{BoZh18}) $\Mas\{\la(s),\mu(s)\}=\Mas_{+}\{\la(s),\mu(s)\}$ and $\Mas_{-}\{\la(s),\mu(s)\}$ are well-defined integers.

Our first main result explicitly express the local Maslov index by a Maslov index in finite dimensional case without symplectic reduction.

\begin{theorem}\label{t:localization}Let $\varepsilon>0$ be a positive number. Let $X$ be a (complex) Banach space with continuously varying symplectic structure $\w(s)$, $s\in(-\varepsilon,\varepsilon)$. Let $(\lambda(s),\mu(s))$, $s\in(-\varepsilon,\varepsilon)$ be a path of Fredholm pairs of Lagrangian subspaces of $(X,\omega(s))$ of index $0$. Let $V(s)$, $s\in(-\varepsilon,\varepsilon)$ be a path of finite dimensional linear subspaces of $X$ such that $X=V(0)\oplus(\lambda(0)+\mu(0))$. Denote by
\begin{align}\label{e:lambdas}\la_0(s):&=\la(s)\cap(V(s)+\mu(s)),\; \la_1(s):=V(s)^{\omega(s)}\cap\la(s),\\
\label{e:mus}\mu_0(s):&=\mu(s)\cap(V(s)+\la(s)),\; \mu_1(s):=V(s)^{\omega(s)}\cap\mu(s),\\
\label{e:xs}X_0(s):&=(V(s)+\la(s))\cap(V(s)+\mu(s)),\;X_1(s):=\la_1(s)+\mu_1(s).
\end{align}
Then there exists a $\delta\in (0,\varepsilon)$ such that the following hold for all $s\in(-\delta,\delta)$ and $[s_1,s_2]\subset(-\delta,\delta)$.
\newline (a) $X_0(s)=X_1(s)^{\omega(s)}$, $X_1(s)=X_0(s)^{\omega(s)}$, and $X_0(s)$, $X_1(s)$ are continuously varying on $s$.
\newline (b) $X_0(s)=V(s)\oplus\la_0(s)=V(s)\oplus\mu_0(s)$, $X_1(s)=\la_1(s)\oplus\mu_1(s)$, $\dim V(0)=\dim V(s)=\dim\la_0(s)=\dim\mu_0(s)=(\dim X_0(s))/2$, and $X=X_0(s)\oplus X_1(s)$.
\newline (c) $\la(s)=\la_0(s)\oplus\la_1(s)$, $\mu(s)=\mu_0(s)\oplus\mu_1(s)$, $\la_0(s)$ and $\mu_0(s)$ are paths of Lagrangian subspaces of $(X_0(s),\w(s)|_{X_0(s)})$, $\la_1(s)$ and $\mu_1(s)$ are paths of Lagrangian subspaces of $(X_1(s),\w(s)|_{X_1(s)})$.
\newline (d) We have
\begin{equation}\label{e:local1}
\Mas_{\pm}\{\la(s),\mu(s);\;s\in[s_1,s_2]\}=\Mas_{\pm}\{\la_0(s),\mu_0(s);\;s\in[s_1,s_2]\}.
\end{equation}
\newline (e) Assume that there is a path $W(s)$, $s\in(-\delta,\delta)$ of Lagrangian subspaces of $(X_0(s),\w(s)|_{X_0(s)})$ such that $W(s)\cap\la_0(s)=W(s)\cap\mu_0(s)=\{0\}$. Then there is an operator $A(s):\mu_0(s)\rightarrow W(s)$ such that $\lambda_0(s)=\Graph(A(s))$.
Define $Q(s)(x,y)\;:=\omega(x,A(s)y)$ with $x,y\in \mu_0(s)$. Then each form $Q(s)$ is symmetric and we have
\begin{align}
\label{e:local2}\Mas\{\la(s),\mu(s);\;s\in[s_1,s_2]\}&=m^+(Q(s_2))-m^+(Q(s_1)),\\
\label{e:local3}\Mas_-\{\la(s),\mu(s);\;s\in[s_1,s_2]\}&=m^-(Q(s_1))-m^-(Q(s_2)),\\
\label{e:local4}\dim(\la(s)\cap\mu(s))&=m^0(Q(s)),
\end{align}
where $m^{\pm}(Q)$, $m^0(Q)$ denote the positive (negative) Morse index and the nullity of $Q$
respectively for a symmetric form $Q$.
\end{theorem}

Our second main result calculate the Maslov index for the path of pairs of Lagrangian subspaces in triangular form.

Let $(X,\w)$ be a symplectic vector space with three isotropic subspaces $\alpha$, $\beta$, $\gamma$. Then we can define the symmetric form $Q(\alpha,\beta;\gamma)$ on $\gamma\cap(\alpha+\beta)$ (cf. \cite[(2.3)]{Du76}) by
\begin{equation}\label{e:Q-form} Q(\alpha,\beta;\gamma)(z_1,z_2)=\omega(x_1,y_2)
\end{equation}
for all $z_j=x_j+y_j\in\gamma$, $x_j\in\alpha$, $y_j\in\beta$, $j=1,2$.

\begin{theorem}\label{t:main-triangle}
	Let $Z$ be a Banach space with continuously varying symplectic structures $\w(s)$, $s\in[0,1]$.
	Assume that $Z=X(s)\oplus Y(s)$ with two paths $\{X(s)\in\Ll(Z,\w(s));\;s\in[0,1]\}$ and $\{Y(s)\in\Ll(Z,\w(s));\;s\in[0,1]\}$.
	Let $\{(\lambda(s),\mu(s))\in\Ff\Ll_0(Z,\w(s)); \;s\in [0,1]\}$ be a path. Set $A(s)\;:=\lambda(s)\cap Y(s)$, $B(s)\;:=\mu(s)\cap Y(s)$,
	$\gamma(s)\;:=A(s)^{\omega(s)}\cap X(s)\oplus A(s)$, $\delta(s)\;:=B(s)^{\omega(s)}\cap X(s)\oplus B(s)$. Assume that $\{A(s)\;s\in[0,1]\}$, $\{B(s)\;s\in[0,1]\}$
	are two paths, and $\dim Y(s)/(A(s)+B(s))<+\infty$ for each $s\in[0,1]$. For each $s\in[0,1]$, we define $Q(s)\;:=\omega(s)(x_1,(\lambda(s)-\mu(s))x_2)$
	with $x_1,x_2\in \dom(\lambda_s)\cap\dom(\mu(s))$.
	Then we have
	\begin{align}\label{e:main-triangle}
		\begin{aligned}
			\begin{split}
				\Mas_{\pm}&\{\la(s),\mu(s);\;\w(s),s\in[0,1]\}=\pm\dim(\lambda(0)\cap\mu(0)\cap Y(0))\mp\\&\dim(\lambda(1)\cap\mu(1)\cap Y(1))\pm m^{\pm}(Q(1))\mp m^{\pm}(Q(0)).
			\end{split}
		\end{aligned}    
	\end{align}
	Moreover, we have
	\begin{align}\label{e:maslov-nullity}
		\dim(\lambda(s)\cap\mu(s))=\dim\ker Q(s)+\dim A(s)\cap B(s)
	\end{align}
	for each $s\in[0,1]$.
\end{theorem}

In order to prove Theorem \ref{t:main-triangle}, we firstly prove a structure theorem of pairs of Lagrangian subspaces in triangular form (Theorem \ref{t:triangle-pairs} below). Then we study the continuity of families of bounded linear relations and families of bounded linear operators acting on closed linear subspaces as technique preparations (Theorem \ref{t:operator-space} below). Then we can reach our result. The proof can be simplified greatly in the finite dimensional case.

A special case of Theorem \ref{t:main-triangle} is the formula for Maslov-type index of symplectic path in triangular form with diagonal boundary value conditions ( Theorem \ref{t:Maslov-type-index-triangular-form} below). Then we apply Theorem \ref{t:Maslov-type-index-triangular-form} and get the formula of splitting number of symplectic matrix in triangular form. By \cite[Lemma 5.5]{Gi10}, it has no loss of generality if we consider the splitting number of a real symplectic matrix at $1$. Finally we study the dependence of iteration theory on triangular frames (Theorem \ref{t:iteration-frame} below).

The paper is organized as follows. In Section 1, we explain why we make the research and state the main results of the paper. In Section 2, we systematically study the $\omega$-closed linear subspaces of symplectic vector spaces. In Section 3, we explicitly express the local Maslov index by a Maslov index in finite dimensional case without symplectic reduction and prove Theorem \ref{t:localization}. In Section 4, we study Fredholm pairs of index $0$ of Lagrangian subspaces in triangular form and prove Theorem \ref{t:triangle-pairs}. In Section 5, we calculate the Maslov index for the path of pairs of Lagrangian subspaces in triangular form and prove Theorem \ref{t:main-triangle}. In Section 6, we recall the notion of Maslov-type index and prove Theorems \ref{t:Maslov-type-index-triangular-form}, \ref{t:splitting-numbers-up-triangular}, \ref{t:iteration-frame} and \ref{t:Maslov-type-index-mod-2}. In Appendix A, we study the continuity of families of bounded linear relations and families of bounded linear operators acting on closed linear subspaces and prove Theorems \ref{t:closed-domain} and \ref{t:operator-space}.

We would like to thank the referees of this paper for their critical reading and
very helpful comments and suggestions.

%% file: s7.tex
\section {The $\omega$-closed linear subspaces}\label{s:omega-closed}
In this section we discuss the $\omega$-closed linear subspaces of a symplectic vector space.

\begin{definition}\label{d:omega-closed}
Let $(X,\omega)$ be a symplectic vector space.
A linear subspace $V$ of $X$ is called {\em $\omega$-closed}, if $V^{\omega\omega}=V$.
\end{definition}

Let $(X,\omega)$ be a symplectic Banach space. Let $V$ be an $\omega$-closed linear subspace. Then $V$ is closed.

By \cite[(1.3)]{BoZh18}, we have

\begin{lemma}\label{l:triple-omega}
Let $(X,\omega)$ be a symplectic vector space.
Let $V$ be a linear subspace of X.
Then $V^\omega$ is $\omega$-closed, i.e. $V=V^{\w\w\w}$.
\end{lemma}

\begin{lemma}\label{l:finite-extension-closed}
Let $(X,\omega)$ be a symplectic vector space.
Let $V$ be an $\omega$-closed linear subspace of $X$.
Let $V_1\supset V$ be a linear subspace of $X$ such that $\dim V_1/V <+\infty$.
Then there holds $\dim V_1/V=\dim V^\omega/V_1^\omega$ and $V_1$ is $\omega$-closed.
Specially, if $V$ is a finite dimensional linear subspace of $X$, $V$ is $\omega$-closed.
\end{lemma}

\begin{proof} Our conditions imply that $V^{\w\w}\cap V_1=V$ and $V_1^{\w\w\w}\cap V^{\w}=V_1^{\w}$.
By \cite[Lemma 1.1.2.b]{BoZh18} we have
$$\dim V_1/V=\dim V^\omega/V_1^\omega=\dim V_1^{\w\w}/V^{\w\w}=\dim V_1^{\w\w}/V.$$
It follows that $V_1^{\w\w}=V_1$ and $V_1$ is $\w$-closed.
\end{proof}

\begin{lemma}\label{l:intersection-omega}
 Let $(X,\omega)$ be a symplectic vector space.
 Let $U_i,i\in I$ be a family of  $\omega$-closed linear subspaces of $X$. Then
 we have
\newline (a)$(\sum_{i\in I}U_i)^\omega =\bigcap_{i\in I}U_i^\omega$,
\newline (b) $\bigcap_{i\in I}U_i$ is $\omega$-closed, and
\newline (c) $(\bigcap_{i\in I}U_i)^\omega=\sum_{i\in I}U_i^\omega$  if and only if  $\sum_{i\in I}U_i^\omega$ is $\omega$-closed.
\end{lemma}

\begin{proof}
 It is trivial that $(\sum_{i\in I}U_i)^\omega=\bigcap_{i\in I}U_i^{\omega}$ and $(\bigcap_{i\in I})^{\omega}\supset \sum_{i\in I}U_i^\omega $.
We have
\begin{align*}
 &(\sum_{i\in I}U_i^\omega)^\omega=\bigcap_{i\in I} U_i^{\omega\omega}=\bigcap_{i\in I}U_i,\text{ and }\\
 &(\sum_{i\in I}U_i^\omega)^{\omega\omega}=(\bigcap_{i\in I}U_i)^\omega.
\end{align*}
So $\bigcap_{i\in I}U_i$ is $\omega$-closed, and we have $(\bigcap_{i\in I}U_i)^\omega=\sum_{i\in I}U_i^\omega$  if and only if  $(\sum_{i\in I}U_i^\omega)$ is $\omega$-closed.
\end{proof}

\begin{lemma}\label{l:closed-sum}
 Let $(Z,\omega)$ be a symplectic vector space such that $Z=X+ Y=X^\omega+ Y^\omega$. Then the following hold.
 \newline (a) We have $Z=X\oplus Y=X^\omega\oplus Y^\omega$, $X=X^{\w\w}$ and $Y=Y^{\w\w}$.
 \newline (b) Let $V\subset X$, $W\subset Y$ be two $\w$-closed linear subspaces of $Z$. Then $V+W$ is an $\omega$-closed subspace of $Z$.
 \end{lemma}

\begin{proof} (a)
We have $X^\omega\cap Y^\omega=(X+Y)^\omega=Z^\omega=\{0\}$ and
$X\cap Y\subset X^{\omega\omega}\cap Y\subset X^{\omega\omega}\cap Y^{\omega\omega}=(X^\omega+Y^\omega)^\omega=\{0\}$.
It follows that
\[Z=X\oplus Y=X^\omega\oplus Y^\omega. \]
By \cite[Lemma A.1.1]{BoZh18}, we have
\[X^{\w\w}=X^{\w\w}\cap(X+Y)=X+X^{\w\w}\cap Y=X.\]
Similarly we have $Y^{\w\w}=Y$.
\newline (b) Set $V':=(V+Y)^{\omega\omega}\cap X$. Since $Y\subset V+Y\subset (V+Y)^{\omega\omega} $, by \cite[Lemma A.1.1]{BoZh18} we have 
\[(V+Y)^{\omega\omega}=(V+Y)^{\omega\omega}\cap(X+Y)=V'\oplus Y .
\]
By Lemma \ref{l:intersection-omega}, we have $V'^{\omega\omega}=V'$.
By Lemma \ref{l:triple-omega}, we have
\begin{equation}
 \begin{cases}
(V+Y)^{\omega\omega\omega}=(V'+Y)^\omega =V'^\omega\cap Y^\omega\\
(V+Y)^{\omega\omega\omega}=(V+Y)^\omega=V^\omega\cap Y^\omega
 \end{cases}.
\end{equation}
It follows that $V'^\omega\cap Y^\omega=V^\omega\cap Y^\omega$.
Since $V',V\subset X$, we have $V\subset V^\prime$ and $X^\omega\subset V'^\omega\subset V^\omega$.
By \cite[Lemma A.1.1]{BoZh18} it follows that
\[
 V'^\omega=(V'^\omega\cap Y^\omega)\oplus X^\omega=(V^\omega\cap Y^\omega)\oplus X^\omega=V^\omega.
\]
Then we have
\[
 V'=V'^{\omega\omega}=V^{\omega\omega}=V.
\]
It follows that
\[
 V+Y=V'+Y=(V+Y)^{\omega\omega}.
\]
By exchanging the role of $(X,V)$ and $(Y, W)$, we have $X+W=(X+W)^{\w\w}$. Since $V+W=(V+Y)\cap(X+W)$, by Lemma \ref{l:intersection-omega}, $V+W$ is $\w$-closed.

\end{proof}

\begin{lemma}\label{l:closed-sum2}
 Let $(X,\omega)$ be a symplectic vector space.
 Let $\lambda,\mu\in \Ll(X)$ be such that $(\lambda,\mu)\in \flg{0}(X)$.
 Let $V$ be a linear subspace of $\lambda$ and $W$ be a linear subspace of $\mu$. Assume that $V, W$ are $\omega$-closed.
 Then $V+W$ is $\omega$-closed.
\end{lemma}

\begin{proof}
 By \cite[Corollary 1.3.4]{BoZh18}, there is an $\alpha\in \Ll(X)$ such that
 \[\dim \mu/(\mu\cap\alpha)<+\infty,\text{ and }\lambda\oplus\alpha=X.\]
By Lemma \ref{l:intersection-omega}, $W\cap\alpha$ is $\omega$-closed.
Since $X=\lambda\oplus\alpha$, by Lemma \ref{l:closed-sum}, $V+\alpha\cap W$ is $\omega$-closed.
Since $W\supset \alpha\cap W$, there is a linear subspace $U\subset W$ such that $W=\alpha\cap W\oplus U$. Then we have
\[V+W=V+\alpha\cap W+U.\]
It follows that
\begin{align*}
 \dim (V+W)/(V+\alpha\cap W)&\leq \dim W/(\alpha\cap W)=\dim (W+\alpha)/\alpha\\
 &\leq \dim(\mu+\alpha)/\alpha=\dim\mu/(\mu\cap\alpha)\\
 &<+\infty.
\end{align*}
By Lemma \ref{l:finite-extension-closed}, $V+W$ is an $\omega$-closed linear subspace.
\end{proof}

\section{A formula of local Maslov index without symplectic reduction}
In this section we shall give a explicit construction of the path $\la_0(s)$ in \cite[Theorem 3.2.1]{BoZh18}. Then we get a very nice formula to calculate the local Maslov index.

Firstly we recall the notion of linear relations.  A {\em linear relation} $A$ between two linear spaces $X$ and $Y$ is a linear subspace of $X\times Y$. For the detailed concept of linear relation, we use Cross \cite{Cr98} and G. W. Whitehead \cite[B.2]{GWh78} as references.

For a linear relation $A\subset X\times Y$, we define the {\em domain}, the {\em range}, the {\em kernel} and the {\em indeterminacy} of $A$ are
\begin{equation*}
	\begin{split}
		\dom(A)&\ := \{x\in X;  \text{ there is a }y\in Y
		\text{ such that } (x,y)\in A\}, \\
		\ran A&\ :=\{y\in Y;  \text{ there is a }x\in X \text{ such that } (x,y)\in A\},\\
		\ker A&\ :=\{x\in X; \;(x,0)\in A\},\\
		A0&\ :=\{y\in Y; \;(0,y)\in A\}.
	\end{split}
\end{equation*}
respectively.

We have some algebraic observations.

\begin{lemma}\label{l:dim-linear-relation}
	Let $A\subset X\times Y$ be a  linear relation. Then we have
	\begin{align}\label{e:dim-linear-relation}
		\dim A=\dim A0+\dim\dom(A)=\dim\ker A+\dim\image A. 
	\end{align}
\end{lemma}

\begin{proof}
	The linear relation $A$ induces a linear map
	$\tilde A=A/(\{0\}\times A0):\dom(A)\to Y$. 
	Then the map $\varphi:\dom(A)\to\tilde A$ defined by $\varphi(x)=(x,\tilde Ax)$ is a linear isomorphism.
	Thus we have
	\begin{align*}
		\dim\dom(A)=\dim\tilde A=\dim A-\dim A0
	\end{align*}
	if $\dim A<+\infty$, and $\dim A=+\infty$ if and only if $\dim A0+\dim\dom(A)=+\infty$. 
	Then we obtain 
	\begin{align*}
		\dim A=\dim A0+\dim\dom(A)=\dim A^{-1}=\dim\ker A+\dim\image A. 
	\end{align*}
\end{proof}

\begin{lemma}\label{l:algebra} Let $V$ be a linear space with three linear subspaces $V_1$, $V_2$, $V_3$. Then we have
\newline (a) $(V_1+V_2)\cap(V_1+V_3)=V_1+V_2\cap(V_1+V_3)=V_1+V_3\cap(V_1+V_2)$, and
\newline (b) The linear relation 
\begin{align}
	A\ :=\{(v_2,v_3)\in V_2\times V_3;\; v_2+v_3\in V_1\}\subset V_2\times V_3
\end{align} induces a linear isomorphism
\begin{equation}\label{e:relation-iso}\tilde A: \frac{V_2\cap(V_1+V_3)}{V_1\cap V_2}\to\frac{V_3\cap(V_1+V_2)}{V_1\cap V_3}.
\end{equation}
\newline (c) Denote by $B:=A+\{0\}\times V_2$ to be a liinear subspace of $V_2\times(V_3+V_2)$. Then the linear relation $B\subset V_2\times(V_3+V_2)$ induces a linear isomorphism
\begin{equation}\label{e:relation-iso2}\tilde B: \frac{V_2\cap(V_1+V_3)}{V_1\cap V_2+V_2\cap V_3}\to\frac{V_3\cap(V_1+V_2)+V_2}{V_1\cap V_3+V_2}.
\end{equation}
\newline (d) If $\dim V_1+\dim V_2\cap V_3<+\infty$, we have $\dim (V_2\cap(V_1+V_3))+\dim(V_1\cap V_3)=\dim V_1\cap(V_2+V_3)+\dim(V_2\cap V_3)=\dim V_1-\dim(V_1+V_2+V_3)/(V_2+V_3)+\dim(V_2\cap V_3)$.
\end{lemma}

\begin{proof} (a) Since $V_1\subset V_1+V_2$ and $V_1\subset V_1+V_3$, (a) follows from \cite[Lemma A.1.1]{BoZh18}.
\newline (b) The linear relation $A$ induces a well-defined surjective linear map $\dom(A)\to\ran(A)/\ker A^{-1}$. So the induced map
\begin{equation}\label{e:relation-isomorphism}
\tilde A:\dom(A)/\ker A\to\ran A/\ker A^{-1}
\end{equation}
is a linear isomorphism. Since we have
\begin{align*}
&\dom(A)=V_2\cap(V_1+V_3),\quad\ran A=V_3\cap(V_1+V_2),\\
&\ker A=V_1\cap V_2,\text{ and }\ker A^{-1}=V_1\cap V_3,
\end{align*}
by \eqref{e:relation-isomorphism} we obtain (b).
\newline (c) Since we have
\begin{align*}
	&\dom(B)=V_2\cap(V_1+V_3),\quad\ran B=V_3\cap(V_1+V_2)+V_2,\\
	&\ker B=V_1\cap V_2+V_2\cap V_3,\text{ and }\ker B^{-1}=V_1\cap V_3+V_2,
\end{align*}
by \eqref{e:relation-isomorphism} we obtain (c).
\newline (d) Since $(V_1+V_2+V_3)/(V_2+V_3)\cong V_1/(V_1\cap(V_2+V_3))$, (d) follows from (b).
\end{proof}

\begin{definition}\label{d:finite-perturb} Let $X$ be a vector space and $M$, $N$ be linear subspaces of $X$.
We define $M\sim^f N$ if $\dim M/(M\cap N), \dim N/(M\cap N)<+\infty$, and call $N$ a
{\em finite change} of $M$ (see \cite[p. 273]{Ne68}). In this case we define the {\em
relative dimension} $[M-N]:=\dim M/(M\cap N)-\dim (N/M\cap N)$.
\end{definition}

\begin{lemma}\label{l:index-finite-perturb} Let $X$ be a vector space and $M_1$, $M_2$, $N_1$, $N_2$ be linear subspaces of $X$.
Assume that $M_1\sim^f M_2$ and $N_1\sim^f N_2$. Then we have
\begin{equation}\label{e:index-finite-perturb}\Index(M_1,N_1)=\Index(M_2,N_2)+[M_1-M_2]+[N_1-N_2]
\end{equation}
if one of the two sides is well-defined.
\end{lemma}

\begin{proof} Firstly we consider the case that $M_1\supset M_2$ and $N_1=N_2$. It is \cite[Problem IV.4.6]{Ka95}. Here we give a complete proof.

Let $V$ be a vector subspace of $M_1$ such that $M_1=V\oplus M_2$.
Then $\dim V=[M_1-M_2]<+\infty$. By Lemma \ref{l:algebra} c., we have
\begin{align*}
\dim &(N_1\cap(V+M_2))+\dim (V\cap M_2)=\dim V\\
&-\dim (V+M_2+N_1)/(M_2+N_1)+\dim(M_2\cap N_1),\text{ and}\\
\dim &(M_1\cap N_1)=[M_1-M_2]-\dim (M_1+N_1)/(M_2+N_1)\\
&+\dim(M_2\cap N_1).
\end{align*}
Since $(X/(M_2+N_1))/((M_1+N_1)/(M_2+N_1))\cong X/(M_1+N_1)$, our result follows.

For the general case, we have
\begin{align*}
\Index(M_1,N_1)&=\Index(M_1\cap M_2,N_1)+[M_1-M_1\cap M_2]\\
&=\Index(M_2,N_1)+[M_1\cap M_2-M_2]+[M_1-M_1\cap M_2]\\
&=\Index(M_2,N_1)+[M_1-M_2],
\end{align*}
and
\begin{align*}
\Index(M_1,N_1)&=\Index(M_2,N_1)+[M_1-M_2]\\
&=\Index(M_2,N_2)+[M_1-M_2]+[N_1-N_2].
\end{align*}
\end{proof}

We have the following lemma (cf. \cite[Lemma 1.4.9]{BoZh18}).

\begin{lemma}\label{l:V-omega-la} Let $(X,\omega)$ be a symplectic vector space. Let $\la$ be a Lagrangian subspace of $X$ and $V$ be a finite
dimensional subspace of $X$. Assume that $V\cap\la=\{0\}$. Then we have $\dim\la/(V^{\w}\cap\la)=\dim V$ and $V^{\w}+\la=X$.
\end{lemma}

\begin{proof} Since $\la=\la^{\w}$, by Lemma \ref{l:finite-extension-closed} we have $V+\la$ is $\w$-closed, and
\[\dim\la/(V^{\w}\cap\la)=\dim (V^{\w}\cap\la)^{\w}/\la=\dim(V+\la)/\la=\dim V.\]
By \cite[Lemma 1.3.2]{BoZh18}, we have $V^{\w}+\la=X$.
\end{proof}

Firstly we discuss some equivalent conditions.

\begin{lemma}\label{l:equivalent-conditions}
Let $(X,\omega)$ be a symplectic vector space with a pair $(\la,\mu)\in\Ff\Ll_0(X)$. Let $V$ be a finite dimensional subspace of $X$ such
that $V\cap\la=V\cap\mu=\{0\}$. Then the following four conditions are equivalent:
\begin{itemize}
\item[(i)] $V+V^{\w}\cap\la+\mu=X$,
\item[(ii)] $V^{\w}\cap(V+\la)\cap\mu=\{0\}$,
\item[(iii)] $V+\la+V^{\w}\cap\mu=X$, and
\item[(iv)] $V^{\w}\cap\la\cap(V+\mu)=\{0\}$.
\end{itemize}
\end{lemma}

\begin{proof} Clearly we have (i)$\Rightarrow$(ii) and (iii)$\Rightarrow$(iv).

By Lemma \ref{l:index-finite-perturb} and Lemma \ref{l:V-omega-la}, we have
\[\Index(V+\la,V^{\w}\cap\mu)=\Index(\la,\mu)+[(V+\la)-\la]+[V^{\w}\cap\mu-\mu]=0.\]
So we have (ii)$\Leftrightarrow$(iii). Similarly we have (i)$\Leftrightarrow$(iv).

Combine the implications together, we have (i)$\Rightarrow$(ii)$\Leftrightarrow$(iii)$\Rightarrow$(iv)$\Leftrightarrow$(i).
Our results then follows.
\end{proof}

By \cite[Proposition 1.3.3]{BoZh18} we have

\begin{lemma}\label{l:condition} Let $(X,\omega)$ be a symplectic vector space with a pair $(\la,\mu)\in\Ff\Ll_0(X)$. Let $V$ be a finite dimensional subspace of $X$. Assume that $X=V\oplus(\la+\mu)$.
Then we have
\begin{equation}\label{e:condition}
X=V\oplus V^{\w}\cap\la\oplus\mu=V\oplus\la\oplus V^{\w}\cap\mu.
\end{equation}
\end{lemma}

Our key observation is the following (cf. \cite[Proposition 1.3.3]{BoZh18}).

\begin{proposition}\label{p:decomposition} Let $(X,\omega)$ be a symplectic vector space with a pair $(\la,\mu)\in\Ff\Ll_0(X)$. Let $V$ be a finite dimensional linear subspace of $X$. Denote by
\begin{align}\label{e:lambda}\la_0\ :&=\la\cap(V+\mu),\; \la_1\ :=V^{\omega}\cap\la,\\
\label{e:mu}\mu_0\ :&=\mu\cap(V+\la),\; \mu_1\ :=V^{\omega}\cap\mu,\\
\label{e:x}X_0\ :&=(V+\la)\cap(V+\mu),\;X_1\ :=\la_1+\mu_1.
\end{align}
Assume that $X=V+\la_1+\mu$ and $V\cap\la=V\cap\mu=\{0\}$.
Then the following hold.
\newline (a) $X_0=X_1^{\omega}$, $X_1=X_0^{\omega}$.
\newline (b) $X_0=V\oplus\la_0=V\oplus\mu_0$, $X_1=\la_1\oplus\mu_1$, $\dim V=\dim\la_0=\dim\mu_0=(\dim X_0)/2$, and $X=X_0\oplus X_1$.
\newline (c) $\la=\la_0\oplus\la_1$, $\mu=\mu_0\oplus\mu_1$, $\la_0$ and $\mu_0$ are Lagrangian subspaces of $X_0$, $\la_1$ and $\mu_1$ are Lagrangian subspaces of $X_1$.
\end{proposition}

\begin{proof} By Lemma \ref{l:triple-omega}, $V^{\w}$ is $\w$-closed. Since $\dim V<+\infty$, by Lemma \ref{l:finite-extension-closed}, $V$, $V+\la$ and $V+\mu$ are $\w$-closed, and $\dim\la/\la_1=\dim\mu/\mu_1=\dim V$. By Lemma \ref{l:intersection-omega}, $\la_0$, $\la_1$, $\mu_0$, $\mu_1$, $X_0$ are $\w$-closed. By Lemma \ref{l:closed-sum2}, $\la_1+\mu$, $\la+\mu_1$, $X_1$ are $\w$-closed.
\newline (a) Direct computaion.
\newline (b) By Lemma \ref{l:algebra}.a, we have $X_0=V+\la_0=V+\mu_0$. Since $V\cap\la=V\cap\mu=\{0\}$, we have $V\cap\la_0=V\cap\mu_0=\{0\}$ and $X_0=V\oplus\la_0=V\oplus\mu_0$.

Since $V+\la+\mu=X$, we have $\la_1\cap\mu_1=V^{\w}\cap\la\cap\mu=\{0\}$ and $X_1=\la_1\oplus\mu_1$.

Since $V+\la+\mu=X$ and $\Index(\la,\mu)=0$, by Lemma \ref{l:algebra}.d, we have $\dim V=\dim\la_0=\dim\mu_0$. Since $X_0=V\oplus\la_0$, we have $\dim V=(\dim X_0)/2$.

By Lemma \ref{l:finite-extension-closed}, $X_0+X_1$ is $\w$-closed. Since $X=V+\la_1+\mu$ and $V\cap\la=V\cap\mu=\{0\}$, by Lemma \ref{l:equivalent-conditions} and \cite[Lemma A.1.1]{BoZh18}, we have
\begin{eqnarray*}(X_0+X_1)^{\w}&=&X_1\cap X_0=(\la_1+\mu_1)\cap(V+\la)\cap(V+\mu)\\
&=&(\la_1+\mu_1\cap(V+\la))\cap(V+\mu)\\
&=&\la_1\cap(V+\mu)=\{0\}.
\end{eqnarray*}
So we have $X=X_0\oplus X_1$.
\newline (c) By Lemma \ref{l:equivalent-conditions} we have $\la_0\cap\la_1=\la\cap(V+\mu)\cap V^{\w}\cap\la=\{0\}$. Clearly we have $\la\supset \la_0+\la_1$. Since $\dim\la/\la_1=\dim V=\dim\la_0$, we have $\la=\la_0\oplus\la_1$. Similarly we have $\mu=\mu_0\oplus\mu_1$.

Since $\la_0$, $\mu_0$ are isotropic subspaces of $X_0$ and $\dim\la_0=\dim\mu_0=(\dim X_0)/2$, we have $\la_0,\mu_0\in\Ll(X_0)$.

Since $X_1=\la_1\oplus\mu_1$ and $\la_1$, $\mu_1$ are isotropic, by \cite[Proposition 1]{BoZh13} we have $\la_1,\mu_1\in\Ll(X_1)$.
\end{proof}

With the above preparations, we can now prove Theorem \ref{t:localization}.

\begin{proof}[Proof of Theorem \ref{t:localization}] (a), (b), (c) Since $\dim V(s)<+\infty$ and $V(s)$ is a path, by \cite[Lemma I.4.10]{Ka95}, we have $\dim V(s)=\dim V(0)$. Since $V(0)\cap\la(0)=\{0\}$, by \cite[Lemma III.1.9]{Ka95}, $V(0)+\la(0)$ is closed. By \cite[Corollary A.3.14]{BoZh18}, there exists a $\delta_1\in(0,\varepsilon)$ such that $V(s)\cap\la(s)=V(s)\cap\mu(s)=\{0\}$ for each $s\in(-\delta_1,\delta_1)$ and $\{V(s)+\la(s);\;s\in(-\delta_1,\delta_1)\}$ is a path.

By \cite[Lemma 3.1.1]{BoZh18}, $\{V(s)^{\w(s)};\;s\in(-\varepsilon,\varepsilon)\}$ is a path. By Lemma \ref{l:V-omega-la}, $V(s)^{\w(s)}+\la(s)=V(s)^{\w(s)}+\mu(s)=X$ for each $s\in(-\delta_1,\delta_1)$. By \cite[Corollary A.3.14]{BoZh18}, $\{\la_1(s);\;s\in(-\delta_1,\delta_1)\}$ and $\{\mu_1(s);\;s\in(-\delta_1,\delta_1)\}$ are two paths. Since $V(0)+\la_1(0)+\mu(0)=X$, by \cite[Proposition A.3.5.c]{BoZh18}, there exists a $\delta\in(0,\delta_1)$ such that $V(s)+\la_1(s)+\mu(s)=X$.

By Proposition \ref{p:decomposition}, (a), (b), (c) follows except the continuity of $X_1(s)$ on $s$, which follows from \cite[Corollary A.3.14]{BoZh18}.

(d) follows from (c) and \cite[Proposition 2.3.1.c, Theorem 3.1.5]{BoZh18}.

(e) follows from (d) and \cite[Lemma 2.3.2]{BoZh18}. Note that there always exists a Lagrangian subspace $W(0)$ of $X_0(0)$ such that $W(0)\cap \la_0(0)=W(0)\cap\mu_0(0)=\{0\}$. Then the existence of the path $W(s)$ for $|s|\ll 1$ follows from \cite[Lemma 1.3.6]{BoZh18},
the existence of $W(0)$, and \cite[Proposition A.3.5.c]{BoZh18}.
\end{proof}

\section{Fredholm pairs of Lagrangian subspaces}\label{s:Fredhom-pair}
In this section we study Fredholm pairs of index $0$ of Lagrangian subspaces in triangular form.

Firstly we study the Fredhlom pairs of Lagrangian subspaces in diagonal forms.

\begin{definition}\label{d:diagonal} Let $X$ be a vector space with a direct sum decomposition $X=\bigoplus_{i\in I}X_i$. A linear subspace $\lambda$ is called to be in {\em diagonal form} with respect to the decomposition, if $\lambda=\bigoplus_{i\in I}(\lambda\cap X_i)$.
\end{definition}

The subspaces in diagonal forms has the following nice properties.

\begin{lemma}\label{l:diagonal-form-properties}
	Let $X$ be a vector space with a direct sum decomposition $X=\bigoplus_{i\in I}X_i$. Let $\{\lambda_j;\;j\in J\}$ is a family of linear subspaces in diagonal form with respect to the decomposition. Then $\sum_{j\in J}\lambda_j$ and $\bigcap_{j\in J}\lambda_j$ are in diagonal form with respect to the decomposition, and we have $(\sum_{j\in J}\lambda_j)\cap X_i=\sum_{j\in J}\lambda_j\cap X_i$.	
\end{lemma}

\begin{proof}
	We have 
	\begin{align*}
		\sum_{j\in J}\lambda_j=&\sum_{j\in J}\sum_{i\in I}\lambda_j\cap X_i=\bigoplus_{i\in I}\sum_{j\in J}\lambda_j\cap X_i,\\
		\bigcap_{j\in J}\lambda_j=&\bigcap_{j\in J}\sum_{i\in I}\lambda_j\cap X_i=\bigoplus_{i\in I}\bigcap_{j\in J}\lambda_j\cap X_i.
	\end{align*}
	Then our lemma follows.
\end{proof}

Then we have the following.

\begin{corollary}\label{c:diagonal-index}
	Let $X$ be a vector space with a direct sum decomposition $X=\bigoplus_{i\in I}X_i$. Let $\lambda$, $\mu$ be two linear subspaces in diagonal form with respect to the decomposition. Then we have
	\begin{align}\label{e:diagonal-index}
		\Index(\lambda,\mu)=\sum_{i\in I}\Index(\lambda\cap X_i,\mu\cap X_i;X_i),
	\end{align}
	where $\Index(M,N;Y)$ denotes the Fredholm index of the pair $(M,N)$ in $Y$.
\end{corollary}

\begin{lemma}\label{l:diag-omega} Let $(Z,\omega)$ be a symplectic vector space with a linear subspace $V$. Assume that $Z=X\oplus Y$ with $X,Y\in \Ll(Z)$. If $V=V\cap X+V\cap Y$, we have 
	\begin{align*}
		V^{\w}=&(V\cap Y)^{\w}\cap X+(V\cap X)^{\w}\cap Y,\\
		V^{\w}\cap X=&(V\cap Y)^{\w}\cap X,\\
		V^{\w}\cap Y=&(V\cap X)^{\w}\cap Y,\\
		V^{\w}+X=&X+(V\cap X)^{\w}\cap Y=(V\cap X)^{\w},\\
		V^{\w}+Y=&(V\cap Y)^{\w}\cap X+Y=(V\cap Y)^{\w}.
	\end{align*}
\end{lemma}

\begin{proof} Direct computation shows that  \[V^{\w}\cap X=(V\cap X)^{\w}\cap(V\cap Y)^{\w}\cap X=(V\cap Y)^{\w}\cap X.\]
Similarly we have $V^{\w}\cap Y=(V\cap X)^{\w}\cap Y$. By \cite[Lemma A.1.1]{BoZh18}, we have 
\begin{equation*}
	\begin{split}
        V^{\w}&=(V\cap X)^{\w}\cap(V\cap Y)^{\w}\cap(X+Y)\\
        &=(V\cap X)^{\w}\cap((V\cap Y)^{\w}\cap X+Y)\\
        &=(V\cap Y)^{\w}\cap X+(V\cap X)^{\w}\cap Y.
    \end{split}
\end{equation*}
By Lemma \ref{l:closed-sum},  $V^{\w}+X=X+(V\cap X)^{\w}\cap Y$ and $V^{\w}+Y=(V\cap Y)^{\w}\cap X+Y$ are closed.
So we have $V^{\w}+X=(V\cap X)^{\w}$ and $V^{\w}+Y=(V\cap Y)^{\w}$.
\end{proof}

We then have the following criterion of Lagrangian subspace in diagonal form.

\begin{corollary}
	\label{l:diag-Lagrangian} Let $(Z,\omega)$ be a symplectic vector space with a linear subspace $V$. Assume that $Z=X\oplus Y$ with $X,Y\in \Ll(Z)$ and $V=V\cap X+V\cap Y$. Then $V$ is a Lagrangian subspace of $Z$ if and only if $V^\w\cap X=V\cap Y$ and $V^\w\cap Y=V\cap X$.
\end{corollary}

We have the following formula of the dimension of intersection between a linear subspace and a Lagrangian subspace in diagonal form.

\begin{lemma}\label{l:dim-triangular-diagonal}
	Let $(Z,\omega)$ be a symplectic vector space such that $Z=X\oplus Y$ with $X,Y\in \Ll(Z)$. We view a linear subspace of $Z$ as a linear relation from $X$ to $Y$. Let $\lambda$ be a linear subspace. Let $\mu\in\Ll(Z)$ be a Lagrangian subspace in diagonal form. Then we have
	\begin{align}\label{etriangular-diagonal-domain}
		\dom(\lambda\cap\mu)=&\mu\cap X\cap\lambda^{-1}(\mu\cap X)^\omega,\\
		\label{e:triangular-diagonal-image}
		\image(\lambda\cap\mu)=&\mu\cap Y\cap\lambda(\mu\cap Y)^\omega,
		\\
		\label{e:dim-triangular-diagonal-domain}
		\dim(\lambda\cap\mu)=&\dim(\mu\cap X\cap\lambda^{-1}(\mu\cap X)^\omega)+\dim(\lambda\cap\mu\cap Y)\\
		\label{e:dim-triangular-diagonal-image}
		=&\dim(\mu\cap Y\cap\lambda(\mu\cap Y)^\omega)+\dim(\lambda\cap\mu\cap X).
	\end{align}
	Especially if $\lambda$ is isotropic and $\dom(\lambda)\supset\mu\cap X$, we have $\mu\cap X\cap\lambda^{-1}(\mu\cap X)^\omega=\ker Q(X,Y;\lambda)|_{\mu\cap X}$.
	If $\lambda$ is isotropic and $\image(\lambda)\supset\mu\cap Y$, we have $\mu\cap Y\cap\lambda^{-1}(\mu\cap Y)^\omega=\ker Q(Y,X;\lambda)|_{\mu\cap Y}$. 
\end{lemma}

\begin{proof} 
	Since $\mu\in\Ll(Z)$ is a Lagrangian subspace in diagonal form, by Lemma \ref{l:diag-omega}, we have $\mu\cap Y=(\mu\cap X)^\omega\cap Y$. Then we have
	\begin{align*}
		\dom(\lambda\cap\mu)
		=&\{x\in\mu\cap X;\;\exists y\in\mu\cap Y\text{ such that }x+y\in \lambda\}\\
		=&\{x\in\mu\cap X;\;\exists y\in(\mu\cap X)^\omega\cap Y\text{ such that }x+y\in \lambda\}\\		
		=&\mu\cap X\cap\lambda^{-1}(\mu\cap X)^\omega,\\
		\image(\lambda\cap\mu)=&\dom(\lambda^{-1}\cap\mu^{-1})=\mu\cap Y\cap\lambda(\mu\cap Y)^\omega.
	\end{align*}
	By Lemma \ref{l:dim-linear-relation}, we obtain \eqref{e:dim-triangular-diagonal-domain} and \eqref{e:dim-triangular-diagonal-image}.
	By definition, we obtain the two special cases.
\end{proof}

\begin{lemma}\label{l:diag-decomposition} Let $(Z,\omega)$ be a symplectic vector space. Assume that $Z=X\oplus Y$ with $X,Y\in \Ll(Z)$.
Let $(\alpha,\beta)$ be in $\Ff\Ll_0(Z)$ such that $\alpha=\alpha\cap X\oplus \alpha\cap Y$ and $\beta=\beta\cap X\oplus \beta\cap Y$.
Then there is an isotropic subspace $V$ and symplectic subspaces $Z_0$, $Z_1$ of $Z$ such that the following hold.
\newline (a) $V=V\cap X\oplus V\cap Y$, $V^{\w}=V^{\w}\cap X\oplus V^{\w}\cap Y$ and $Z=V\oplus(\alpha+\beta)$.
\newline (b) $Z= Z_0\oplus Z_1$, $Z_0=Z_1^{\w}$, $Z_1=Z_0^{\w}$, $Z_0=V\oplus \alpha\cap \beta$, $Z_1=V^\omega\cap \alpha\oplus V^\omega\cap \beta$.
\newline (c) Set $\beta':=V+V^\omega\cap \beta$. Then we have $Z=\alpha\oplus\beta'$, $\beta'\in \Ll(Z)$, and $\beta'=\beta'\cap X\oplus\beta'\cap Y$.
\end{lemma}

\begin{proof} Since $(\alpha,\beta)\in\Ff\Ll_0(Z)$ and $\alpha$, $\beta$ are in diagonal form, there are finite dimensional subspaces $W_1$ of $X$ and $W_2$ of $Y$ respectively such that
\[X=W_1\oplus(\alpha\cap X+\beta\cap X),\quad Y=W_2\oplus(\alpha\cap Y+\beta\cap Y).\]
Set $W\ :=W_1+W_2$. Then we have $Z=W\oplus(\alpha+\beta)$. Set $Z_0\ :=W+\alpha\cap \beta$ and $Z_1\ :=W^\omega\cap \alpha+W^\omega\cap \beta$.
By Lemma \ref{l:diag-omega} and \cite[Proposition 1.3.3]{BoZh18}, (a)--(c) holds if we replace $V$ by $W$ and $\dim Z_0=2\dim W=2\dim(\alpha\cap\beta)$. Moreover, $Z_0$ and $Z_1$ are symplectic. Note that $\alpha\cap\beta=\alpha\cap\beta\cap X\oplus \alpha\cap\beta\cap Y$ and $Z_0=Z_0\cap X\oplus Z_0\cap Y$. By \cite[Proposition 1]{BoZh13}, $Z_0\cap X$ and $Z_0\cap Y$ are Lagrangian subspaces of $Z_0$.

By taking bases in $Z_0\cap X$ and $Z_0\cap Y$, we can assume that
\[(\omega|_{Z_0})(u,v)=\lla f(x_1),y_2\rra-\overline{\lla f(x_2),y_1\rra}\]
for each $u=x_1+y_1$, $v=x_2+y_2$, $x_1,x_2\in Z_0\cap X$, $y_1,y_2\in Z_0\cap Y$ under a linear isomorphism $f:Z_0\cap X\to Z_0\cap Y$, where $\lla\cdot,\cdot\rra$ is the standard inner product.
\newline (a) Since $\alpha\cap\beta\in\Ll(Z_0)$, we have $f(\alpha\cap\beta\cap X)=(\alpha\cap\beta\cap Y)^{\bot_2}$, where $\bot_2$ denotes the orthogonal complement in $Z_0\cap Y$. Set $V\ :=f^{-1}(\alpha\cap\beta\cap Y)\oplus f(\alpha\cap\beta\cap X)$. Then we have $V=V\cap X\oplus V\cap Y\subset Z_0\cap V^\w$. By Lemma \ref{l:diag-omega}, we have $V^{\w}=V^{\w}\cap X\oplus V^{\w}\cap Y$. 
By definition, we have $V\in\Ll(Z_0)$ and $Z_0=V\oplus\alpha\cap\beta$. By \cite[Proposition 1.3.3]{BoZh18}, we have
\[Z=Z_0\oplus Z_1=V\oplus\alpha\cap\beta\oplus Z_1=V\oplus(\alpha+\beta).\]
\newline (b) Since $V\subset Z_0$, we have 
\[V^{\w}\supset Z_0^{\w}=Z_1\supset W^{\w}\cap\alpha.\]
Note that 
\[V^{\w}\cap\alpha\cap\beta=(V+\alpha+\beta)^{\w}=\{0\}.\]
By \cite[Lemma A.1.1]{BoZh18} and \cite[Proposition 1.3.3]{BoZh18}, there holds that
\[V^{\w}\cap\alpha=V^{\w}\cap(\alpha\cap\beta+W^{\w}\cap\alpha)=W^{\w}\cap\alpha.\]
Similarly we have $V^{\w}\cap\beta=W^{\w}\cap\beta$. Hence (b) follows from \cite[Proposition 1.3.3]{BoZh18}.
\newline (c) By \cite[Proposition 1.3.3]{BoZh18}, we have $Z=\alpha\oplus\beta'$, $V\in\Ll(Z_0)$ and $V^\w\cap\beta\in\Ll(Z_1)$. Then we have $\beta'\in\Ll(Z)$.
By (a), we have $\beta'=\beta'\cap X\oplus\beta'\cap Y$.
\end{proof}

By Proposition \ref{p:decomposition} and the above proof, we have

\begin{corollary}\label{c:diag-decomposition}Let $(Z,\omega)$ be a symplectic vector space. Assume that $Z=X\oplus Y$ with $X,Y\in \Ll(Z)$.
Let $(\alpha,\beta)$ be in $\Ff\Ll_0(Z)$ such that $\alpha=\alpha\cap X\oplus \alpha\cap Y$ and $\beta=\beta\cap X\oplus \beta\cap Y$. Let $V$ be an isotropic subspace of $X$ such that $X=V+\alpha+V^{\omega}\cap\beta$, $V\cap\alpha=V\cap\beta=\{0\}$ and $V=V\cap X\oplus V\cap Y$. Then the following hold.
\newline (a) $Z= Z_0\oplus Z_1$, $Z_0=Z_1^{\w}$, $Z_1=Z_0^{\w}$, $Z_0=V\oplus \alpha\cap(V+\beta)=V\oplus\beta\cap(V+\alpha)$, $Z_1=V^\omega\cap \alpha\oplus V^\omega\cap \beta$.
\newline (b) Set $\beta':=V+V^\omega\cap \beta$. Then we have $Z=\alpha\oplus\beta'$, $\beta'\in \Ll(Z)$, and $\beta'=\beta'\cap X\oplus\beta'\cap Y$.
\end{corollary}

We recall the notion of symplectic reduction (cf. \cite[Definition 1.4.1]{BoZh18}.

\begin{definition}\label{d:sympl-red}
Let $(X,\omega)$ be a symplectic vector space with a coisotropic subspace $W$.
\newline (a) The space $W/W^{\omega}$ is a symplectic vector space with induced symplectic structure
\begin{equation}\label{e:red-structure}
\wt{\omega}(x+W^{\omega},y+W^{\omega})\ :=\ \omega(x,y) \text{ for all $x,y\in W$}.
\end{equation}
We call $(W/W^{\omega}, \wt{\omega})$ the {\em
symplectic reduction} of $X$ via $W$.
\newline (b) Let $\lambda$ be a linear subspace of $X$.
The {\em symplectic reduction} of $\lambda$ via $W$ is defined by%
\begin{equation}\label{e:red-subspace}
R_W(\lambda)=R_W^{\w}(\la)\ :=\ \bigl((\lambda+W^{\omega})\cap W\bigr)/W^{\omega}\ =\ \bigl(\lambda\cap  W+W^{\omega}\bigr)/W^{\omega}.%
\end{equation}
We denote by $\pi_w\ :=R_{w^{\w}}$ for an $\w$-closed isotropic subspace $w$.
\end{definition}

Clearly, $R_W(\lambda)$ is isotropic if $\lambda$ is isotropic.

The following lemma generalizes \cite[Lemma 1.4.6]{BoZh18}.

\begin{lemma}\label{l:g-index-reduction}
	Let $X$ be a vector space with  four linear subspaces $W_1$, $W_2$, $\lambda$, $\mu$. Assume that $W_1\subset W_2$.
	For each linear subspace $V$, set $R(V)\defeq (V\cap W_2+W_1)/W_1=((V+W_1)\cap W_2)/W_1$.
	Then we have the following.
	\newline (a) There are linear isomorphisms
	\begin{align}\label{e:reduction-iso-cap1}
		\frac{R(\lambda)\cap R(\mu)}{R(\lambda\cap\mu)}&\cong \frac{(\lambda+\mu\cap W_2)\cap W_1}{\lambda\cap W_1+\mu\cap W_1},\\
		\label{e:reduction-iso-cap2}
		R(\lambda\cap \mu)&\cong \frac{\lambda\cap\mu\cap W_2}{\lambda\cap \mu\cap W_1}.
	\end{align}
	In particular, we have $R(\lambda\cap\mu)=R(\lambda)\cap R(\mu)$ if $W_1\subset\mu\subset W_2$.
	\newline (b) There are linear isomorphisms
	\begin{align}\label{e:reduction-iso-plus1}
		\frac{R(\lambda+\mu)}{R(\lambda)+ R(\mu)}&\cong \frac{(\lambda+\mu)\cap W_2+W_1}{\lambda\cap W_2+\mu\cap W_2+W_1},\\
		\label{e:reduction-iso-plus2}
	    \frac{R(W_2)}{R(\lambda+\mu)}&\cong \frac{\lambda+\mu+W_2}{\lambda+\mu+W_1}.
	\end{align}
	In particular, we have $R(\lambda+\mu)=R(\lambda)+R(\mu)$ if $W_1\subset\mu\subset W_2$.
\end{lemma}

\begin{proof} (a) Since $W_1\subset W_2$, by  \cite[Lemma A.1.1]{BoZh18}, we have
	\begin{align*}
		R(\lambda)\cap R(\mu)&=((\lambda+W_1)\cap(\mu+W_1)\cap W_2)/W_1\\
		&=((\lambda\cap(\mu+W_1)+W_1)\cap W_2)/W_1\\
		&=(\lambda\cap(\mu+W_1)\cap W_2+W_1)/W_1\\
		&=(\lambda\cap(\mu\cap W_2+W_1)+W_1)/W_1.
	\end{align*}
	If $W_1\subset\mu\subset W_2$, we have $\mu\cap W_2+W_1=\mu$ and $R(\lambda\cap\mu)=R(\lambda)\cap R(\mu)$.
	
    By Lemma \ref{l:algebra}.c, there is a linear isomorphism
    \begin{align*}
    	f:\frac{(\lambda+\mu\cap W_2)\cap W_1}{\lambda\cap W_1+\mu\cap W_1}\to\frac{\lambda\cap(\mu\cap W_2+W_1)+ W_1}{\lambda\cap\mu\cap W_2+W_1}.
    \end{align*}
    Then \eqref{e:reduction-iso-cap1} follows.
    
    Since $W_1\subset W_2$, we have
    \begin{align*}
    	R(\lambda\cap\mu)&=\frac{\lambda\cap\mu\cap W_2+W_1}{W_1}\\
    	&\cong\frac{\lambda\cap\mu\cap W_2}{\lambda\cap \mu\cap W_1}.
    \end{align*}
	\newline (b) By definition we have \eqref{e:reduction-iso-plus1}.
	Since $W_1\subset W_2$, we have
	\begin{align*}
		\frac{R(W_2)}{R(\lambda+\mu)}
		&\cong \frac{W_2}{(\lambda+\mu+W_1)\cap W_2}\\
		&\cong\frac{\lambda+\mu+W_2}{\lambda+\mu+W_1}.			
	\end{align*}
	Then \eqref{e:reduction-iso-plus2} follows.
	
	If $W_1\subset\mu\subset W_2$, by \cite[Lemma A.1.1]{BoZh18} and \eqref{e:reduction-iso-plus1}, we have $(\lambda+\mu)\cap W_2=\lambda\cap W_2+\mu=\lambda\cap W_2+\mu\cap W_2$ and $R(\lambda+\mu)=R(\lambda)+R(\mu)$.
\end{proof}

We have the following algebraic fact.

\begin{lemma}\label{l:separate-space} Let $X$ be a vector space with four linear subspaces $W_1$, $W_2$, $\lambda$, $\mu$. Assume that $W_1\supset\lambda$ and $W_2\supset\mu$. Set $W:=W_1\cap W_2$. Then we have $(\lambda+\mu)\cap W=\lambda\cap W+\mu\cap W$.
\end{lemma}

\begin{proof} By \cite[Lemma A.1.1]{BoZh18}, we have
\begin{align*}
(\lambda+\mu)\cap W &=(\lambda+\mu)\cap W_1\cap W_2=(\lambda+\mu\cap W_1)\cap W_2\\
 &=\lambda\cap W_2+\mu\cap W_1=\lambda\cap W+\mu\cap W.
\end{align*}
\end{proof}

Then we have the following corollary (cf. \cite[Proposition 1.4.8]{BoZh18}.

\begin{corollary}\label{c:Lagrange-preserved}
Let $(X,\omega)$ be a symplectic vector space with $(\lambda,\mu)\in\Ff\Ll_0(X)$ and two linear subspaces $W_1$, $W_2$.
Assume that $W:=W_1\cap W_2$ is a coisotropic subspace of $X$, $W_1\supset\lambda$ and $W_2\supset\mu$.
Then $R_W(\lambda),R_W(\mu)$ are both Lagrangian subspaces of $W/W^{\w}$ and there hold
\[\Index(R_W(\lambda),R_W(\mu))=0 ,\quad \lambda\cap W^{\omega}+\mu\cap W^{\omega}=W^{\omega}.\]
\end{corollary}

\begin{proof} 
Since $(\lambda,\mu)\in\Ff\Ll_0(X)$, by \cite[Lemma 1.2.8]{BoZh18}, we have $\lambda+\mu=(\lambda\cap\mu)^{\w}$. Since $W_1\supset \lambda$ and $W_2\supset \mu$, we have $W_1^\omega\subset \lambda$ and $W_2^\omega\subset \mu$. Then we have $W\supset \lambda\cap\mu$ and $W^{\w}\subset\lambda+\mu$. 
By Lemma \ref{l:separate-space}, we have \[(\lambda+\mu)\cap W=\lambda\cap W+\mu\cap W.\]
Then we have
\begin{align*}
	W^{\omega}&\supset(\lambda+\mu\cap W)\cap W^{\omega}\\
	&\supset(\lambda\cap W+\mu\cap W)\cap W^{\omega}\\
	&=(\lambda+\mu)\cap W\cap W^{\omega}=W^{\omega}.
\end{align*}
Thus we obtain
\[(\lambda+\mu\cap W)\cap W^{\omega}=W^{\omega}.\]
Since $R_W(\lambda)$, $R_W(\mu)$ are isotropic subspaces of $W/W^{\omega}$, by Lemma \ref{l:g-index-reduction} and \cite[Lemma1.1.2]{BoZh18}, we have
\begin{align*}
0\ge &\Index(R_W(\lambda),R_W(\mu))\\
\ge&\Index(\lambda,\mu)-\dim(\lambda\cap\mu\cap W^{\w})+\dim X/(\lambda+\mu+W)\ge 0.
\end{align*}
It follows that $\Index(R_W(\lambda),R_W(\mu))=0$. By Lemma \ref{l:g-index-reduction}, we have 
\begin{align*}
	\lambda\cap W^{\omega}+\mu\cap W^{\omega}=(\lambda+\mu\cap W)\cap W^{\omega}=W^{\omega}.
\end{align*}
By \cite[Proposition 1]{BoZh13}, $R_W(\lambda),R_W(\mu)$ are both Lagrangian subspaces of $W/W^{\w}$.
\end{proof}

We have the following algebraic observation. 

\begin{lemma}\label{l:reduction-separate} Let $Z$ be a vector space with five linear subspaces $X$, $Y$, $W_1$, $W_2$ and $\lambda$ such that $Z=X\oplus Y$ and $W_1\subset W_2$. Define $R(V)\ :=(V\cap W_2+W_1)/W_1$ for each linear subspace $V$ of $Z$. Assume that $W_2=W_2\cap X+W_2\cap Y$ and $\lambda=\lambda\cap X+\lambda\cap Y$. Then we have $\lambda\cap W_2=\lambda\cap W_2\cap X\oplus\lambda\cap W_2\cap Y$ and $R(\lambda)=R(\lambda\cap X)\oplus R(\lambda\cap Y)$.
\end{lemma}

\begin{proof} Clearly we have
\begin{align*}
\lambda\cap W_2&=(\lambda\cap X)\cap(W_2\cap X)\oplus(\lambda\cap Y)\cap (W_2\cap Y)\\
&=\lambda\cap W_2\cap X\oplus\lambda\cap W_2\cap Y.
\end{align*}
By definition, we have $R(\lambda)=R(\lambda\cap X)\oplus R(\lambda\cap Y)$.
\end{proof}

We make some preparations for the proof of Theorem \ref{t:triangle-pairs} below.

\begin{lemma}\label{l:new-Lag} Let $(X,\omega)$ be a symplectic vector space with a Lagrangian subspace $\lambda$ and an $\omega$-closed isotropic subspace $V$. Assume that $V+\lambda$ is $\omega$-closed. Then $V+V^{\omega}\cap\lambda$ is a Lagrangian subspace of $X$.
\end{lemma}

\begin{proof} By \cite[Lemma A.1.1]{BoZh18}, we have
\[(V+V^{\omega}\cap\lambda)^{\omega}=V^{\omega}\cap(V+\lambda)=V+V^{\omega}\cap \lambda.\]
\end{proof}

The following lemma shows that under natural conditions, a Lagrangian subspace is in triangular form.

\begin{lemma}\label{l:Lag-to-triangle}
Let $(Z,\omega)$ be a symplectic vector space such that $Z=X\oplus Y$ with $X,Y\in \Ll(Z)$. We view a linear subspace of $Z$ as a linear relation from $X$ to $Y$. Let $\lambda \in\Ll(Z)$ be a Lagrangian subspace. Set $W\;:=\lambda\cap Y$. Assume that $\lambda+Y$ is $\omega$-closed. For each $t\in\R$, we define the linear operator $P(t):Z\to Z$ by $P(t)(x+y)=x+ty$ for each $x\in X$ and $y\in Y$, and the subspaces $\alpha(t)$ of $Z$ by $\alpha(t)\;:=P(t)(\lambda)+W$.
Then for each $t\in\R$ we have $\alpha(0)=W^\omega\cap X+W$, $\alpha(1)=\lambda$, $\alpha(t)\in \Ll(Z)$, $\dom(\alpha(t))=W^\omega\cap X$, and $\alpha(t)\cap Y=W$. We call $\alpha(0)$ the {\em diagonal part} of $\lambda$. 
\end{lemma}

\begin{proof}
Note that 
\begin{align*}
	W^{\omega}&=\lambda+Y,\text{ and }\\
	\dom(\lambda)&=(\lambda+Y)\cap X=(\lambda\cap Y)^\omega\cap X=W^\omega\cap X.
\end{align*}
Then we have 
\begin{align*}
	\alpha(0)&=P(0)\lambda+W=\dom(\lambda)+W=W^\omega\cap X+W,\\
	\alpha(1)&=P(1)\lambda+W=\lambda+W=\lambda, \text{ and }\\
	\dom(\alpha(t))&=\dom(P(t)(\lambda))=\dom(\lambda)=W^\omega\cap X.
\end{align*}
Note that $P(t)(\lambda)\cap Y=W$. By Lemma \cite[A.1.1]{BoZh18}, we have 
\begin{align*}
	\alpha(t)\cap Y&=(P(t)(\lambda)+W)\cap Y=W,\text{ and }\\
	\lambda+Y&=(\lambda+Y)\cap(X+Y)=\dom(\lambda)+Y.
\end{align*}
Then there holds that
\begin{align}\label{e: W-dom}
	W=\lambda\cap Y=(\dom(\lambda))^{\omega}\cap Y.
\end{align}

Now we prove that $\alpha(t)\in\Ll(Z)$ for $t\in\R$. Note that we have
\begin{align*}
	0=\omega(x+y_2,x_1+y_1)=\omega(x,y_1)+\omega(y_2,x_1)
\end{align*}
for $x+y_2,x_1+y_1\in\lambda$ with $x,x_1\in X$ and $y_1,y_2\in Y$.
By our assumptions and \eqref{e: W-dom}, we have 
\begin{align*}
	\begin{aligned}
		\alpha(t)^{\omega}=&(P(t)(\lambda))^{\omega}\cap(\lambda+Y)\\
		=&\left\{x+y;\;x\in \dom(\lambda), y\in Y,\omega(x+y,x_1+ty_1)=0\right.\\
		&\left.\text{ for all }x_1\in \dom(\lambda),y_1\in Y, x_1+y_1\in\lambda\right\}\\
		=&\left\{x+y;\;x\in \dom(\lambda), y\in Y,
		\exists y_2\in Y\right.\\
		&\left.\text{ such that } x+y_2\in\lambda,y-ty_2\in W\right\}\\
		=&P(t)(\lambda)+W=\alpha(t).
	\end{aligned}
\end{align*}
\end{proof}

\begin{lemma}\label{l:compare-subspaces} Let $X$ be a vector space with four linear subspaces $V_1$, $V_2$, $W_1$, $W_2$. Assume that $V_1\subset W_1$, $V_2\subset W_2$, $V_1+V_2=W_1+W_2$, $V_1\cap V_2=W_1\cap W_2$. Then we have $V_1=W_1$ and $V_2=W_2$.
\end{lemma}

\begin{proof} Set $V:=V_1\cap V_2=W_1\cap W_2$. By our assumptions, we have
\[V_1/V\subset W_1/V,\; V_2/V\subset W_2/V,\text{ and }V_1/V\oplus V_2/V=W_1/V\oplus W_2/V.\]
It follows that $V_1/V=W_1/V$, $V_2/V=W_2/V$. Hence there hold $V_1=W_1$ and $V_2=W_2$.
\end{proof}

\begin{lemma}\label{l:space-distribution}
	Let $Z$ be a linear space with linear subspaces $X$ and $Y$ such that $Z=X+Y$. Let $\{V_i\}_{i\in A}$ are linear subspaces of $Z$ with $V_i\supset Y$ for each $i\in A$. Then we have 
    \begin{align}\label{e:space-distribution}
    	\sum_{i\in A}V_i=\sum_{i\in A}V_i\cap X+Y.
    \end{align}
\end{lemma}

\begin{proof}
	By \cite[Lemma A.1.1]{BoZh18}, for each $i\in A$ we have 
	\[V_i=V_i\cap(X+Y)=V_i\cap X+Y.\] 
	Then \eqref{e:space-distribution} follows.
\end{proof}

Now we can prove the main result of this section.

\begin{theorem} \label{t:triangle-pairs}
Let $(Z,\omega)$ be a symplectic vector space with $X,Y\in \Ll(Z)$ such that $Z=X\oplus Y$.
Let $(\alpha,\beta)$ be in $\Ff\Ll_0(X)$. We view $\alpha,\beta$ as linear relations from $X$ to $Y$.
Set $W_1\;:=\alpha\cap Y$, $W_2:=\beta\cap Y$, $\gamma\;:=W_1^\omega\cap X+W_1$, and $\delta\;:=W_2^\omega\cap X+W_2$.
Assume that $\dim Y/(W_1+W_2)<+\infty$. Then the following hold.
\newline (a) We have $\omega$-closed linear subspaces $W_1$, $W_2$, $X+W_1$, $X+W_2$, $W_1+W_2$, $W_1^\omega\cap X+W_2^\omega\cap X$, $\dom(\alpha)=W_1^\omega\cap X $, $\dom(\beta)=W_2^\omega\cap X$, $\alpha+Y=W_1^\omega$, $\beta+Y=W_2^\omega$. Moreover, we have $(\gamma,\delta)\in \Ff\Ll_0(Z)$, and
\begin{equation}\label{e:Fredholm-X-Y}
\begin{cases}
\dim Y/(W_1+W_2)=\dim(W_1^\omega\cap W_2^\omega\cap X),\\
\dim (W_1\cap W_2)=\dim X/(W_1^\omega\cap X+W_2^\omega\cap X).
\end{cases}
\end{equation}
\newline (b) We have $W_1^{\w}\cap X,W_2^{\w}\cap X,W_1^{\w}\cap W_2^{\w}\cap X,W_1^\omega\cap X+W_2^\omega\cap X\in\Ss^c(X)$,
and $W_1,W_2,W_1\cap W_2,W_1+W_2\in\Ss^c(Y)$.
\newline (c) Set $W\ :=(W_1+W_2)^{\w}$. Let $U$, $V$ be two isotropic subspace of $Z$. Assume that $\dom (U)=W_1^\omega\cap X$, $ \dom (V)=W_2^\omega\cap X$. Then we have $W_1^{\w}=U+Y$, $W_2^{\w}=V+Y$, $(U+V)\cap W=U\cap W+V\cap W$, $U\cap Y\subset W_1+W_2$, $V\cap Y\subset W_1+W_2$, $\dim R_W(Z)=2\dim Y/(W_1+W_2)$, $R_W(U)$, $R_W(V)$ are Lagrangian subspaces of $R_W(Z)$, $U+V+W=W_1^\omega\cap X+W_2^\omega\cap X+Y$, and 
\begin{align}\label{e:dim-z-uvw}
	\dim Z/(U+V+W)=&\dim(W_1\cap W_2)\ge\dim(U\cap V\cap Y).
\end{align}. 

If we assume in addition that $U\cap Y+V\cap Y\supset W_1+W_2$ and $\dim(U\cap V\cap Y)=\dim(W_1\cap W_2)$, we have $U\cap Y+V\cap Y=W_1+W_2$ and
$(U,V)\in\Ff\Ll_0(Z)$.
\end{theorem}

\begin{proof} (a) Since $\dom(\alpha)=(\alpha+Y)\cap X$ and $\dom(\beta)=(\beta+Y)\cap X$, by Lemma \ref{l:space-distribution}, we have
\begin{align}\label{e:dom-alpha-beta}
 \dom(\alpha)+\dom(\beta)=(\alpha+\beta+Y)\cap X.
\end{align}
Since $(\alpha,\beta)\in\Ff\Ll_0(Z)$, by \cite[Lemma 1.2.8]{BoZh18}, $\alpha+\beta$ is $\omega$-closed and there holds
\[\dim(\alpha+\beta+Y)/(\alpha+\beta)\le\dim Z/(\alpha+\beta)<+\infty .\] 
By Lemma \ref{l:finite-extension-closed} and \eqref{e:dom-alpha-beta}, $\alpha+\beta+V$ and  $\dom(\alpha)+\dom(\beta)$ are $\omega$-closed. 
By Lemma \ref{l:finite-extension-closed}, we have
\begin{align}
\dim X/(\dom(\alpha)+\dom(\beta))&=\dim(\alpha+\beta+Y+X)/(\alpha+\beta+Y)\notag\\
&=\dim Z/(\alpha+\beta+Y)=\dim(\alpha+\beta+Y)^\omega\notag\\
&=\dim (\alpha\cap \beta\cap Y)=\dim (W_1\cap W_2) \label{eq:1}.
\end{align}

Since $W_1,W_2\subset Y$, they are both isotropic subspaces of $Z$.
By Lemma \ref{l:intersection-omega}, $W_1$ and $W_2$ are $\omega$-closed. By Lemma \ref{l:closed-sum2}, $W_1+X$ and $W_2+X$ are $\omega$-closed.
We also have
\[\dom (\alpha)\subset((\alpha+Y)\cap X)^{\omega\omega}\subset(X+Y\cap\alpha)^\omega=W_1^\omega\cap X.\]
Similarly we have $\dom(\beta)\subset W_2^{\omega}\cap X$.

Note that we have $\dim Y/(W_1+W_2)<+\infty$.
Since $(\alpha,\beta)\in\Ff\Ll_0(Z)$, by Lemma \ref{l:closed-sum2}, $W_1+W_2$ is $\omega$-closed. Since $Z=X\oplus Y$,
by Lemma \ref{l:closed-sum2} again, $X+W_1+W_2$ is $\omega$-closed.
By \cite[Lemma 1.1.2]{BoZh18}, we have
\begin{align*}
\dim (W_1^\omega\cap W_2^\omega\cap X )= \dim Z/(X+W_1+W_2)=\dim Y/(W_1+W_2).
\end{align*}
Note that $X=X^{\omega}=X^{\omega\omega}$. By \cite[Lemma 1.1.2]{BoZh18} we have
\begin{align*}
 \dim X/(W_1^\omega\cap X+W_2^\omega\cap X)&\geq \dim ((X+W_1)\cap(X+W_2))/X\\
 &=(X+W_1\cap W_2)/X=\dim(W_1\cap W_2)\\
 &=\dim X/(\dom(\alpha)+\dom(\beta))\\
 &\geq \dim X/(W_1^\omega\cap X+W_2^\omega\cap X).
\end{align*}
It follows that
\begin{equation}\label{e:dom-index}
\dim X/(W_1^\omega\cap X+W_2^\omega\cap X)=\dim (W_1\cap W_2).
\end{equation}
Then we have
\begin{align*}
 \Index(\gamma,\delta)=&\dim(W_1\cap W_2)+\dim(W_1^\omega\cap W_2^\omega\cap X)\\
 &-\dim Y/(W_1+W_2)-\dim X/(W_1^\omega\cap X+W_2^\omega\cap X)\\
 = &0.
\end{align*}
Since $X$ and $Y$ are Lagrangian subspaces of $Z$, $\gamma$ and $\delta$ are both isotropic subspaces of $Z$.
By \cite[Proposition 1]{BoZh13}, $\gamma$ and $\delta$ are Lagrangian subspaces of $Z$ and $\gamma+\delta$ is $\omega$-closed.

Since $X+W_1+W_2=X+(\gamma+\delta)\cap Y$, and
$W_1+W_2$ are $\omega$-closed, by Lemma \cite[Lemma 1.1.2]{BoZh18} and \eqref{eq:1}, we have
\begin{align*}
 \dim Y/(W_1+W_2)&=Z/(X+W_1+W_2)\\
 &\geq\dim(X\cap W_1^{\omega}\cap W_2^{\omega})\\
 &\ge\dim(\dom(\alpha)\cap\dom(\beta))\\
 &=\dim (\alpha+Y)\cap X\cap(\beta+Y)\\
 &=\dim Z/((\alpha+Y)\cap X\cap(\beta+Y))^\omega\\
 &= \dim Z/(X+W_1+W_2)\\
&=\dim Y/(W_1+W_2).
 \end{align*}

Then we conclude that
\begin{equation*}
 \begin{cases}
 W_1^\omega\cap X\supset \dom(\alpha),\quad W_2^\omega\cap X\supset \dom(\beta),\\
W_1^\omega\cap X +W_2^\omega\cap X= \dom(\alpha)+\dom(\beta),\\
 W_1^\omega\cap X \cap W_2^\omega\cap X=\dom(\alpha)\cap\dom(\beta).
 \end{cases}
\end{equation*}
By Lemma \ref{l:compare-subspaces}, we have
\begin{equation}\label{e:domain-ab}
 \begin{cases}
W_1^\omega\cap X=\dom (\alpha),\\
W_2^\omega\cap X=\dom (\beta).
 \end{cases}
\end{equation}
By Lemma \ref{l:intersection-omega}, $W_1^\omega\cap X$ and $W_2^\omega\cap X$ is $\omega$-closed.
By \cite[Lemma A.1.1]{BoZh18} we have
\begin{align*}
	\alpha+Y&=(\alpha+Y)\cap X+Y=\dom(\alpha)+Y,\\
	W_1^{\w}&=W_1^{\w}\cap X+Y.
\end{align*}
So we have $\alpha+Y=W_1^{\w}$. Similarly we have
$\beta+Y=W_2^{\w}$. By Lemma \ref{l:triple-omega}, $W_1^{\w}$ and $W_2^{\w}$ are $\omega$-closed.
\newline (b) By (\ref{e:Fredholm-X-Y}) and \cite[Lemma A.2.6]{BoZh18}.
\newline (c) 1. By (a), we have $W^{\w}=W_1+W_2$.
Since $W_1\subset Y\subset W_1^{\w}$ and $W_2\subset Y\subset W_2^{\w}$, by \cite[Lemma A.1.1]{BoZh18}, we have
\begin{align*}
W_1^{\w}&=W_1^{\w}\cap(X+Y)=W_1^{\w}\cap X+Y\\
&=\dom(U)+Y=(U+Y)\cap X+Y\\
&=(U+Y)\cap(X+Y)=U+Y.
\end{align*}
Similarly we have $W_2^{\w}=V+Y$. So we have $W=W_1^{\w}\cap W_2^{\w}=(U+Y)\cap(V+Y)$.
By Lemma \ref{l:separate-space}, we have $(U+V)\cap W=U\cap W+V\cap W$.

Since $Y\subset W$, by \cite[Lemma A.1.1]{BoZh18}, we have $W=W\cap(X+Y)=W\cap X+Y$. Since $W$ and $W_1+W_2$ are $\omega$-closed and $W\supset Y$, by Lemma \ref{l:V-omega-la} and \cite[Lemma 1.1.2]{BoZh18}, we have
\begin{align*}
 \dim R_W(Z)&=\dim W/W^\omega =\dim W/Y +\dim Y/W^{\omega}\\
 &=2\dim W/Y=2\dim Y/(W_1+W_2)\\
 &=2\dim (W\cap X).
\end{align*}
Since $W^\omega\subset Y\subset W$ and $W\cap X\subset U+Y$, by Lemma \ref{l:algebra}.b, we have
\begin{align}
\dim R_W(U)&=\dim (U\cap W+W^{\w})/W^{\w}\notag\\
&=\dim (U\cap W)/(U\cap W^{\w})\notag\\
&=\dim (U\cap W\cap(X+Y))/(U\cap W^{\w})\notag\\
&\ge\dim (U\cap(W\cap X+Y))/(U\cap Y)\notag\\
&=\dim ((W\cap X)\cap(U+Y))/(W\cap X\cap Y)\notag\\
&=\dim (W\cap X)\label{e:index-reduction}.
\end{align}
Since $U$ is isotropic, $R_W(U)$ is an isotropic subspace of $R_W(Z)$. So we have $\dim R_W(U)\le\dim(W\cap X)$.
By (\ref{e:index-reduction}), we have $\dim R_W(U)=\dim(W\cap X)$ and $R_W(U)\in\Ll(R_W(Z))$. The identity in (\ref{e:index-reduction})
shows that 
\[U\cap Y=U\cap W^{\w}\subset W^{\w}=W_1+W_2.\] Similarly we have $R_W(V)\in\Ll(R_W(Z))$ and $V\cap Y=V\cap W^{\w}\subset W_1+W_2$.
So we get $\Index(R_W(U),R_W(V))=0$.

Note that by Lemma \ref{l:space-distribution}, we have
\begin{align*}
U+V+W=&U+V+(U+Y)\cap(V+Y)\\
=&U+V+Y\\
=&\dom (U)+\dom (V)+Y\\
=&W_1^\omega\cap X+W_2^\omega\cap X+Y.
\end{align*}
By (a) and Lemma \ref{l:closed-sum}, $U+V+W$ is $\omega$-closed.
It follows that
\begin{align*}
\dim Z/(U+V+W)&=\dim X/(W_1^\omega\cap X+W_2^\omega\cap X)
\end{align*}
By \eqref{e:dom-index} and \cite[Lemma 1.1.2.b]{BoZh18}, we have
\begin{align*}
	\dim(W_1\cap W_2)=&\dim Z/(U+V+W)\\
	\ge&\dim(U^\omega\cap V^\omega\cap W^{\w})\ge\dim(U\cap V\cap W^{\w})\\
	=&\dim(U\cap V\cap Y).
\end{align*}

2. Now we assume that $U\cap Y+V\cap Y \supset W_1+W_2=W^{\omega}$ and $\dim(U\cap V\cap W^{\omega})=\dim(U\cap V\cap Y)=\dim(W_1\cap W_2)$. Then we have 
\begin{align*}
	W^\omega\supset U\cap W^\omega+V\cap W^\omega=U\cap Y+V\cap Y\supset W^\omega.
\end{align*}
Thus we have $W^\omega=U\cap W^\omega+V\cap W^\omega$.
Since $U\subset W_1^{\w}$ and $V\subset W_2^{\w}$, we have $W\supset U\cap V$. Then we have 
\begin{align*}
	W^\omega&\supset(U+V\cap W)\cap W^\omega\\
	&\supset (U\cap W+V\cap W)\cap W^\omega\\
	&=(U+V)\cap W\cap W^\omega=W^\omega.	
\end{align*}
Thus we obtain $(U+V\cap W)\cap W^\omega=W^\omega$.
By Lemma \ref{l:g-index-reduction}, we have
\begin{align*}
    \Index(U,V)=&\Index(R_W(U),R_W(V)) \\
    &+\dim(U\cap V\cap W^{\w})-\dim Z/(W+U+V)\\
    =&0.
\end{align*}
Since $U,V$ are both isotropic subspaces, by \cite[Proposition 1]{BoZh13}, we have $(U,V)\in\Ff\Ll_0(X)$.
\end{proof}

\section{The Maslov index for the path of pairs of Lagrangian subspaces in triangular form}\label{s:maslov-triangle}
We fix our data and choices as following.
\begin{data}
Let $Z$ be a Banach space with continuously varying symplectic structures $\w(s)$, $s\in[0,1]$.
Assume that $Z=X(s)\oplus Y(s)$ with two paths $\{X(s)\in\Ll(Z,\w(s));\;s\in[0,1]\}$ and $\{Y(s)\in\Ll(Z,\w(s));\;s\in[0,1]\}$.
Let $\{(\lambda(s),\mu(s))\in\Ff\Ll_0(Z,\w(s)); \;s\in [0,1]\}$ be a path. Set $A(s)\;:=\lambda(s)\cap Y(s)$, $B(s)\;:=\mu(s)\cap Y(s)$,
$\gamma(s)\;:=A(s)^{\omega(s)}\cap X(s)\oplus A(s)$, $\delta(s)\;:=B(s)^{\omega(s)}\cap X(s)\oplus B(s)$. Assume that
\begin{equation}\label{e:data}
 \begin{cases}
\lambda(s)\cap Y(s),\mu(s)\cap Y(s) \text{ are two paths, and}\\
\dim Y(s)/(A(s)+B(s))<+\infty.
 \end{cases}.
\end{equation}
\end{data}

By \cite[Lemma I.4.10]{Ka95}, we can fix $X(s)$ and $Y(s)$ locally in the considerations. Note that each $\omega(s)$-closed linear subspace is a closed subspace of $Z$.
By Theorem \ref{t:triangle-pairs}.b, we have
\begin{equation}\label{e:data1}
\begin{cases}
\lambda(s)+Y(s)=A(s)^{\omega(s)}\cap X(s)+Y(s)\in \Ss^c(Z),\\
\mu(s)+Y(s)=B(s)^{\omega(s)}\cap X(s)+Y(s)\in \Ss^c(Z).
\end{cases}.
\end{equation}
Since $\{\lambda(s)\cap Y(s);\;s\in[0,1]\}$ and $\{\mu(s)\cap Y(s);\;s\in[0,1]\}$ are both paths, by \cite[Proposition A.3.13]{BoZh18},
$\{\mu(s)+Y(s);\;s\in[0,1]\}$ and $\{\mu(s)+Y(s);\;s\in[0,1]\}$ are both paths in $\Ss(Z)$.
Since $\lambda(s)+Y(s)+X(s)=\mu(s)+Y(s)+X(s)=Z$, by \cite[Proposition A.3.13]{BoZh18},
$\{(\lambda(s)+Y(s))\cap X(s);\;s\in[0,1]\}$, $\{(\mu(s)+Y(s))\cap X(s);\;s\in[0,1]\}$ are both paths in $\Ss(Z)$.

By Theorem \ref{t:triangle-pairs}, we have $\lambda(s)+Y(s)=A(s)^\{\omega(s)\}$, $\mu(s)+Y(s)=B(s)^\{\omega(s)\}$, and $(\gamma(s),\delta(s))\in \Ff\Ll_0(Z,\w(s))$.

We need some preparations for the calculation of the Maslov index of the path $(\lambda(s),\mu(s))$, $s\in[0,1]$.

\begin{lemma}\label{l:compare-relation} Let $Z$ be a vector space with linear subspaces $X$ and $Y$ such that $Z=X\oplus Y$. A linear subspace of $Z$ is viewed as linear relation from $X$ to $Y$. Let $\alpha$, $\beta$, $M$ be linear subspaces in $Z$. Then the following hold.
\newline (a) Assume that $\alpha\supset\beta$, $\dom(\alpha)\subset\dom(\beta)$ and $\alpha\cap Y\subset\beta\cap Y$. Then we have $\alpha=\beta$.
\newline (b) Assume that $Y=\alpha\cap Y+M\cap Y$. Then we have $\alpha=\alpha\cap(M+X)+\alpha\cap Y$.
\end{lemma}

\begin{proof} (a) Since $\alpha\supset\beta$, we have $\dom(\alpha)\supset\dom(\beta)$ and $\alpha\cap Y\supset\beta\cap Y$. So there hold $\dom(\alpha)=\dom(\beta)$ and $\alpha\cap Y=\beta\cap Y$. For each $x\in\dom(\alpha)=\dom(\beta)$, take a vector $y\in\beta x$. Since $\alpha\supset\beta$, we have $y\in\alpha x$. Then there holds $\alpha x=y+\alpha\cap Y=y+\beta\cap Y=\beta x$. Hence we have $\alpha=\beta$.
\newline (b) Set $\beta:=\alpha\cap(M+X)+\alpha\cap Y$. Then we have $\alpha\supset\beta$ and $\alpha\cap Y\subset\beta\cap Y$. By \cite[Lemma A.1.1]{BoZh18}, we have 
\begin{align*}
	M+X=(M+X)\cap(X+Y)=X+(M+X)\cap Y\supset X+M\cap Y.
\end{align*}
Since $Y=\alpha\cap Y+M\cap Y$, by \cite[Lemma A.1.1]{BoZh18}, we have
\begin{align*}
\dom(\beta)&=\dom (\alpha\cap(M+X)+\alpha\cap Y)\\
&=(\alpha\cap(M+X)+\alpha\cap Y+Y)\cap X \\
&\supset(\alpha\cap (X+M\cap Y)+\alpha\cap Y+Y)\cap X\\
&=(\alpha\cap (X+M\cap Y+\alpha\cap Y)+Y)\cap X\\
&=(\alpha\cap (X+Y)+Y)\cap X\\
&=(\alpha+Y)\cap X\\
&=\dom(\alpha).
\end{align*}
By (a), we have $\alpha=\beta$.
\end{proof}

\begin{lemma}\label{l:two-finite} Let $(X,\omega)$ be a symplectic vector space with two finite dimensional linear subspaces $V$, $W$ such that $\dim V=\dim W$. Then the following four conditions are equivalent:
\begin{itemize}
\item[(i)] $V\cap W^{\w}=\{0\}$,
\item[(ii)] $V+W^{\w}=X$,
\item[(iii)] $V^{\w}\cap W=\{0\}$,
\item[(iv)] $V^{\w}+W=X$.
\end{itemize}
\end{lemma}

\begin{proof} By \cite[Lemma 1.1.2.b]{BoZh18} we have $\dim X/V^{\w}=\dim V=\dim W=\dim X/W^{\w}$. So we have (i)$\Leftrightarrow$(ii)$\Rightarrow$(iii)$\Leftrightarrow$(iv)$\Rightarrow$(i).
\end{proof}

The following simple lemma establish some vector space is symplectic.

\begin{lemma}\label{l:iso-to-sym}
	Let $Z$ be a vector space with a skew symmetric form $\omega$ and two linear subspaces $X$, $Y$. Assume that $X\subset X^{\omega}$, $Y\subset Y^{\omega}$, $X^{\omega}\cap Y=Y^{\omega}\cap X=\{0\}$, and $Z=X+Y$. Then the space $(Z,\omega)$ is symplectic, $X$ and $Y$ are Lagrangian subspaces of $Z$, and we have $Z=X\oplus Y$.
\end{lemma}

\begin{proof}
	By \cite[Lemma A.1.1]{BoZh18}, we have 
	\begin{align*}
		X^{\omega}=X^{\omega}\cap(X+Y)=X+X^{\omega}\cap Y=X.
	\end{align*}
    Similarly we have $Y^{\omega}=Y$. Then we have
    	\begin{align*}
    	Z^{\omega}=X^{\omega}\cap Y^{\omega}=X\cap Y=X^{\omega}\cap Y=\{0\}.
    \end{align*}
    Thus the space $(Z,\omega)$ is symplectic, and we have $Z=X\oplus Y$.
    By \cite[Proposition 1]{BoZh13}, $X$ and $Y$ are Lagrangian subspaces of $Z$.
\end{proof}

Given a finite dimensional isotropic subspace of a sympectic Banach space in diagonal form, we can always parametrize them locally in a continuous family of sympactic Banach spaces.  
 
\begin{lemma}\label{l:family-isotropic-separate}
	Let $B$ be a topological space. Let $Z$ be a Banach space with continuously varying symplectic structures $\{\w(b);\; b\in B\}$.
    Assume that $Z=X(b)\oplus Y(b)$ with two continuous families $\{X(b)\in \Ll(Z,\w(b));\; b\in B\}$ and $\{Y(b)\in \Ll(Z,\w(b));\; b\in B\}$. Fix a $b_0\in B$. Let $V(b_0)$ be a finite dimensional isotropic subspace of $(Z,\w(b_0))$ with $V(b_0)=V(b_0)\cap X(b_0)+V(b_0)\cap Y(b_0)$. Then there exist an neighborhood $U$ of $b_0$ and a continuous family $\{V(b)\}_{b\in U}$ such that, for each $b\in U$, $V(b)$ is an isotropic subspace of $(Z,\w(b))$, and $V(b)=V(b)\cap X(b)\oplus V(b)\cap Y(b)$.
\end{lemma}

\begin{proof} 
	By \cite[Lemma I.4.10]{Ka95}, for each $b\in B$ close to $b_0$, there is a Linear isomorphism $L(b)\in\Bb(Z)$ such that, $L(b_0)=I$, the map $b\mapsto L(b)$ is continuous, $L(b)X(b_0)=X(b)$, and $L(b)Y(b_0)=Y(b)$. Then we have $X(b_0),\;Y(b_0)\in \Ll(Z,L(b)^*\w(b))$. Then we can assume that $X(b)=X$ and $Y(b)=Y$. 
	
	Since $(V(b_0)\cap X)^{\w(b_0)}\supset X$, by \cite[Lemma A.1.1]{BoZh18}, we have $(V(b_0)\cap X)^{\w(b_0)}=X+(V(b_0)\cap X)^{\w(b_0)}\cap Y$. Similarly we have $(V(b_0)\cap Y)^{\w(b_0)}=(V(b_0)\cap Y)^{\w(b_0)}\cap X+Y$.

    Since $V(b_0)$ is finite dimensional and $V(b_0)\cap X\cap Y=\{0\}$, by Lemma \ref{l:V-omega-la}, we have $\dim X/((V(b_0)\cap Y)^{\w(b_0)}\cap X)=\dim(V(b_0)\cap Y)$. Then there is a linear subspace $W_1$ of $X$ such that $X=W_1\oplus (V(b_0)\cap Y)^{\w(b_0)}\cap X$ and $\dim W_1=\dim(V(b_0)\cap Y)$. Similarly there is a linear subspace $W_2$ of $Y$ such that $Y=W_2\oplus (V(b_0)\cap X)^{\w(b_0)}\cap Y$ and $\dim W_2=\dim(V(b_0)\cap X)$.
    So we have
    \[W_1\cap(V(b_0)\cap Y)^{\w(b_0)}=W_1\cap(V(b_0)\cap Y)^{\w(b_0)}\cap X=\{0\}.\]
    Similarly we have $W_2\cap(V(b_0)\cap X)^{\w(b_0)}=\{0\}$. By Lemma \ref{l:two-finite}, we have $W_1^{\w(b_0)}\cap V(b_0)\cap Y=W_2^{\w(b_0)}\cap V(b_0)\cap X=\{0\}$.

    Since $V(b_0)$ is an isotropic subspace of $(Z,\w(b_0))$, we have $V(b_0)\cap X\subset (V(b_0)\cap Y)^{\w(b_0)}\cap X$ and $V(b_0)\cap Y\subset (V(b_0)\cap X)^{\w(b_0)}\cap Y$. So we have $V(b_0)\cap X\cap W_1=V(b_0)\cap Y\cap W_2=\{0\}$.

    Set $\tilde X\;:=V(b_0)\cap X\oplus W_1$, $\tilde Y\;:=V(b_0)\cap Y\oplus W_2$, and $\tilde Z\;:=\tilde X\oplus\tilde Y$. Then we have
    \begin{align*}
	    \dim\tilde Z=2\dim (V(b_0)\cap X)+2\dim (V(b_0)\cap Y)=2\dim V(b_0).
    \end{align*} 

    Since $\tilde X\subset X$ and $\tilde Y\subset Y$ are isotropic subspaces of $(Z,\omega(b_0))$, $\tilde Y^{\w(b_0)}\cap \tilde X\subset(V(b_0)\cap Y)^{\w(b_0)}\cap W_1=\{0\}$ and 
    $\tilde X^{\w(b_0)}\cap \tilde Y\subset(V(b_0)\cap X)^{\w(b_0)}\cap W_2=\{0\}$, by Lemma \ref{l:iso-to-sym}, 
    the spacce $(\tilde Z,\w(b_0)|_{\tilde Z})$ is symplectic.
    So we have $V(b_0)\in\Ll(\tilde Z,\w(b_0)|_{\tilde Z})$.

    Since $\dim\tilde Z<+\infty$, there is a neighborhood $U$ of $b_0$ such that $(\tilde Z,\w(b)|_{\tilde Z})$ is symplectic for each $b\in U$. By taking a suitable basis in $\tilde X$ and $\tilde Y$, we can assume that
    \[(\omega(b)|_{\tilde Z})(u,v)=\lla f(b)(x_1),y_2\rra-\overline{\lla f(b)(x_2),y_1\rra}\]
    for each $u=x_1+y_1$, $v=x_2+y_2$, $x_1,x_2\in \tilde X$, $y_1,y_2\in \tilde Y$, $b\in U$ under a continuous family of linear isomorphisms $f(b):\tilde X\to \tilde Y$, $b\in U$, where $\lla\cdot,\cdot\rra$ denotes the standard inner product. Since $V(b_0)\in\Ll(\tilde Z,\w(b_0)|_{\tilde Z})$, we have $f(b_0)(V(b_0)\cap X)=(V(b_0)\cap Y)^{\bot_2}$, where $\bot_2$ denotes the orthogonal complement in $\tilde Y$. Define
    \[V(b):=V(b_0)\cap X\oplus (f(b)(V(b_0)\cap X))^{\bot_2}.\]
    Then the family $\{V(b)\}_{b\in U}$ satisfies the desired properties.
\end{proof}

The following proposition make our path in triangular form homotopic rel. endpoints to the composition of a diagonal path and two path with fixed diagonal part.

\begin{proposition}\label{p:continuity-data} For the Data given above, the following hold.
\newline (a) For each $s_0\in[0,1]$, there exists a $d>0$ such that there are two paths $E(s),F(s)\in\Ll(Z,\w(s))$, $s\in (s_0-d,s_0+d)\cap[0,1]$ which satisfy
\begin{equation}\label{e:continuity-data}
\begin{cases}
E(s)\oplus\gamma(s)=F(s)\oplus\delta(s)= Z,\\
E(s)=E(s)\cap X(s)\oplus E(s)\cap Y(s),\\ F(s)=F(s)\cap X(s)\oplus F(s)\cap Y(s),\\
\dim\frac{\gamma(s)}{\gamma(s)\cap F(s)}=\dim\frac{\delta(s)}{\delta(s)\cap E(s)}=\dim \frac{Z}{\gamma(s_0)+\delta(s_0)}.
\end{cases}
\end{equation}
\newline (b) For each $s,t\in [0,1]$, we define the operator $P(s,t):Z\rightarrow Z$ by
\[P(s,t)(x+y)\defeq x+ty,\quad x\in X(s),y\in Y(s).\]
For each $s\in (s_0-d,s_0+d)\cap[0,1]$ and $t\in[0,1]$, we define the subspaces of $Z$ by
\begin{align*}M(s,t)\;:&=P(s,t)((E(s)\cap X(s))^{\omega(s)}\cap \lambda(s))+A(s),\\
N(s,t)\;:&=P(s,t)((F(s)\cap X(s))^{\omega(s)}\cap \mu(s))+B(s).
\end{align*}
Then we have
\begin{itemize}
\item[(i)] The two families $(E(s)\cap X(s))^{\omega(s)}\cap \lambda(s)$ and $(F(s)\cap X(s))^{\omega(s)}\cap\mu(s)$, $s\in(s_0-d,s_0+d)\cap[0,1]$ are both paths in $\Ss(Z)$,
\item[(ii)] $\lambda(s)=M(s,1)=(E(s)\cap X(s))^{\omega(s)}\cap \lambda(s)\oplus A(s)$, $\mu(s)=N(s,1)=(F(s)\cap X(s))^{\omega(s)}\cap \mu(s)\oplus B(s)$,
\item[(iii)] $\dom (M(s,t))=\dom (\lambda(s))$, $
M(s,t)\cap Y=A(s)$, $\dom (N(s,t))=\dom (\mu(s))$, $N(s,t)\cap Y=B(s)$, $\gamma(s)=M(s,0)$, $\delta(s)=N(s,0)$, and
\item[(iv)] the maps $M(\cdot,\cdot)$ and $N(\cdot,\cdot)$ are both continuous maps from $((s_0-d,s_0+d)\cap[0,1])\times [0,1]$ to $\Ss(Z)$,
and for each $(s,t)\in((s_0-d,s_0+d)\cap[0,1])\times [0,1]$, we have $(M(s,t), N(s,t))\in\Ff\Ll_0(Z,\w(s))$.
\end{itemize}
\end{proposition}

\begin{proof} (a) By \cite[Lemma I.4.10]{Ka95} and the first paragraph of the proof of Lemma \ref{l:family-isotropic-separate}, we can assume that $X(s)=X$ and $Y(s)=Y$.
	
By Lemma \ref{l:diag-decomposition}, there is a finite dimensional isotropic subspace $V(s_0)$ of $(Z,\w(s_0))$ such that $Z=V(s_0)\oplus(\gamma(s_0)+\delta(s_0))$ and $V(s_0)=V(s_0)\cap X\oplus V(s_0)\cap Y$. By Lemma \ref{l:family-isotropic-separate}, there exist a $d_1>0$ and a continuous family $V(s)$, $s\in(s_0-d_1,s_0+d_1)\cap[0,1]$ of isotropic subspace of $(Z,\w(s))$ such that $V(s)=V(s)\cap X\oplus V(s)\cap Y$. Set $E(s)\ :=V(s)+ V(s)^{\omega(s)}\cap \delta(s)$ and $F(s) :=V(s)+V(s)^{\w(s)}\cap\gamma(s)$. By Lemma \ref{l:condition}, we have $Z=E(s_0)\oplus\gamma(s_0)=F(s_0)\oplus\delta(s_0)$. By \cite[Proposition A.3.5]{BoZh18}, there exists a $d_2\in(0,d_1)$ such that for each $s\in(s_0-d_2,s_0+d_2)\cap[0,1]$, we have $V(s)\cap\gamma(s)=V(s)\cap\delta(s)=\{0\}$. By \cite[Lemma 3.1.1]{BoZh18}, the family $\{V(s)^{\omega(s)};\;s\in(s_0-d_1,s_0+d_1)\cap[0,1]\}$ is continuous. By Lemma \ref{l:V-omega-la}, we have 
\begin{align*}
	V(s)^{\omega(s)}+\gamma(s)=V(s)^{\omega(s)}+\delta(s)=Z
\end{align*}
for each $s\in(s_0-d_2,s_0+d_2)\cap[0,1]$.
By \cite[Proposition A.3.13]{BoZh18}, the families $\{V(s)^{\omega(s)}\cap\gamma(s);\;s\in(s_0-d_2,s_0+d_2)\cap[0,1]\}$ and $\{V(s)^{\omega(s)}\cap\delta(s);\;s\in(s_0-d_2,s_0+d_2)\cap[0,1]\}$ are continuous. By \cite[Proposition A.3.13]{BoZh18} again, the families $\{E(s);\;s\in(s_0-d_2,s_0+d_2)\cap[0,1]\}$ and $\{F(s);\;s\in(s_0-d_2,s_0+d_2)\cap[0,1]\}$ are continuous.
Then there exists a $d\in(0,d_2)$ such that for all $s\in(s_0-d,s_0+d)\cap[0,1]$, we have $Z=E(s)\oplus\gamma(s)=F(s)\oplus\delta(s)$.
By Corollary \ref{c:diag-decomposition}, $E(s),F(s)\in\Ll(Z,\w(s))$ and (\ref{e:continuity-data}) hold.
\newline (b) (i) By \eqref{e:continuity-data}, we have $Y=E(s)\cap Y\oplus \gamma(s)\cap Y=E(s)\cap Y\oplus A(s)$. Since $E(s)\in \Ll(Z,\w(s))$, we have
\begin{align*}
 (E(s)\cap X)^{\omega(s)}+\lambda(s)&\supset E(s)+X+\lambda(s)\\
 &\supset E(s)\cap Y+X+A(s)\\
 &=X+Y=Z.
\end{align*}
Then by \cite[Proposition A.3.13]{BoZh18}, the family $(E(s)\cap X(s))^{\omega(s)}\cap \lambda(s)$, $s\in(s_0-d,s_0+d)\cap[0,1]$ is a path in $\Ss(Z)$.
Similarly, the family $(F(s)\cap X(s))^{\omega(s)}\cap\mu(s)$, $s\in(s_0-d,s_0+d)\cap[0,1]$ is a path in $\Ss(Z)$.

(ii) Since $E(s)=E(s)\cap X+E(s)\cap Y$, we have $E(s)+X=X+E(s)\cap Y$. By Lemma \ref{l:diag-omega}, we have $(E(s)\cap X)^{\omega(s)}= E(s)+X=X+E(s)\cap Y$. By the proof of (i), we have $Y=E(s)\cap Y\oplus A(s)=E(s)\cap Y\oplus \lambda(s)\cap Y$
Then we have
\begin{align*}
(E(s)\cap X)^{\omega(s)}\cap \lambda(s)\cap A(s)&=(X+E(s)\cap Y)\cap A(s)\\
&=E(s)\cap Y\cap A(s)=\{0\}.
\end{align*}
By Lemma \ref{l:compare-relation}.b, we have $(E(s)+X)\cap\lambda(s)+A(s)=\lambda(s)$. Then we have $\lambda(s)=((E(s)\cap X(s))^{\omega(s)}\cap \lambda(s)\oplus A(s)=M(s,1)$. Similarly we have $\mu(s)=(F(s)\cap X(s))^{\omega(s)}\cap \mu(s)\oplus B(s)=N(s,1)$.

(iii) By (ii), we have
\begin{align*}
\dom (M(s,t))&=\dom(P(s,t)((E(s)\cap X)^{\omega(s)}\cap \lambda(s)))\\
&=\dom((E(s)\cap X)^{\omega(s)}\cap \lambda(s))=\dom (\lambda(s)),\text{ and}\\
M(s,t)\cap Y &=A(s)+Y\cap P(s,t)((E(s)\cap X)^{\omega(s)}\cap \lambda(s)) \\
&=\begin{cases}A(s)+Y\cap (E(s)\cap X)^{\omega(s)}\cap \lambda(s)& \text{if }t\in(0,1],\\
A(s)&\text{if }t=0
\end{cases}\\
&=A(s).
\end{align*}
Similarly we have $\dom (N(s,t))=\dom (\mu(s))$, and $N(s,t)\cap Y=B(s)$. By Theorem \ref{t:triangle-pairs}, we have $\dom(\lambda(s))=A(s)^{\omega(s)}\cap X$ and
$\dom(\mu(s))=B(s)^{\omega(s)}\cap X$. Then we have $\gamma(s)=M(s,0)$ and $\delta(s)=N(s,0)$.

(iv) We claim that $M(s,t)$ and $N(s,t)$ are isotropic subspaces of $(Z,\w(s))$. $(s,t)\in((s_0-d,s_0+d)\cap[0,1])\times [0,1]$.
Let $x_1+ty_1, x_2+ty_2\in P(s,t)((E(s)\cap X)^{\omega(s)}\cap \lambda(s))$, $u_1,u_2\in A(s)$ be four vectors, where $x_1,x_2\in X$ and $y_1,y_2\in Y$. Then we have $x_1+y_1,x_2+y_2\in \lambda(s)$.
Note that $\omega(s)(u_1,y_2)=\omega(s)(y_1,u_2)=0$ and $\omega(s)(u_1,x_2+y_2)=\omega(s)(x_1+y_1,u_2)=0$.
It follows that
\begin{align*}
\omega(s)&(x_1+ty_1+u_1,x_2+ty_2+u_2)=t\omega(s)(y_1,x_2)+t\omega(s)(x_1,y_2)\\
&\quad+\omega(s)(u_1,x_2)+\omega(s)(x_1,u_2)+\omega(s)(u_1,y_2)+\omega(s)(y_1,u_2)\\
&=t\omega(s)(x_1+y_1,x_2+y_2)+\omega(s)(u_1,x_2+y_2)+\omega(s)(x_1+y_1,u_2)\\
&=0.
\end{align*}
So $M(s,t)$ is an isotropic subspace of $(Z,\w(s))$. Similarly, $N(s,t)$ is an isotropic subspace of $(Z,\w(s))$.

By Theorem \ref{t:triangle-pairs}.c, we have we have $(M(s,t), N(s,t))\in\Ff\Ll_0(Z,\w(s))$.
By \cite[Lemma A.2.6]{BoZh18}, we have $M(s,t),N(s,t)\in \Ss(Z)$.

By (ii), we have $\dom((E(s)\cap X)^{\omega(s)}\cap \lambda(s))=\dom(\lambda(s))$. Since $(E(s)\cap X)^{\omega(s)}\cap \lambda(s)\cap Y=\{0\}$, there is an operator $C(s):\dom(\lambda(s))\rightarrow Y$ such that
$\Graph(C(s))=(E(s)\cap X)^{\omega(s)}\cap \lambda(s)$.
Since $(E(s)\cap X)^{\omega(s)}\cap \lambda(s)$ and $\dom(C(s))= \dom(\lambda_s)$ are both closed, $C(s)$ is a bounded linear operator.

By the proof of (ii), we have $(E(s)\cap X)^{\omega(s)}=E(s)+X=X+E(s)\cap Y$. By By the proof of (i), we have
\begin{align*}
	E(s)\cap X)^{\omega(s)}+ \lambda(s)\supset X+E(s)\cap Y+A(s)=Z.
\end{align*}
Since $\{\lambda(s);\;s\in(s_0-d,s_0+d)\cap[0,1]\}$ and $\{E(s)\cap Y;\;s\in(s_0-d,s_0+d)\cap[0,1]\}$ are paths, by \cite[Proposition A.3.13]{BoZh18}, $\{(E(s)\cap X)^{\omega(s)}\cap \lambda(s);\;s\in(s_0-d,s_0+d)\cap[0,1]$ is a path in $\Ss(Z)$.
By \cite[Lemma 0.1]{Ne68}, the map $(s,t)\mapsto P(s,t)((E(s)\cap X)^{\omega(s)}\cap \lambda(s))$ is a continuous family from $((s_0-d,s_0+d)\cap[0,1])\times [0,1]$ to $\Ss(Z)$.
Since $M(s,t)=P(s,t)((E(s)\cap X)^{\omega(s)}\cap \lambda(s))+A(s)$ is closed, by Theorem \ref{t:operator-space}, we conclude that the family $\{M(s,t)\}_{(s,t)\in((s_0-d,s_0+d)\cap[0,1])\times [0,1]}$ is a continuous one in $\Ss(Z)$. By the same argument we conclude that the family $\{N(s,t)\}_{(s,t)\in((s_0-d,s_0+d)\cap[0,1])\times [0,1]}$ is a continuous one in $\Ss(Z)$.  
\end{proof}

We have the following formula of triple index for the Lagrangian subspaces in diagonal form in a finite dimensional symplectic vector space.

\begin{lemma}\label{lm:Q_diagonal} Let $(Z,\omega)$ be a symplectic vector space such that $Z=X\oplus Y$ with $X,Y\in \Ll(Z)$.
Assume that $\dim Z=2n$. Let $\lambda,\mu, V\in \Ll(Z)$.
Assume that $\lambda=\lambda\cap X+\lambda\cap Y$, $\mu=\mu\cap X+\mu\cap Y$, $V=V\cap X+V\cap Y$, and $Z=\lambda\oplus V=\mu\oplus V$. Then the following hold.
\newline (a) We have $\dim(\lambda\cap X)=\dim(\mu\cap X)=n-\dim(V\cap X)$ and $\dim(\lambda\cap Y)=\dim(\mu\cap Y)=n-\dim(V\cap Y)$.
\newline (b) Set $Q\ :=Q(\mu,V;\lambda)$. Then we have
\begin{align}
\label{e:m-pm-Q1}m^{\pm}(Q)&=\frac{1}{2}(n-\dim(\lambda\cap\mu))\\
\label{e:m-pm-Q2}&=\dim(\la\cap X)-\dim(\la\cap\mu\cap X)\\
\label{e:m-pm-Q3}&=\dim(\la\cap Y)-\dim(\la\cap\mu\cap Y).
\end{align}
\end{lemma}

\begin{proof} (a) Since $X=\la\cap X\oplus V\cap X=\mu\cap X\oplus V\cap X$ and $Y=\la\cap Y\oplus V\cap Y=\mu\cap Y\oplus V\cap Y$, (a) follows.
\newline (b) Set $w:=\lambda\cap \mu $. Let $\pi_w\ :=R_{w^{\w}}$ denote the symplectic reduction map. Set $Q_1\ :=Q(\pi_w(\mu),\pi_w(V);\pi_w(\lambda))$.
Note that $\mu\cap(V+\lambda)=\mu$. By \cite[Lemma 3.3]{ZhWuZh18} and \cite[Lemma A.1.1]{BoZh18}, we have
\begin{align*}
	\ker(Q)=\mu\cap(V+\lambda\cap\mu^\omega)=\mu\cap V+\lambda\cap\mu=\lambda\cap\mu.
\end{align*} 
By \cite[(15)]{ZhWuZh18}, we have
$m^{\pm}(Q)=m^{\pm}(Q_1)$
and $\ker(Q_1)=\{0\}$.
Note that $w=w\cap X \oplus w\cap Y$ and $w^{\w}=\lambda+\mu=w^{\w}\cap X \oplus w^{\w}\cap Y$. By Lemma \ref{l:reduction-separate}, we have
\begin{equation}\label{e:three-reduction-separate}
\begin{cases}
\pi_w(\lambda)=\pi_w(\lambda)\cap\pi_w(X) +\pi_w(V)\cap\pi_w(Y),\\
\pi_w(\mu)=\pi_w(\mu)\cap\pi_w(X) +\pi_w(\mu)\cap\pi_w(Y),\\
\pi_w(V)=\pi_w(V)\cap\pi_w(X) +\pi_w(V)\cap\pi_w(Y).
\end{cases}
\end{equation}
By \cite[Lemma A.1.1]{BoZh18}, we have $(V+\la)\cap(\lambda+\mu)=V\cap(\lambda+\mu)+\la$.  
By Lemma \ref{l:reduction-separate}.b, we have
\begin{align*}
	\pi_w(Z)=\pi_w(V+\la)=\pi_w(V)+\pi_w(\la). 
\end{align*}
Since $(X+w)\cap(Y+w)=w\cap X+w\cap Y=w$, by definition, we have $\pi_w(X)\cap\pi_w(Y)=\{0\}$. By \cite[Lemma A.1.1]{BoZh18} and \eqref{e:three-reduction-separate}, we have \begin{align*}
	\pi_w(\mu)\cap \pi_w(X)\subset& \pi_w(X)=
	\pi_w(X)\cap(\pi_w(V)+\pi_w(\la))\\
	=&\pi_w(V)\cap \pi_w(X)+\pi_w(\la)\cap \pi_w(X).
\end{align*} 
Let $z_1,z_2\in \pi_w(\mu)\cap \pi_w(X)$. We have $z_j=-x_j+y_j$ with $x_j \in \pi_w(V)\cap \pi_w(X)$, $y_j\in \pi_w(\mu)\cap \pi_w(X)$, $j=1,2$.
Note that $x_2,y_1\in \pi_w(X)$ and $\pi_w(X)$ is isotropic. Denote by $\tilde \omega$ the induced symplectic structure on $w^\w/w$. Then we have
\[Q_1(z_1,z_2)=\tilde\omega(z_1,x_2)=\tilde\omega(y_1,x_2)=0.\]

By exchanging $X$ and $Y$, we have
\begin{equation}\label{l:zero-Q1}
\begin{cases}
Q_1|_{\pi_w(\mu)\cap\pi_w(X)}=0,\\
Q_1|_{\pi_w(\mu)\cap\pi_w(Y)}=0.
 \end{cases}
\end{equation}
Set $c:=\dim\pi_w(\mu)$, $a:=\dim (\pi_w(\mu)\cap\pi_w(X))$ and $b:=\dim (\pi_w(\mu)\cap\pi_w(Y))$. Then $a+b=c$.
Since $\ker Q_1=\{0\}$, we have $m^{\pm}(Q_1)\ge a$ and $m^{\pm}(Q_1)\ge b$. So we have
\[c=m^+(Q_1)+m^-(Q_1)\ge a+b=c.\]
Hence there holds that $m^{\pm}(Q_1)=a=b=\frac{c}{2}$.

By Lemma \ref{l:g-index-reduction}, we have $\pi_w(\mu)\cap\pi_w(X)\cong(\mu\cap X)/(\la\cap\mu\cap X)$ and $\pi_w(\mu)\cap\pi_w(Y)\cong(\mu\cap Y)/(\la\cap\mu\cap Y)$. Since $\dim(\la\cap X)=\dim(\mu\cap X)$, we have
\begin{align*}
	a=b&=\dim(\mu\cap X)-\dim(\la\cap\mu\cap X)=\dim(\mu\cap Y)-\dim(\la\cap\mu\cap Y)\\
	&=\dim(\la\cap X)-\dim(\la\cap\mu\cap X)=\dim(\la\cap Y)-\dim(\la\cap\mu\cap Y).
\end{align*}

The number $c$ is given by
\begin{align*}
c=\dim\pi_w(\mu)=\dim\mu -\dim(\lambda\cap\mu)=n-\dim(\lambda\cap\mu).
\end{align*}
Therefore (b) holds.
\end{proof}

With the preparations above, we can calculate the Maslov index for a path of pairs of Lagrangian subspaces in diagonal form.

\begin{proposition}\label{p:maslov-diagonal} Let $Z$ be a Banach space with continuously varying symplectic structures $\w(s)$, $s\in[0,1]$.
Assume that $Z=X(s)\oplus Y(s)$ with two paths $\{X(s)\in\Ll(Z,\w(s));\;s\in[0,1]\}$ and $\{Y(s)\in\Ll(Z,\w(s));\;s\in[0,1]\}$.
Let $\{(\lambda(s),\mu(s))\in\Ff\Ll_0(Z,\w(s)); \;s\in [0,1]\}$ be a path.
Assume that $\lambda(s)=\lambda(s)\cap X(s)+\mu(s)\cap X(s)$ and $\mu(s)=\mu(s)\cap X(s)+\mu(s)\cap Y(s)$ hold for each $s\in [0,1]$. Then the following hold.
\newline (a) Denote by $a(s)\ :=\Index(\la(s)\cap X(s),\mu(s)\cap X(s))$ and $b(s)\ :=\Index(\la(s)\cap Y(s),\mu(s)\cap Y(s))$. Then $a(s)$ and $b(s)$ are well-defined constant integers, and we have
\begin{align}\label{e:index-a(s)-b(s)}
	a(s)=-b(s)=\dim(\la(s)\cap\mu(s)\cap X(s))-\dim(\la(s)\cap\mu(s)\cap Y(s)).
\end{align}
\newline (b) The Maslov index of the path $\{(\lambda(s),\mu(s)); \;s\in [0,1]\}$ is given by
\begin{align}
\label{e:maslov-diagonal1}\Mas_{\pm}&\{\la(s),\mu(s)\}=\pm\frac{1}{2}(\dim(\la(0)\cap\mu(0))-\dim(\la(1)\cap\mu(1)))\\
\label{e:maslov-diagonal2}&=\pm\dim(\la(0)\cap\mu(0)\cap X(0))\mp\dim(\la(1)\cap\mu(1)\cap X(1))\\
\label{e:maslov-diagonal3}&=\pm\dim(\la(0)\cap\mu(0)\cap Y(0))\mp\dim(\la(1)\cap\mu(1)\cap Y(1)).
\end{align}
\end{proposition}

\begin{proof} (a) By Theorem \ref{t:triangle-pairs}.a, $a(s)$ and $b(s)$ are well-defined integers, and $a(s)=-b(s)=\dim(\la(s)\cap\mu(s)\cap X(s))-\dim(\la(s)\cap\mu(s)\cap Y(s))$. By \cite[Remark IV.4.31]{Ka95}, $a(s)$ and $b(s)$ are constant integers.
\newline (b) Fix an $s_0\in[0,1]$. Set $n\ :=\dim(\la(s_0)\cap\mu(s_0))$. By Lemma \ref{l:diag-decomposition}, there is an isotropic subspace $V_{s_0}$ such that
$V_{s_0}=V_{s_0}\cap X(s_0)\oplus V_{s_0}\cap Y(s_0)$, and $Z=V_{s_0}\oplus(\la(s_0)+\mu(s_0)$. By Lemma \ref{l:family-isotropic-separate}, Lemma \ref{l:condition} and \cite[Lemma 3.1.1]{BoZh18}, there is a connected open neighborhood $U_{s_0}$ of $s_0$ in $[0,1]$ and a continuous family $V_{s_0}: U(s_0)\to\Ss(Z)$ of isotropic subspaces of $(Z,\w(s))$ such that for each $s\in U_{s_0}$, there hold that $V_{s_0}(s_0)=V_{s_0}$, $V_{s_0}(s)=V_{s_0}(s)\cap X(s)\oplus V_{s_0}(s)\cap Y(s)$, $V_{s_0}(s)\cap\la(s)=V_{s_0}(s)\cap\mu(s)=\{0\}$, and
\begin{align*}
Z&=V_{s_0}(s)\oplus (V_{s_0}(s))^{\w(s)}\cap\la(s)\oplus\mu(s)\\
&=V_{s_0}(s)\oplus\la(s)\oplus (V_{s_0}(s))^{\w(s)}\cap\mu(s).
\end{align*}

Set $\la_0(s)\ :=\la(s)\cap(V_{s_0}(s)+\mu(s))$, $\mu_0(s)\ :=\mu(s)\cap(V_{s_0}(s)+\la(s))$ and $Z_0(s)\ :=V_{s_0}(s)+\la_0(s)$ for each $s\in U_{s_0}$.
By Proposition \ref{p:decomposition} and Corollary \ref{c:diag-decomposition}, there hold that $\la(s)\cap\mu(s)=\la_0(s)\cap\mu_0(s)$, and $Z_0(s)$ is a symplectic subspace of $(Z,\w(s))$ such that
\begin{align*}
Z_0(s)&=V_{s_0}(s)\oplus\la_0(s)=V_{s_0}(s)\oplus\mu_0(s)\\
&=Z_0(s)\cap X(s)\oplus Z_0(s)\cap Y(s).
\end{align*}
By Theorem \ref{t:localization}, there hold that $\{Z_0(s)\}_{s\in U_{s_0}}$ is a continuous family, and the families $\{\la_0(s)\in\Ll(Z_0(s),\w(s)|_{Z_0(s)})\}_{s\in U_{s_0}}$, $\{\mu_0(s)\in\Ll(Z_0(s),\w(s)|_{Z_0(s)})\}_{s\in U_{s_0}}$ are continuous ones. By Corollary \ref{c:diag-decomposition}, there holds
$\la_0(s)=\la_0(s)\cap X(s)\oplus \la_0(s)\cap Y(s)$. Then
the families $\{\la_0(s)\cap X(s)\}_{s\in U_{s_0}}$ and $\{\la_0(s)\cap Y(s)\}_{s\in U_{s_0}}$ are continuous ones of linear subspaces of $Z$. By \cite[Corollary IV.2.6]{Ka95}, the integers $\dim(\la_0(s)\cap X(s))$ and $\dim(\la_0(s)\cap Y(s))$ are constants for $s\in U_{s_0}$.
By Theorem \ref{t:localization} and Lemma \ref{lm:Q_diagonal}, there holds that 
\begin{align*}
\Mas_{\pm}&\{\la(s),\mu(s);s\in[s_1,s_2]\}=\pm(\frac{1}{2}(n-\dim(\la_0(s_2)\cap\mu_0(s_2)))\\
&\mp(\frac{1}{2}(n-\dim(\la_0(s_1)\cap\mu_0(s_1)))\\
=&\pm\frac{1}{2}(\dim(\la_0(s_1)\cap\mu_0(s_1))-\dim(\la_0(s_2)\cap\mu_0(s_2)))\\
=&\pm\frac{1}{2}(\dim(\la(s_1)\cap\mu(s_1))-\dim(\la(s_2)\cap\mu(s_2)))\\
=&\pm(\dim(\la_0(s_1)\cap\mu_0(s_1)\cap X(s_1))-\dim(\la_0(s_2)\cap\mu_0(s_2)\cap X(s_2)))\\
=&\pm(\dim(\la(s_1)\cap\mu(s_1)\cap X(s_1))-\dim(\la(s_2)\cap\mu(s_2)\cap X(s_2))).
\end{align*}
By the compactness of $[0,1]$ and the path additivity of the Maslov index, we get (\ref{e:maslov-diagonal1}) and (\ref{e:maslov-diagonal2}).
The equation (\ref{e:maslov-diagonal3}) follows from (a) and \eqref{e:maslov-diagonal2}. 
\end{proof}

We now consider the Maslov index for the path with fixed diagonal part.

\begin{proposition}\label{p:fixed_intersection}
Let $(Z,\omega)$ be a symplectic Banach space such that $Z=X\oplus Y$ with $X,Y\in \Ll(Z)$.
Let $\{(\lambda(s),\mu(s))\}_{s\in[0,1]}$ be a path in $\Ff\Ll_0(Z)$.
Assume that for each $s\in[0,1]$, there hold that $\lambda(s)\cap Y= W_1$, $\mu(s)\cap Y=W_2$,  $\lambda(s)+Y=W_1^{\omega}$ and $\mu(s)+Y=W_2^{\omega}$.
Set $W\ :=W_1+W_2$. Set $Q(s)\ :=Q(\pi_{W}(\mu(s)),\pi_{W}(Y);\pi_{W}(\lambda(s)))$. Then the following hold.
\newline (a) We have $\dim W^{\w}/W=2\dim Y/W<+\infty$, $Y/W\in\Ll(W^{\w}/W)$, $W$ is $\omega$-closed, and
\begin{align}
\label{e:maslov-fixed-diagonal-s2} \Mas_{\pm}\{\lambda(s),\mu(s);\;s\in[0,1]\}=&\pm m^\pm(Q(1))\mp m^\pm(Q(0)),\\
\label{e:maslov-fixed-diagonal-s3} \dim(\lambda(s)\cap\mu(s))=&\dim\ker Q(s)+\dim W_1\cap W_2
\end{align}
for each $s\in[0,1]$.
\newline (b) We have $\dom(\lambda(s)\cap W^{\w})=\dom(\mu(s)\cap W^{\w})=W^{\w}\cap X=\dom(\lambda(s))\cap\dom(\mu(s))$, $\pi_W(\lambda(s))$ and $\pi_W(\mu(s))$ are the graph of some operator $G(s)$ and $H(s)$ respectively, where $G(s),H(s)\in\Hom(\pi_W(X), \pi_W(Y))$, and for each $s\in[0,1]$,
\begin{align}\label{e:maslov-fixed-diagonal-s4}
	Q(s)&([x_1]+H(s)[x_1],[x_2]+H(s)[x_2])=\omega(s)([x_1],(G(s)-H(s))[x_2])\\
	\label{e:maslov-fixed-diagonal-s5}
	&=\omega(s)(x_1,(\lambda(s)-\mu(s))x_2)
\end{align}
with $x_1,x_2\in \dom(\lambda(s))\cap\dom(\mu(s))$.
\end{proposition}

\begin{proof} (a) 1. By Lemma \ref{l:Lag-to-triangle}, we have $(W_1^\omega\cap X+W_1,W_2^\omega\cap X+W_2)\in\Ff\Ll_0(Z)$. By Corollary \ref{c:diagonal-index}, we have $\dim Y/W<+\infty$. By Theorem \ref{t:triangle-pairs}, $W$ is an $\omega$-closed subspace.
	
2. We have $W_1\cap \mu(s)=W_1\cap Y\cap\mu(s)=W_1\cap W_2$. Since $W_1$ is $\omega$-closed and $W_1\subset\lambda\subset W_1^\omega$, by \cite[Theorem 3.3.3]{BoZh18}, we have
\begin{equation}\label{e:maslov-r-s}
\Mas_{\pm}\{\lambda(s),\mu(s)\}=\Mas_{\pm}\{\pi_{W_1}(\lambda(s)),\pi_{W_1}(\mu(s))\}.
\end{equation}

By Lemmas \ref{l:diag-omega} and \ref{l:reduction-separate}, we have $\pi_{W}(Z)=\pi_{W}(X)\oplus\pi_{W}(Y)$. Since $W$ is $\omega$-closed, we have $\pi_W(Y)=Y/W$ and $\pi_W(Z)=W^{\w}/W$. By \cite[Proposition 1]{BoZh13}, we have $Y/W\in\Ll(W^{\w}/W)$. 
Since $\lambda(s)\cap Y= W_1$, $\mu(s)\cap Y=W_2$,  $\lambda(s)+Y=W_1^{\omega}$ and $\mu(s)+Y=W_2^{\omega}$, by Lemma \ref{l:g-index-reduction}, we have
\[W^{\w}/W=\pi_W(\la(s))\oplus \pi_W(Y)=\pi_W(\mu(s))\oplus\pi_W(Y).\]
Therefore we can view $\pi_W(\la(s))$ as the graph of a linear operator $L(s)$, where $L(s)\in\Hom(\pi_W(\mu(s)),\pi_W(Y))$.

Since $W_1\subset\lambda(s)\subset W_1^\omega$, by Lemma \ref{l:g-index-reduction}.a, we have
\[\pi_{W_1}(\lambda(s))\cap\pi_{W_1}(W_2)=\pi_{W_1}(\lambda(s)\cap W_2)=\pi_{W_1}(W_1\cap W_2)=\{0\}.\]
Since $W$ is an $\omega$-closed subspace, $\pi_{W_1}(W_2)=W/W_1$ is an $\omega$-closed subspace of $W_1^{\w}/W_1$. Since $\pi_{W_1}(W_2)\subset \pi_{W_1}(\mu(s))\subset\pi_{W_1}(W_2^\omega)$, by \cite[Theorem 3.3.3, Lemma 1.4.3]{BoZh18} and (\ref{e:maslov-r-s}), we have
\begin{align*}
\Mas_{\pm}\{\lambda_s,\mu_s\}&=\Mas_{\pm}\{\pi_{W_1}(\lambda(s)),\pi_{W_1}(\mu(s))\}\\
&=\Mas_{\pm}\{\pi_{\pi_{W_1}(W_2)}(\pi_{W_1}(\lambda(s))),\pi_{\pi_{W_1}(W_2)}(\pi_{W_1}(\mu(s)))\}\\
&=\Mas_{\pm}\{\pi_W(\lambda(s)),\pi_W(\mu(s))\}.
\end{align*}
By \cite[Lemma 1.4.6]{BoZh18}, for each $s\in[0,1]$, we have
\begin{align*}
	\dim(\lambda(s)\cap\mu(s))=&\dim(\pi_{W_1}(\lambda(s))\cap\pi_{W_1}(\mu(s)))+\dim W_1\cap W_2\\
	=&\dim(\pi_W(\lambda(s))\cap\pi_W(\mu(s)))+\dim W_1\cap W_2\\
	=&\dim\ker Q(s)+\dim W_1\cap W_2.	
\end{align*}

3. Since $\pi_W(Y)=Y/W$, by Step 1 and Step 2, we have $Y/W\in\Ll(W^{\w}/W)$. So we have $\dim W^{\w}/W=2\dim Y/W$. 
By \cite[Lemma 2.3.2]{BoZh18} and Step 2, we have
\begin{align*}
	\Mas_{\pm}\{\lambda(s),\mu(s)\}
	&=\Mas_{\pm}\{\pi_W(\lambda(s)),\pi_W(\mu(s))\}\\
	&=\pm(m^\pm(Q(1))-m^\pm(Q(0))).
\end{align*}
\newline (b) We have
\begin{align*}
	W^\omega\cap X
	&=W_1^\omega\cap X\cap W_2^\omega\cap X=\dom(\lambda(s))\cap\dom(\mu(s)),\\
    \pi_W(X)&=(W^\omega\cap X+W)/W=\pi_W(\dom(\lambda(s))\cap\dom(\mu(s))).
\end{align*}
By the proof of (a), we have
\[\pi_W(Z)=\pi_W(X)\oplus \pi_W(Y)=\pi_W(\la(s))\oplus \pi_W(Y)= \pi_W(\mu(s))\oplus\pi_W(Y).\]
So $\pi_W(\lambda(s))$ and $\pi_W(\lambda(s))$ are the graph of some operators $G(s)$ and $H(s)$ respectively, where $G(s),H(s)\in\Hom(\pi_W(X), \pi_W(Y))$.
Since $W\subset Y\subset W^{\w}$, by \cite[Lemma A.1.1]{BoZh18}, we have
\begin{align*}
\dom(\lambda(s)\cap W^{\w})+Y&=(\la(s)\cap W^{\w}+Y)\cap X+Y\\
&=(\la(s)\cap W^{\w}+Y)\cap (X+Y)\\
&=(\la(s)+Y)\cap W^{\w}\\
&=W_1^{\w}\cap W^{\w}=W^{\w}.
\end{align*}
So we have $\dom(\lambda(s)\cap W^{\w})=W^{\w}\cap X=\dom(\lambda(s))\cap\dom(\mu(s))$. Similarly we have $\dom(\mu(s)\cap W^{\w})=W^{\w}\cap X=\dom(\lambda(s))\cap\dom(\mu(s))$.

Let $x_1,x_2\in W^{\w}\cap X$ be two vectors. By the definition of $Q(s)$, we have
\begin{align*}
Q(s)&([x_1]+H(s)[x_1],[x_2]+H(s)[x_2])=\omega(s)([x_1],(G(s)-H(s))[x_2])\\
&=\omega(s)(x_1,(\lambda(s)-\mu(s))x_2).
\end{align*}
\end{proof}

Now we prove Theorem \ref{t:main-triangle}.

\begin{proof}[Proof of Theorem \ref{t:main-triangle}] By Proposition \ref{p:continuity-data}, for each $s_0\in[0,1]$, there exists a $d>0$ and a continuous family 
\begin{align*}
	\{(M(s,t), N(s,t))\in\Ff\Ll_0(Z,\w(s))\}_{(s,t)\in((s_0-d,s_0+d)\cap[0,1])\times [0,1]}
\end{align*}
such that
$\lambda(s)=M(s,1)$, $\mu(s)=N(s,1)$, $\dom (M(s,t))=\dom (\lambda(s))$, $
M(s,t)\cap Y(s)=A(s)$, $\dom (N(s,t))=\dom (\mu(s))$, $N(s,t)\cap Y(s)=B(s)$, $\gamma(s)=M(s,0)$, $\delta(s)=N(s,0)$. 
Let $s_1,s_2\in(s_0-d,s_0+d)\cap[0,1]$ be two real numbers satisfying $s_1\le s_2$. By \cite[Theorem 3.1.5]{BoZh18}, we have
\begin{align*}
	\Mas_{\pm}&\{\la(s),\mu(s);\;\w(s),s\in[s_1,s_2]\}=\Mas_{\pm}\{\gamma(s),\delta(s);\;\w(s),s\in[s_1,s_2]\}\\&+\Mas_{\pm}\{M(s_2,t),N(s_2,t);\;\w(s_2),t\in[0,1]\}\\&-\Mas_{\pm}\{M(s_1,t),N(s_1,t);\;\w(s_1),t\in[0,1]\}.
\end{align*}
By Proposition \ref{p:maslov-diagonal}, we have
\begin{align*}
	\Mas_{\pm}&\{\gamma(s),\delta(s);\w(s);\;s\in[s_1,s_2]\}\\
	&=\pm\dim(\gamma(s_1)\cap\delta(s_1)\cap Y(s_1))\mp\dim(\gamma(s_2)\cap\delta(s_2)\cap Y(s_2))\\
	&=\pm\dim(\lambda(s_1)\cap\mu(s_1)\cap Y(s_1))\mp\dim(\lambda(s_2)\cap\mu(s_2)\cap Y(s_2)).
\end{align*}
For each pair $x_1,x_2\in\dom(\gamma(s))\cap\dom(\delta(s))$ and each $s\in[0,1]$, we have 
\begin{align*}
	\omega(s)(x_1,(\gamma(s)-\delta(s))x_2)=0.
\end{align*}
By Theorem \ref{t:triangle-pairs}.a, we have $\lambda(s)+Y(s)=A(s)^\omega(s)$ and $\mu(s)+Y(s)=B(s)^\omega(s)$. By Proposition \ref{p:fixed_intersection}, we have
\begin{align*}
	\Mas_{\pm}\{M(s,t),N(s,t);\w(s);\;t\in[0,1]\}= \pm m^{\pm}(Q(s)).
\end{align*}
Then we have
\begin{align*}
	\Mas_{\pm}&\{\la(s),\mu(s);\w(s);\;s\in[s_1,s_2]\}=\pm\dim(\lambda(s_1)\cap\mu(s_1)\cap Y(s_1))\mp\\&\dim(\lambda(s_2)\cap\mu(s_2)\cap Y(s_2))\pm m^{\pm}(Q(s_2))\mp m^{\pm}(Q(s_1)).
\end{align*}
By the compactness of $[0,1]$ and the path additivity of the Maslov index, we get \eqref{e:main-triangle}.
Since $\dom(\lambda(s))\cap\dom(\mu(s))\cap (A(s)+B(s))=0$ and $\dom(\lambda(s))\cap\dom(\mu(s))\subset (A(s)+B(s))^\omega(s)$, we have
$\dim\pi_{A(s)+B(s)}\ker Q(s)=\dim\ker Q(s)$. By Proposition \ref{p:fixed_intersection}, we have \eqref{e:maslov-nullity}.
\end{proof}

%% file: applications.tex
\section{Applications}

\subsection{Maslov type index}

Firstly we recall notion of the Maslov-type index in the finite dimensional case.

\begin{definition}\label{d:symplectic-linear-map}
	Let $(X_l,\omega_l)$, $l=1,2$ be symplectic vector spaces. We define 
	\begin{align}\label{e:linear-symplectic-action}
		\Sp(X_1,X_2)\ :=&\{M\in\Hom(X_1,X_2);\;M^*\omega_2=\omega_1, LX_1=X_2\},\\
		\label{e:linear-symplectic-group}
		\Sp(X_1)\ :=&\Sp(X_1,X_1).
	\end{align}
\end{definition}

For each $M\in\Sp(X_1,X_2)$, $L$ is an isomorphism. Then we have $\Gr(M)\oplus\Gr(-M)=X_1\times X_2$. 
By \cite[Proposition 1]{BoZh13}, we have $M\in\Ll(X_1\times X_2, (-\omega_1)\oplus\omega_2)$.

For each $\tau >0$ and a finite dimensional symplectic vector space $(X,\omega)$, we define
\begin{equation}\label{e:symplectic-path-I}
	\mathcal{P}_{\tau}(X)\ :=\{\gamma\in C([0,\tau],\Sp(X));
	\gamma(0)=I_X\}.
\end{equation}

\begin{definition}\label{d:Maslov-type-index}(cf. \cite[Definition 4.6]{Zhu06})
	Let $(V_l,\omega_l)$, $l=1,2$ be two finite-dimensional symplectic
	vector spaces. Then $(V=V_1\times V_2, (-\omega_1)\oplus\omega_2)$ is
	a symplectic vector space. Let $W\in \mathcal{L}(V)$. Let
	$\gamma(t)$, $0\leq t\leq\tau$ be a path in $\Sp(V_1,V_2)$. The
	{\em Maslov-type index} $i_W(\gamma)$ is defined to be
	$\Mas\{\Gr\circ\gamma, W\}$. If $P\in \Sp(V_1,V_2)$, we define
	$i_P(\gamma)\ :=i_{\Gr(P)}(\gamma)$. If
	$(V_1,\omega_1)=(V_2,\omega_2)$ and
	$\gamma\in\mathcal{P}_{\tau}(V_1)$, we denote by
	$i_z(\gamma)\ :=i_{zI_{V_1}}(\gamma)$, $\nu_z(\gamma)\ :=\dim
	\ker(\gamma(\tau)-zI_{V_1})$, $\nu_{W}(\gamma)\ :=\dim
	(\Graph(\gamma(\tau))\cap W)$ and $\nu_{P}(\gamma)\ :=\dim \ker
	(\gamma(\tau)-P)$ for $z\in S^1$ and $P\in\Sp(V_1)$.
\end{definition}

Similarly, we can define the Maslov-type index in symplectic Banach space.

\begin{definition}\label{d:Maslov-type-index-Banach}
	Let $a<b$ be two real numbers. Let $(X_l,\omega_l(s))$, $l=1,2$ be symplectic Banach spaces with continuously varying symplectic structures $\{\omega_l(s);\;s\in[a,b]\}$. Then $(X_1\times X_2, (-\omega_1(s))\oplus\omega_2(s))$ is
	a symplectic Banach space for each $s\in[a,b]$. Let $c\: =\{c(s)\in \Ll(X_1\times X_2, (-\omega_1(s))\oplus\omega_2(s));\;s\in[a,,b]\}$ be a path. Let
	$\{\gamma(t)\in\Sp((X_1,\omega_1(s)),(X_2,\omega_2(s)));\;s\in[a,b]\}$ be a path. Assume that $\Index(\Gr(\gamma(s)),W(s))=0$ for each $s\in[a,b]$. The
	{\em Maslov-type index} $i_{\pm,c}(\gamma)$ is defined to be
	$\Mas_\pm\{\Graph\circ\gamma, c\}$, $i_c(\gamma)$ is defined to be
	$\Mas\{\Graph\circ\gamma, c\}$ respectively. We define $\nu_c(\gamma)(s)\ :=\dim(\Gr(\gamma(s))\cap c(s))$ for each $s\in[a,b]$. If $c(s)=W$ for each $s\in[a,b]$,  we define $i_{\pm,W}(\gamma)\ :=i_{\pm,c}(\gamma)$, 
	$i_W(\gamma)\ :=i_c(\gamma)$ and $\nu_W(\gamma)(s)\ :=\nu_c(\gamma)(s)$.
\end{definition}

\subsection{A formula of Maslov-type index} Let $X,Y$ be two Banach spaces over $\KK$, $\KK$ is $\R$ or $\C$. A bounded sesquilinear form $\Omega\colon X\times Y\to \KK$ define two bounded linear maps $L_\Omega\colon X\to Y^*$ and $R_\Omega\colon Y\to X^*$ by 
\begin{align}\label{e:form-operator}
	\Omega(x,y)=(L_\Omega x)y=\overline{(R_\Omega y)x}\quad\text{ for all }x\in X\text{ and }y\in Y.
\end{align} 

Let $(X_l,\omega_l)$, $l=1,2$ be symplectic Banach spaces. Then we have
\begin{align}\label{e:linear-symplectic-map}
	M^*\omega_2=\omega_1\Leftrightarrow M^*L_{\omega_2}M=L_{\omega_1}
\end{align}
for $M\in\Bb(X_1,X_2)$.

Let $X$ and $Y$ be two vector spaces over $\KK$. Given a nondegenerate sesquilinear form $\Omega\colon X\times Y\to\KK$, there is a natural symplectic structure $\omega$ on $Z\ :=X\times Y$ defined by (cf. \cite{Ben72}, \cite[Definition 4.2.1]{BoZh18}, \cite{Ek90})
\begin{align}\label{e:natural-symplectic-structure}
	\omega((x_1,y_1),(x_2,y_2))\ :=\Omega(x_1,y_2)-\overline{\Omega(x_2,y_1)}
\end{align}
for all $x_1,x_2\in X$ and $y_1,y_2\in Y$. For each subsets $\lambda$ of $X$ and $\mu$ 
of $Y$, the {\em annihilator} 
$\lambda^{\Omega,r}$ and $\mu^{\Omega,l}$ (cf. \cite[Definition 1.1.1.a]{BoZh18}) are defined by 
\begin{align}\label{e:right-side-annihilator}
	\lambda^{\Omega,r}\ :=&\{y\in Y;\ \Omega(x,y)=0\text{ for all }x\in\lambda\},\\
	\label{e:left-side-annihilator}
	\mu^{\Omega,l}\ :=&\{x\in X;\ \Omega(x,y)=0\text{ for all }y\in\mu\}.
\end{align}
We call $\omega$ defined by \eqref{e:natural-symplectic-structure} {\em the symplectic structure associated to $\Omega$}. 
Then there are two nature Lagrangian subspaces $X\times\{0\}$ and $\{0\}\times Y$ of $(Z,\omega)$.

If $X, Y$ are Banach spaces and $\Omega$ is bounded, by \eqref{e:form-operator}, the induced operator of $\omega$ is given by
\begin{align}\label{e:L-omega}
	L_\omega=\left(\begin{array}{cc}
		0&-R_\Omega  \\
		L_\Omega&0 
	\end{array}\right).
\end{align}

By definition, we have the following fact(cf. \cite[Lemma 1.1.2]{Lo02}). 

\begin{lemma}\label{l:symplectic-block-form}
	Let $X_l$ and $Y_l$, $l=1,2$ be four Banach spaces over $\KK$, $\KK$ is $\R$ or $\C$. Let $\Omega_l\colon X_l\times Y_l\to\KK$, $l=1,2$ be nondegenerate sesquilinear forms. Let $\omega_l$ be the symplectic structure of $Z_l\ :=X_l\times Y_l$ associated to $\Omega_l$ for each $l=1,2$. Suppose that $M\in\Bb(Z_1,Z_2)$ has the block form
	\begin{align*}
		M=\left(\begin{array}{cc}
		A & B \\
		C & D
		\end{array}\right).
	\end{align*}
	\newline (a) There holds $M^*\omega_2=\omega_1$ if and only if $A^*R_{\Omega_2}C=C^*L_{\Omega_2}A$, $B^*R_{\Omega_2}D=D^*L_{\Omega_2}B$, and $A^*R_{\Omega_2}D-C^*L_{\Omega_2}B=R_{\Omega_1}$.
	\newline (b) Especially when $C=0$, there hold $M^*\omega_2=\omega_1$ if and only if $B^*R_{\Omega_2}D=D^*L_{\Omega_2}B$ and $A^*R_{\Omega_2}D=R_{\Omega_1}$, $M$ is an isomorphism if and only if both $A$ and $D$ are isomorphisms. When $B=0$, there hold $M^*\omega_2=\omega_1$ if and only if $A^*R_{\Omega_2}C=C^*L_{\Omega_2}A$ and $A^*R_{\Omega_2}D=R_{\Omega_1}$, $M$ is an isomorphism if and only if both $A$ and $D$ are isomorphisms.
\end{lemma}

We have the following generalization of \cite[Theorem 2.2]{Zhu06}.

\begin{theorem}\label{t:Maslov-type-index-triangular-form}
	Let $a<b$ be two real numbers. Let $X_l$ and $Y_l$, $l=1,2$ be four Banach spaces over $\KK$, $\KK$ is $\R$ or $\C$. Let $\{\Omega_l(s)\colon X_l\times Y_l\to\KK;\;s\in[a,b]\}$, $l\in\{1,2\}$ be two paths of nondegenerate sesquilinear forms. Let $\omega_l(s)$, $l=1,2$ be the symplectic structure of $Z_l\ :=X_l\oplus Y_l$ associated to $\Omega_l(s)$ for each $l=1,2$ and $s\in[a,b]$. Suppose that two paths $\{M(s)\in\Sp((Z_1,\omega_1(s)),(Z_2,\omega_2(s)));\;s\in[a,b]\}$
	and $\{M_0(s)\in\Hom((Z_1,\omega_1(s)),(Z_2,\omega_2(s)));\;s\in[a,b]\}$ has the block form
	\begin{align*}
		M(s)=\left(\begin{array}{cc}
			A(s) & B(s) \\
			C(s) & D(s)
		\end{array}\right),\quad
		M_0(s)=\left(\begin{array}{cc}
			A(s) & 0 \\
			0 & D(s)
		\end{array}\right),\quad s\in[a,b].
	\end{align*}
	Let $\{R_1(s)\in \Ss(X_1(s)\times X_2(s));\;s\in[a,b]\}$ and $\{R_2(s)\in \Ss(Y_1(s)\times Y_2(s));\;s\in[a,b]\}$ be two paths such that 
	\begin{align*}
		W(s)\ :=R_1(s)\oplus R_2(s)\in\Ll(Z(s))
	\end{align*}
	for each $s\in[a,b]$, where $Z(s)\ :=(Z_1\times Z_2,(-\omega_1(s))\oplus\omega_2(s))$. Assume that $\Index(\Gr(M(s)),W(s))=0$ for each $s\in[a,b]$.
	\newline (a) Assume that $C(s)=0$ for each $s\in[a,b]$. For each $s\in[a,b]$, we define the symmetric form $g(s)$ by 
	\begin{align*}
		g(s)(y_1,y_2)=\overline{\Omega_2(B(s)y_2,D(s)y_1)}		
	\end{align*}
	for $(y_l,D(s)y_l)\in R_2(s)$, $l\in[1,2]$.
	Then we have		
	\begin{align}\nonumber
		i_{\pm,W}(M)=&\pm\dim(\Gr(A(a))\cap R_1(a))\mp\\
		\label{e:Maslov-type-index-uptriangular}
		&\dim(\Gr(A(b))\cap R_1(b))\pm m^{\mp}(g(b))\mp m^{\mp}(g(a)),\\
		\label{e:Maslov-type-nullity-uptriangular}
		\nu_W(\gamma)(s)=&\dim\ker g(s)+\dim(\Gr(A(s))\cap R_1(s)).		
	\end{align}
	\newline (b) Assume that $B(s)=0$ for each $s\in[a,b]$. For each $s\in[a,b]$, we define the symmetric form $h(s)$ by 
	\begin{align*}
		h(s)(x_1,x_2)=\Omega_2(A(s)x_1,C(s)x_2)		
	\end{align*}
	for $(x_l,A(s)x_l)\in R_1(s)$, $l\in[1,2]$.
	Then we have		
	\begin{align}\nonumber
		i_{\pm,W}(M)=&\pm\dim(\Gr(D(a))\cap R_2(a))\mp\\
		\label{e:Maslov-type-index-downtriangular}
		&\dim(\Gr(D(b))\cap R_2(b))\pm m^{\pm}(h(b))\mp m^{\pm}(h(a)),\\
		\label{e:Maslov-type-nullity-downtriangular}
		\nu_W(\gamma)(s)=&\dim\ker h(s)+\dim(\Gr(D(s))\cap R_2(s)).		
	\end{align}
\end{theorem}

\begin{proof} 
	Let $s\in[a,b]$ be a real number.
	We have two paths $\{ X_1(s)\times X_2(s)\in\Ll;\;s\in[a,b]\in\Ll(Z(s))\}$ and $\{R_2(s)\in Y_1(s)\times Y_2(s);\;s\in[a,b]\in\Ll(Z(s))\}$, and $Z(s)=X_1(s)\times X_2(s)\oplus Y_1(s)\times Y_2(s)$. By Corollary \ref{c:diagonal-index}, we have
	\begin{align}\label{e:graph-sum-index}
		\Index(\Gr(A(s)),R_1(s))+\Index(\Gr(D(s)),R_2(s))=0.
	\end{align}
	\newline (a) We have 
	\begin{align*}
		\Gr(M(s))\cap W(s)=&\{(x,y)\in Z(s);\;
		(x,A(s)x+B(s)y)\in R_1(s),\\
		&(y,Dy)\in R_2(s)\}.
	\end{align*}
	Then we have 
	\begin{align*}
		\Gr(M(s))\cap W(s)\cap(X_1(s)\times X_2(s))&=\Gr(A(s))\cap R_1(s),\\
		\Gr(M(s))\cap W(s)\cap(Y_1(s)\times Y_2(s))&=\Gr(D(s))\cap R_2(s).
	\end{align*}
	By Lemma \ref{l:symplectic-block-form}, the form $g(s)$ is symmetric, and the form $Q(s)$ on $\Gr(D(s))\cap R_2(s)$ defined in Theorem \ref{t:main-triangle} is calculated by
	\begin{align*}
		&Q(s)((y_1,D(s)y_1),(y_2,D(s)y_2))\\
		=&
		(-\omega_1(s)\oplus\omega_2(s))(0,y_1,0,D(s)y_1),(x_2,y_2,A(s)x_2+B(s)x_2,D(s)y_2)\\
		=&\overline{\Omega_1(x_2,y_1)}-\overline{\Omega_2(A(s)x_2+B(s)y_2,D(s)y_1)}\\
		=&-\overline{\Omega_2(B(s)y_2,D(s)y_1)}=-g(s)(y_1,y_2)
	\end{align*}
	for $(y_l,D(s)y_l)\in R_2(s)$, $l\in\{1,2\}$ and $x_2\in X_1(s)$.
	Then (a) follows from Theorem \ref{t:main-triangle}.
	\newline (b) follows from (a) by interchanging the spaces $X(s)$ and $Y(s)$ and replacing $\Omega_l(s)(x_l,y_l)$ by $\overline{\Omega_l(s)(x_l,y_l)}$
	 for $s\in[a,b]$, $l\in\{1,2\}$, $x_l\in X_l(s)$, $y_l\in Y_l(s)$.  
\end{proof}

\subsection{Splitting numbers}

Let $X$ and $Y$ be two Banach spaces over $\C$. Given a nondegenerate sesquilinear form $\Omega\colon X\times Y\to\C$, let $\omega$ be the symplectic structure of $Z\ :=X\times Y$ associated to $\Omega$. 

\begin{definition}\label{d:splitting-numbers}(\cite[Definition 9.1.4]{Lo02}, \cite[Definition 4.1]{LoZh00})
	Let $M\in\Sp(Z)$ be a linear symplectic map, and $z\in S^1$ be a complex number on the unit circle. Assume that $\Index(M-zI)=0$ holds. Denote by $\mi$ the imaginary unit $\sqrt{-1}$. We define 
	\begin{align}\label{e:splitting-numbers-minus}
		S^{\pm}_M(z)=&S^{\pm}_{-,M}(z)\ :=\lim\limits_{t\to 0^+}i_1(Mz^{-1}e^{\mp \mi s};\;s\in[0,t]),\\
		\label{e:splitting-numbers-plus}
		S^{\pm}_{+,M}(z)\ :=&-\lim\limits_{t\to 0^+}i_{-,1}(Mz^{-1}e^{\mp \mi s};\;s\in[0,t]).
	\end{align}
	We call them {\em splitting numbers} of $M$ at $z$.
\end{definition}

By \cite[Theorem 1.4.1]{Lo02}, a real symplectic matrix $M\in\Sp(2n)$ has a upper triangular normal form for eigenvalue $1$. The fact (see \cite[Lemma 5.5]{Gi10}) is used in the proof of Conley conjecture. 

We want to calculate the splitting numbers of a linear symplectic map in triangular form.

\begin{theorem}\label{t:splitting-numbers-up-triangular}
	Let $M\in\Sp(Z)$ be a linear symplectic map, and $z\in S^1$ be a complex number on the unit circle. Assume that $\Index(M-zI)=0$ holds and $M$ has the block form
	\begin{align*}
		M=\left(\begin{array}{cc}
			A & B \\
			C & D
		\end{array}\right).
	\end{align*}
	\newline (a) Assume that $C=0$. We define the symmetric form $g$ by 
	\begin{align*}
		g(y_1,y_2)=\overline{\Omega_2(By_2,Dy_1)}		
	\end{align*}
	for $y_l\in \ker(D-zI_Y)$, $l\in[1,2]$.
	Then we have		
	\begin{align}\label{e:splitting-numbers-up-minus}
		S^{\pm}_M(z)=&\dim\ker(A-zI_X)- m^{-}(g),\\
		\label{e:splitting-numbers-up-plus}
		S^{\pm}_{+,M}(z)=&\dim\ker(A-zI_X)- m^{+}(g),\\
		\label{e:splitting-nullity-uptriangular}
		\nu_W(\gamma)(s)=&\dim\ker g+\dim\ker(A-zI_X).		
	\end{align}
	\newline (b) Assume that $B=0$. We define the symmetric form $h$ by 
	\begin{align*}
		h(x_1,x_2)=\Omega_2(Ax_1,Cx_2)		
	\end{align*}
	for $x_l\in \ker(A-zI_X)$, $l\in[1,2]$.
	Then we have		
	\begin{align}\label{e:splitting-numbers-down-minus}
		S^{\pm}_M(z)=&\dim\ker(D-zI_Y)- m^{+}(h),\\
		\label{e:splitting-numbers-down-plus}
		S^{\pm}_{+,M}(z)=&\dim\ker(D-zI_Y)- m^{-}(h),\\
		\label{e:splitting-nullity-downtriangular}
		\nu_W(\gamma)(s)=&\dim\ker h+\dim\ker(A-zI_Y).		
	\end{align}
\end{theorem}

\begin{proof}
	By Theorem \ref{t:Maslov-type-index-triangular-form}.
\end{proof}  
  
\subsection{Dependence of iteration theory on triangular frames}

The independence of iteration theory on frames for closed geodesics follows from \cite[Theorem 1.5]{DuLo09} and \cite[Section 12.1]{Lo02}. Here we study the dependence of iteration theory on triangular frames. Then we recover the independence of iteration theory on frames for closed orbits by \cite[Corollary 5.1, Corollary 2.1, Theorem 2.3]{Zhu06} (see also \cite[Proposition 3.2, Theorem 1.1]{ZhuGs21} for a special case) and Theorem \ref{t:iteration-frame} below. The idea is used by \cite[Theorem 1.3]{ZhuGs21}, where he separate it into orientable case and nonorientable case.

Let $(V,\omega)$ be a finite dimensional symplectic vector space. 
For each $\gamma\in\Pp_\tau(V)$ and and $P\in\Sp(V)$, we define the {\em iteration path} $\tilde\gamma(P)\colon[0,+\infty)\to\Sp(V)$ (cf. \cite[(4.81)]{Li19})
by
\begin{align}\label{e:iteration}
	\tilde\gamma(P)(t)\ :=P^j\gamma(t-j\tau)(P^{-1}\gamma(\tau))^j,\quad\forall t\le[j\tau,(j+1)\tau],j\in\N.
\end{align}
The {\em $k$-th iteration path} is defined by $\gamma(k,P)\ :=\tilde\gamma(P)|_{[0,k\tau]}$ for $k\in\N$ and $P\in\Sp(V)$. We denote by $i(k,\gamma;P)\ :=i_{P^k}(\gamma(k,P))=i_1(P^{-k}\gamma(k,P))$ (the second equality follows form \cite[Lemma 4.4]{Zhu06}) for $k\in\N$ and $P\in\Sp(V)$. Then we have the following result on the Maslov-type index for the iteration path. 

\begin{theorem}\label{t:iteration-frame}
	Let $(V,\omega)$ be a finite dimensional complex symplectic vector space. 
	Let $\gamma\in\Pp_\tau(V)$ be a symplectic path with $\gamma(\tau)=M$. Let $P\in\Sp(V)$ be a symplectic linear map. 
	\newline (a) The function $f(k,M,P)\ :=i(k,\gamma;P)-ki(1,\gamma;P)$ is well-defined, and
	\begin{align}\label{e:iteration-frame}
		f(k,M,P)=f(k,P^{-1}M,I_V)-f(k,P^{-1},I_V).
	\end{align}
	\newline (b) Let $X$ and $Y$ be two Lagrangian subspaces of $V$ and $\dim V=2n$. Assume that $\omega$ be the symplectic structure of $V$ associated to some nondegenerate sesquilinear form $\Omega\colon X\times Y\to\C$. Let $k\in\N$ be a positive integer. Suppose that $P$ and $P^k$ has the block form
	\begin{align*}
		P=\left(
		\begin{array}{cc}
			A & B \\
			C & D
		\end{array}\right), \quad
		P^k=\left(
		\begin{array}{cc}
			A(k) & B(k) \\
			C(k) & D(k)
		\end{array}\right).
	\end{align*} 
	\begin{itemize}
		\item [(i)] If $C=0$, we define the symmetric form $g(k)$ by 
		\begin{align*}
			g(k)(y_1,y_2)=\overline{\Omega(B(k)y_2,D^ky_1)}		
		\end{align*}
		for $y_l\in \ker(D^k-I_Y)$, $l\in[1,2]$. Then we have
		\begin{align}\label{e:Pk-upstriangular}
			P^k=&\left(
			\begin{array}{cc}
				A^k & \sum_{j=1}^{k-1}\left(\begin{array}{c}
					k-1 \\
					j
				\end{array}\right)A^jBD^{k-j} \\
				0 & D^k
			\end{array}\right)\quad\text{ for }k>1,\\
			\label{e:iteration-index-Pk-upstriangular}
			f(k,P,I_V)=&(1-k)n+k\dim\ker(A-I_X)-\dim\ker(A^k-I_X)\\ \nonumber
			&+m^-(g(k))-km^-(g(1)).
		\end{align}
		\item [(ii)] If $B=0$, we define the symmetric form $h(k)$ by 
		\begin{align*}
			h(k)(x_1,x_2)=\overline{\Omega(A^kx_1,C(k)x_2)}	
		\end{align*}
		for $x_l\in \ker(A^k-I_X)$, $l\in[1,2]$. Then we have
		\begin{align}\label{e:Pk-downstriangular}
			P^k=&\left(
			\begin{array}{cc}
				A^k & 0 \\
				\sum_{j=1}^{k-1}\left(\begin{array}{c}
					k-1 \\
					j
				\end{array}\right)D^jCA^{k-j} & D^k
			\end{array}\right)\quad\text{ for }k>1,\\
			\label{e:iteration-index-Pk-downstriangular}
			f(k,P,I_V)=&(1-k)n+k\dim\ker(A-I_X)-\dim\ker(D^k-I_X)\\ \nonumber
			&+m^+(h(k))-km^+(h(1)).
		\end{align}
	\end{itemize}
\end{theorem}

\begin{proof} (a) Let $\gamma_1\in\Pp_\tau(\C^{2n})$ be a symplectic paths with $\gamma_1(\tau)=P^{-1}$. 
We have a homotopy $H\colon[0,k\tau]^2\to\Sp(V)$ defined by 
	\begin{align*}
		H(s,t)\ := \gamma_1(k,I_V)(s)\gamma(k,P)(t)\text{ for each } (s,t)\in [0,k\tau]^2.
	\end{align*}
Then we have 
\begin{align}\nonumber
	i(k,\gamma_2;I_V)=&i_1(H(s,s);\;s\in[0,k\tau])\\
	\nonumber
	=&i_1(H(s,0);\;s\in[0,k\tau])+i_1(H(k\tau,t);\;t\in[0,k\tau])\\
	\label{e:index-iteration-frame}
	=&i(k,\gamma_1;I_V)+i(k,\gamma;P).
\end{align}
By the proof of \cite[Theorem 4.2]{LoZh00}, the right hand side of \eqref{e:iteration-frame} is well-defined. By \eqref{e:index-iteration-frame} and its special case $k=1$, we obtain (a). 
\newline (b) (i) Assume that $C$=0. By Lemma \ref{l:symplectic-block-form}, there is a symplectic path $\alpha\in\Pp_1(V)$ such that $\alpha(1)=P$ and $\alpha(t)$ has the block form
\begin{align*}
\alpha(t)=\left(
\begin{array}{cc}
	A(t) & B(t) \\
	0 & D(t)
\end{array}\right),\quad\text{for all }t\in[0,1].
\end{align*} 
Then \eqref{e:Pk-upstriangular} holds, and $\alpha(k,I_V)(t)$ has the block form
\begin{align*}
	\alpha(k,I_V)(t)=\left(
	\begin{array}{cc}
		A(t) & B(t) \\
		0 & D(t)
	\end{array}\right),\quad\text{for all }t\in[0,k].
\end{align*} 
By Theorem \ref{t:Maslov-type-index-triangular-form}, we have
\begin{align*}
	f(k,P,I_V)=&i(k,\alpha;I_V)-ki(1,\alpha;I_V)\\
	=&(n-\dim\ker(A^k-I_X)+m^-(g(k)))\\
	&-k(n-\dim\ker(A-I_X)+m^-(g(1)))
	\\
	=&(1-k)n+k\dim\ker(A-I_X)-\dim\ker(A^k-I_X)\\
	&+m^-(g(k))-km^-(g(1)).
\end{align*}

(ii) By interchanging the spaces $X$ and $Y$ and replacing $\Omega(x_l,y_l)$ by $\overline{\Omega(x_l,y_l)}$
for $l\in\{1,2\}$, $x_l\in X$, $y_l\in Y$ in (i).  
\end{proof}

Let $k\in\N$ be a natural numbers. Denote by $S_k\ :=\R/\sim_k$, where $x\sim_k y$ if and only if $k^{-1}(x-y)\in\Z$ for each pair $(x,y)\in\R^2$. Let $\pi\colon E\to S_1$ be a $C^2$ complex vector bundle of dimensional $n$. Then there is a $C^2$ frame $t\mapsto e(t)=\{e_1(t),\ldots,e_n(t)\}$ of $\pi^{-1}(t)$ for $t\in[0,1]$.
Define $a=(a_{ij})_{i,j=1,\ldots,n}\in \GL(n,\C)$ by
\begin{align}\label{e:frame-holonomy}
	e_i(1)=\sum_{j=1}^{n}a_{ji}e_j(0).
\end{align}
Then the frame $e(t)$ can be extended to a locally Lipschitz frame by 
\begin{align}\label{e:frame-holonomy-extended}
	e_i(t+1)=\sum_{j=1}^{n}a_{ji}e_j(t),\quad\forall t\in\R.
\end{align}
We assume our choice of frame make the extended frame defined by \eqref{e:frame-holonomy-extended} is of class $C^2$.

Let $\Ii$ be an $C^1$ index form of order $1$ on the Sobolev space $H^1(S_1;E)$ with positive definite highest order term. By Sobolev embedding theorem, we have $H^1(S_1;E)\subset C(S_1;E)$. For two sections $X,Y\in H^1(S_k;E)$, there are two maps $x=(x_1,\ldots,x_n),y=(y_1,\ldots,y_n)\in H^1(R;\C^n)$ such that 
\begin{align}\label{e:frame-section}
	X(t)=\sum_{i=1}^{n}x_i(t)e_i(t),\quad Y(t)=\sum_{i=1}^{n}y_i(t)e_i(t),\quad\forall t\in\R.
\end{align}
Denote by $\dot x\ :=\frac{dx}{dt}$. Denote by $\lla\cdot,\cdot\rra$ the standard inner product of $\C^n$. The index form $\Ii$ has the form
\begin{align}\label{e:index-form}
	\Ii(X,Y)=\int_{0}^{1}\left(\lla p\dot x+qx,\dot y\rra+\lla q^*\dot x,y\rra+\lla rx,y\rra\right)dt,  
\end{align}
for $X,Y\in H^1(S_1;E)$, where $p\in C^1([0,1],\gl(n,\C))$, $q\in C^1([0,1],\gl(n,\C))$, $r\in C([0,1],\gl(n,\C))$, $p(t)=p(t)^*>0$ and $r(t)=r(t)^*$ for all $t\in[0,1]$.
The {\em $k$-th iteration form} $\Ii_k$ of $\Ii$ has the form
\begin{align}\label{e:index-form-iteration}
	\Ii_k(X,Y)=\int_{0}^{1}\left(\lla p\dot x+qx,\dot y\rra+\lla q^*\dot x,y\rra+\lla rx,y\rra\right)dt,  
\end{align}
for $X,Y\in H^1(S_k;E)$, where $p\in C^1([0,k],\gl(n,\C))$, $q\in C^1([0,k],\gl(n,\C))$, $r\in C([0,k],\gl(n,\C))$, $p(t)=p(t)^*>0$ and $r(t)=r(t)^*$ for all $t\in[0,k]$. Moreover, for all $X_0,Y_0\in H_0^1([l,l+1];E)$ and $l=0,\ldots,k-1$, we have 
\begin{align}\label{e:translation}
	\Ii_k(X,Y)=\Ii(X_0(\cdot+l),Y_0(\cdot+l)),
\end{align}
where $X$, $Y$ are $0$ extension of $X_0$ and $Y_0$ respectively.

By \eqref{e:frame-section}, for all $X_0,Y_0\in H_0^1([l,l+1];E)$, $l=0,\ldots,k-1$ and $t\in[0,1]$, there are $x=(x_1,\ldots,x_n),y=(y_1,\ldots,y_n),\tilde x=(\tilde x_1,\ldots,\tilde x_n), \tilde y=(\tilde y_1,\ldots,\tilde y_n)\in H_0^1([0,1],\C^n)$ such that 
\begin{align*}
	X_0(t+l)=&\sum_{i=1}^{n}x_i(t)e_i(t)=\sum_{i=1}^{n}\tilde x_i(t)e_i(t+l),\\ Y_0(t+l)=&\sum_{i=1}^{n}y_i(t)e_i(t)=\sum_{i=1}^{n}\tilde y_i(t)e_i(t+l).
\end{align*}
By \eqref{e:frame-holonomy} and \eqref{e:frame-holonomy-extended}, we have
\begin{align}\label{e:frame-coordinates}
	x(t)=a^l\tilde x(t),\quad y(t)=a^l\tilde y(t).
\end{align}
By \eqref{e:index-form}-\eqref{e:frame-coordinates}, we have (cf. \cite[(91)]{Zhu06})
\begin{align}\label{e:index-form-iterated-matrix}
	\left(\begin{array}{cc}
		p(t) & q (t)\\
		q(t) & r(t)
	\end{array}\right)
	=\left(\begin{array}{cc}
		a^l & 0\\
		0 & a^l
	\end{array}\right)^*\left(\begin{array}{cc}
		p(t+l) & q (t+l)\\
		q(t+l) & r(t+l)
	\end{array}\right)\left(\begin{array}{cc}
	a^l & 0\\
	0 & a^l
	\end{array}\right).
\end{align}

Denote by $J_{2n}\ :=\left(\begin{array}{cc}
	0 & I_n \\
	I_n & 0
\end{array}\right)$. Denote by $\gamma$ the fundamental solution of 
\begin{align}\label{e:linear-hamiltonian-system}
	\dot u=Jbu,
\end{align}
where (cf. \cite[(14.3.6)]{Lo02})
\begin{align}\label{e:first-second}
	b(t)=\left(\begin{array}{cc}
		p(t)^{-1} & -p(t)^{-1}q(t) \\
		-q(t)^*p(t)^{-1} & q(t)^*p(t)^{-1}q(t)-r(t)
	\end{array}\right),\quad \forall t\in\R.
\end{align}

By \cite[Corollary 5.1]{Zhu06}, we have
\begin{align}\label{e:iterated-gamma-second-order}
    \gamma|_{[0,k]}=\gamma|_{[0,1]}(k,P),
\end{align}
where $P=\left(\begin{array}{cc}
	a^{*-1} & 0 \\
	0 & a
\end{array}\right)$.

Then we have the following independence of iteration theory on frames.

\begin{corollary}\label{c:independence-iteration-frame}
	Denote by $\varphi(k,\Ii)\ :=m^-(\Ii_k)-km^-(\Ii)$. Then we have
	\begin{align}
		\label{e:independence-iteration-frame}
		\varphi(k,\Ii)=f(k,P^{-1}\gamma(1),I_{2n})+(k-1)n.
	\end{align}
\end{corollary}

\begin{proof}
	By \cite[Corollary 2.1, Theorem 2.3]{Zhu06} and Theorem \ref{t:iteration-frame}, we have 
	\begin{align*}
		\varphi(k,\Ii)=&
		f(k,\gamma(1),P)-\dim\ker((a^k)^*-I_n)+k\dim\ker(a^*-I_n)\\
		=&f(k,P^{-1}\gamma(1),I_{2n})-f(k,P^{-1},I_{2n})\\
		&-\dim\ker(a^k-I_n)+k\dim\ker(a-I_n)\\
		=&f(k,P^{-1}\gamma(1),I_{2n})+(k-1)n.
	\end{align*}
\end{proof}

\subsection{Mod 2 Maslov-type index}
 
Assume that $(V,\omega)$ be a real symplectic vector space of dimension $2n$, where $n\in\N$. Let $\tau$ be a positive number. In this subsection, we shall give a formula to calculate $i_1(\gamma) \mod 2$ for $\gamma\in\Pp_\tau(V)$, which is shown to be important in the study of the minimal period problem (\cite{DoLo97}). 

We recall 

\begin{definition}\label{d:hyperbolic-index} (cf. \cite[Definition 1.8.1]{Lo02})
	Let $(V,\omega)$ be a real symplectic vector space of dimension $2n$, where $n\in\N$.  Let $M\in\Sp(V)$ be a symplectic linear map. We define the {\em hyperbolic index} $\alpha(M)$ of $M$ by the mod 2 number of the total algebraic multiplicity of negative eigenvalues of $M$ which are strictly less than $-1$. We define $\tilde\alpha(M)$ of $M$ by
	\begin{align}\label{e:tilde-alpha}
		\tilde\alpha(M)\ :=\alpha(M)+\frac{\dim_{\C} E_{-1}(M)}{2}\mod 2,
	\end{align}
	where $E_1(M){-1}$ denotes the complex root space of $M$ belonging to $-1$.
\end{definition}

Let $(V,\omega)$ be a real symplectic vector space of dimension $2n$. For $z\in S^1$ and $M\in\Sp(V)$, we define (cf. \cite[(1.8.3)]{Lo02})
\begin{align}\label{e:D-z-M}
	D_z(M)\ :=(-1)^{n-1}z^{-n}\det(M-zI_{2n}).
\end{align} 

We have the following result.

\begin{theorem}\label{t:Maslov-type-index-mod-2}
	Let $(V,\omega)$ be a real symplectic vector space of dimension $2n$, where $n\in\N$. Let $\tau$ be a positive number. Let $\gamma\in C([0,\tau],\Sp(V))$ be a symplectic path. Then we have
	\begin{align}\label{e:Maslov-type-index-mod-2}
		i_1(\gamma)=\tilde\alpha(\gamma(\tau))+S^+_{\gamma(\tau)}(1)-\tilde\alpha(\gamma(0))-S^+_{\gamma(0)}(1)\mod 2.
	\end{align}
	Especially if $\gamma\in\Pp_\tau(V)$, we have
	\begin{align}\label{e:Maslov-type-index-mod-2-I}
		i_1(\gamma)=\tilde\alpha(\gamma(\tau))+S^+_{\gamma(\tau)}(1)+n\mod 2.
	\end{align}
\end{theorem}

\begin{proof}
	1. Let $\gamma\in\Pp_\tau(V)$ be a symplectic path with $\gamma(\tau)=M$. Let $\varepsilon>0$ be such that 
	\begin{align*}
		\sigma(M)\cap\{z\in\C;\;|z-1|<|e^{\mi\varepsilon}-1|\}\subset\{1\},
	\end{align*}
	where $\sigma(M)$ denotes the spectrum of $M$. 
	Denote by $\sign(a)$ the sign of a real number $a$. Denote by $E_\lambda(M)$ the complex root space of $M$ belonging to $\lambda$. By \cite[Theorem 1.3.1]{Lo02}, we have
	\begin{align*}
		\sign(D_{e^{\mi\varepsilon}}(M)=&-\sign\left((-1)^ne^{-n\mi\varepsilon}\Pi_{\lambda\in\sigma(M)}(\lambda-e^{\mi\varepsilon})^{\dim_{\C} E_\lambda(M)}\right)\\
		=&-\sign\left(\Pi_{\lambda\in\sigma(M)\cap\{1,-1\}}(-e^{-\mi\varepsilon}(\lambda-e^{\mi\varepsilon})^2)^{\frac{\dim_{\C} E_\lambda(M)}{2}}\right)\\
		&\sign\left(\Pi_{\lambda\in\sigma(M)\cap\R,|\lambda|>1}(-e^{-\mi\varepsilon}(\lambda-e^{\mi\varepsilon})(\lambda^{-1}-e^{\mi\varepsilon}))^{\dim_{\C}E_\lambda(M)}\right)\\
		=&-(-1)^{\frac{\dim_{\C} E_{-1}(M)}{2}}\Pi_{\lambda\in\sigma(M)\cap\R,\lambda<-1}(-1)^{\dim_{\C} E_\lambda(M)}\\
		=&-(-1)^{\tilde\alpha(M)}.
	\end{align*}
	By \cite[Theorem 2.4.1]{Lo02} and the definition of Maslov-type index, we have
	\begin{align*}
		i_{e^{\mi\varepsilon}}(\gamma)=&\frac{\sign(D_{e^{\mi\varepsilon}}(M))+1}{2}\\
		=&\frac{-(-1)^{\tilde\alpha(M)}+1}{2}\\
		=&\tilde\alpha(M)\mod 2.
	\end{align*}
	By the definition of the splitting number, we have
	\begin{align*}
		i_1(\gamma)=&i_1(e^{-\mi\varepsilon}\gamma)-i_1(e^{-it}M;\;0\le t\le\varepsilon)+i_1(e^{-it}I_{2n};\;0\le t\le\varepsilon)\\
		=&i_{e^{\mi\varepsilon}}(\gamma)-S^+_M(1)+n\\
		=&\tilde\alpha(M)+S^+_M(1)+n\mod 2.
	\end{align*}
	\newline 2. In the general case, by \cite[Corollary 2.2.8]{Lo02}, there is a path $\gamma_1\in\Pp_\tau(V)$ with $\gamma_1(\tau)=\gamma(0)$. By \cite[Proposition 2.3.1.b]{BoZh18}, we have
	\begin{align*}
		i_1(\gamma)=&i_1(\gamma_1*\gamma)-i_1(\gamma_1)\\
		=&\tilde\alpha(\gamma(\tau))+S^+_{\gamma(\tau)}(1)-\tilde\alpha(\gamma(0))-S^+_{\gamma(0)}(1)\mod 2.
	\end{align*}
\end{proof}

By \cite[Theorem 1.6]{ZhuGs21}, the following corollary generalize \cite[Theorem 1.1]{HuSu09} and \cite[Theorem 1.3]{ZhuGs21}. 

\begin{corollary}\label{c:mod-2-Morse-index}
	Let $\pi\colon E\to S_1$ be a $C^2$ real vector bundle of dimensional $n$. Let $\Ii$ be an $C^1$ index form of order $1$ on the Sobolev space $H^1(S_1;E)$ with positive definite highest order term. With the notations in Corollary \ref{c:independence-iteration-frame}, we have 
	\begin{align}\label{e:mod-2-Morse-index}
	    m^-(\Ii)=\tilde\alpha(P^{-1}\gamma(1))+S^+_{P^{-1}\gamma(1)}(1)+\frac{\sign(\det a)-1}{2}\mod 2.
	\end{align}
\end{corollary}

\begin{proof}
	By \cite[Corollary 2.1]{Zhu06}, Theorem \ref{t:Maslov-type-index-mod-2} and Theorem \ref{t:splitting-numbers-up-triangular}, we have
	\begin{align*}
		m^-(\Ii)=&i_P(\gamma|_{[0,1]})-\dim\ker(a^*-I)\\
		=&\tilde\alpha(P^{-1}\gamma(1))+S^+_{P^{-1}\gamma(1)}(1)-\tilde\alpha(P^{-1})-S^+_{P^{-1}}(1)
		-\dim\ker(a-I)\\
		=&\tilde\alpha(P^{-1}\gamma(1))+S^+_{P^{-1}\gamma(1)}(1)+\tilde\alpha(P^{-1})
		\\
		=&\tilde\alpha(P^{-1}\gamma(1))+S^+_{P^{-1}\gamma(1)}(1)+\frac{\sign(\det a)-1}{2}\mod 2.
	\end{align*}
\end{proof}

%% file: operator-space.tex
\appendix
\section{Continuous families of bounded linear relations}\label{s:operator-space}

In this appendix we study the continuity of families of bounded linear relations and families of bounded linear operators acting on closed linear subspaces. Then we get a generalization of \cite[Lemma 0.1]{Ne68}.

Denote by $\Ss(X)$ ($\Ss^c(X)$) the
set of all (complemented) closed linear subspaces of a Banach space $X$. Denote by $\Bb(X,Y)$
($\Cc(X,Y)$, $\CLR(X,Y)$) the set of all bounded operators (closed operators, closed linear relations) between
Banach spaces $X$ and $Y$. We equip $\Ss(X)$ and $\Ss^c(X)$ with the gap topology.

Firstly we recall the notion of the gap between closed subspaces in a given Banach space $X$. We denote by $\dist(u,K)$ the distance between $u\in X$ and a subset of $K$ of $X$.

For two subsets $A$, $B$ of $X$, we denote by 

\begin{align}\label{e:distence-banach}
	\dist(A,B)=\inf\limits_{u\in A,v\in B}\|u-v\|.
\end{align}

\begin{definition}[The gap between subspaces]\label{d:closed-distance}	
	(a) We set \begin{multline*}
		d(M,N)\ =\
		d(\romS_M,\romS_N)\\
		:= \begin{cases} \max\left\{\begin{matrix}\sup\limits_{u\in \romS_M}\dist(u,\romS_N),\\ \sup\limits_{u\in
					\romS_N}\dist(u,\romS_M)\end{matrix}\right\},&\text{ if both $M\ne 0$ and $N\ne 0$},\\
			0,&\text{ if $M=N=0$},\\
			2,&\text{ if either $M= 0$ and $N\ne 0$ or vice versa}.
		\end{cases}
	\end{multline*}
	\newline (b) If $M\cap N$ is closed, we set
	\begin{eqnarray*}
		\delta(M,N)\ &:=&\ \begin{cases}\sup\limits_{u\in \romS_M}\dist(u,N),& \text{if $M\ne\{0\}$},\\
			0,& \text{if $M=\{0\}$},\end{cases}\\
		\hat\delta(M,N)\ &:=&\ \max\{\delta(M,N),\delta(N,M)\}.
	\end{eqnarray*}
	$\hat\delta(M,N)$ is called the {\em gap} between $M$ and $N$.
	\newline (c) We set
	\begin{eqnarray*}
		\gamma(M,N)\ :&=&\ \begin{cases} \inf\limits_{u\in M\setminus N}\frac{\dist(u,N)}{\dist(u,M\cap N)}\ (\le
			1), & \text{if $M\nsubseteq N$},\\
			1, & \text{if $M\subset N$},\end{cases}\\
		\hat\gamma(M,N)\ :&=&\ \min\{\gamma(M,N),\gamma(N,M)\}.
	\end{eqnarray*}
	$\hat\gamma(M,N)$ is called the {\em minimum gap} between $M$ and $N$. If $M\cap N=\{0\}$, we have
	\[\gamma(M,N)\ =\ \inf_{u\in \romS_M}\dist(u,N).\]
\end{definition}

We have the following lemma (cf. \cite[Theorem IV.5.2]{Ka95}).

\begin{lemma}\label{l:closed-domain}
	Let $Z=X\oplus Y$ be a Banach space with closed linear subspaces $X$, $Y$. We view a linear subspace $M$ as a linear relation between $X$ and $Y$. Then the following hold.
	\begin{itemize}
		\item[(a)] Assume that $M$ is closed. Then the domain $\dom(M)$ is closed if and only if $M+Y$ is closed.
		\item[(b)] Assume that $M$ is closed and $\dom(M)$ is a nonzero closed linear subspace. Then we have $\gamma(M,Y)>0$, $\gamma(X,Y)>0$, and 
		\begin{align*}
			\max\{\gamma(X,Y)\|x\|,\gamma(M,Y)\dist(x,Mx )\}\le\dist(x,Y)\le\delta(X,Y)\|x\|
		\end{align*}
		for each $x\in\dom(M)$.
	\end{itemize}
\end{lemma}

\begin{proof}
	(a) By \cite[Lemma A.1.1]{BoZh18}, we have $M+Y=\dom(M)\oplus  Y$. Then the domain $\dom(M)$ is closed if and only if $M+Y$ is closed.
	\newline (b) By (a), $M+Y$ is closed . By \cite[Theorem IV.4.2]{Ka95}, we have $\gamma(M,Y)>0$ and $\gamma(X,Y)>0$. Since there hold
	\begin{align*}
		\dist(x,Y)&=\dist(x+y, Y),\\ 
		\dist(x,Mx)&=\dist(x+y, M\cap Y)
	\end{align*}
	for each $x\in\dom(M)$, $y\in Mx$, our inequality follows from Definition \ref{d:closed-distance}.
\end{proof}

We shall define $a(M)$ and $\|M\|$ for a linear relation between two Banach spaces $X$ and $Y$. 

\begin{definition}\label{d:alpha-beta-lr}(cf. \cite[Definitions II.1.3]{Cr98})
	Let $Z=X\oplus Y$ be a Banach space with closed linear subspaces $X$, $Y$. We view a linear subspace $M$ as a linear relatioin between $X$ and $Y$. We define $a(M)$ and $\|M\|$ by
		\begin{align}\label{e:alpha-norm-relation}
			a(M):&=\begin{cases}
				\sup\limits_{x\in\dom(M)\setminus\{0\}}\frac{\dist(x,Mx )}{\|x\|}& \text{ if $\dom(M)\ne\{0\}$},\\
				0& \text{ if $\dom(M)=\{0\}$},
			\end{cases}\\
			\label{e:beta-norm-relation}
			\|M\|:&=\begin{cases}
				\sup\limits_{x\in\dom(M)\setminus\{0\}}\inf\limits_{y\in Mx}\frac{\|y\|}{\|x\|}&\text{ if $\dom(M)\ne\{0\}$},\\
				0& \text{ if $\dom(M)=\{0\}$}
			\end{cases}
		\end{align}
		respectively. We call $M$ {\em bounded} if $\dom(M)$ is closed and $\|M\|<+\infty$.
\end{definition}

We need the notion of distance between  affine spaces. 

\begin{definition}\label{d:distance-affine-space}
	Let $X$ be a Banach space with two closed affine subspaces $M$ and $N$. Let $M_0$ and $N_0$ be the associated linear spaces of $A$ and $B$ respectively. We define $d(M,N)$  by
	\begin{align}\label{e:distance-affine-space}
		d(M,N)\ :=d(M_0,N_0)+\inf_{u\in M,v\in N}\|u-v\|.
	\end{align}  
	Then $d(M,N)$ define a metric on the set of closed affine subspaces.
\end{definition}

We have the following criteria of bounded linear relations.

\begin{lemma}\label{l:criteria-bounded-lr}
	Let $Z=X\oplus Y$ be a Banach space with closed linear subspaces $X$, $Y$. We view a linear subspace $M$ as a linear relation between $X$ and $Y$. Assume that $\dom(M)$ is closed. Then the following three properties are equivalent:
	\begin{itemize}		
		\item[(i)] $M$ is bounded;		 
		\item[(ii)] $M0$ is closed and the family $\{Mx;\;x\in X\}$ is continuous;
		\item[(iii)] $M0$ is closed and the family $\{Mx;\;x\in X\}$ is continuous at $x=0$;
		\item[(iv)] $M$ is closed.
	\end{itemize}
\end{lemma}

\begin{proof} 
   1. (i)$\Rightarrow$(ii)$\Rightarrow$(iii). By definition.
	\newline 2. (iii)$\Rightarrow$(iv). Let $(x_n,y_n)\in M$ be a sequence convergent to $(x,y)\in Z$. Since $\dom(M)$ is closed, we have $x\in\dom(M)$. Since $x_n-x\to 0$ and $x_n-x\in\dom(M)$, we have 
	\begin{align*}
		d(M(x_n-x),M0)=\inf_{u\in M(x_n-x),v\in M0}\|u-v\|=\dist(y_n,Mx)\to 0.
	\end{align*}
	Then we have $\dist(y,Mx)=0$. Since $M0$ is closed, $Mx$ is closed. Thus we obtain $y\in Mx$ and $M$ is closed.
	\newline 3. (iv)$\Rightarrow$(i). Since $M$ is closed, $M0=M\cap Y$ is closed. Then 
	\begin{align}
		M/M0\colon\dom(M)\to Y/M0
	\end{align} 
	is a closed linear operator. By the closed graph theorem, we have $M/M0\in\Bb(\dom(M), Y/M0)$. Then $M$ is bounded.  	
\end{proof}

We have the following lemma (cf. \cite[Proposition II.3.2.a, (II.5.6)]{Cr98}) which shows that a bounded linear relation is always closed.

\begin{lemma}\label{l:bounded-lr}
	Let $Z=X\oplus Y$ be a Banach space with closed linear subspaces $X$, $Y$. We view a linear subspace $M$ as a linear relation between $X$ and $Y$. Then the following hold.
	\begin{itemize}		
		\item[(a)] We have
	    \begin{align}\label{e:lower-bound-alpha}
	    	a(M)
	    	\begin{cases}
	    		\ge\delta(\dom(M),Y)\ge\gamma(X,Y)&\text{ if $\dom(M)\ne\{0\}$},\\
	    		=0& \text{ if $\dom(M)=\{0\}$}
	    	\end{cases}.
	    \end{align}
	    and
	    \begin{align}\label{e:alpha-beta-norm-relation}
	    	|a(M)-\|M\| |\le 1.
	    \end{align}
        \item[(b)] Assume that $M\cap Y$ is closed and $\dom(M)$ is a nonzero closed linear subspace. Assume that $a(M)$ (or $\|M\|$) is finite. 
		Then $M$ is closed, and we have
		\begin{align}\label{e:gamma-norm-relation}
			\gamma(M,Y)\ge\frac{\gamma(X, Y)}{a(M)}.
		\end{align}  
	\end{itemize}
\end{lemma}

\begin{proof}
	 (a) If $\dom(M)=\{0\}$, we have $a(M)=\delta(\dom(M),Y)=\|M\|=0$ and our results hold.
	
	Now we assume that $\dom(M)\ne\{0\}$. Since $\dist(x, Mx)\ge\dist(x,Y)$ if $x\in\dom(M)$ and $y\in Y$, we have
	\begin{align*}
		a(M)\ge\delta(\dom(M),Y)\ge\gamma(X,Y).
	\end{align*}
	
	We have
	\begin{align*}		
		\dist(x,Mx)&=\inf\{\|x-y\|;\;y\in Mx\}\\
		&\in\left[\inf\{\|y\|;\;y\in Mx\}-\|x\|,\inf\{\|y\|;\;y\in Mx\}+\|x\|\right].		
	\end{align*}
	for each $x\in\dom(M)$. By definition we have $\|M\|-1\le a(M)\le\|M\|+1$.
	\newline (b) 	
	Let $\{x_n+y_n\}_{n\in\N}\subset M$ be a sequence of $M$ convergent in $Z$, where $x_n\in X$ and $y_n\in Y$. Then there is an $x\in X$ and a $y\in Y$ such that 
	\begin{align*}
		\lim_{n\to+\infty} x_n=x,\quad \lim_{n\to+\infty} y_n=y.
	\end{align*}
	Since $x_n\in\dom(M)$, $y_n\in Mx_n$ and $\dom(M)$ is closed, we have $x\in\dom(M)$. Then we have
	\begin{align*}
		\dist(y,Mx)&\le\|y-y_n\|+\dist(y_n,Mx)\\
		&\le\|y-y_n\|+\|M\|\|x_n-x\|.
	\end{align*}
    Let $n\to+\infty$, we obtain $\dist(y,Mx)=0$ and $y\in Mx$. Thus $M$ is closed.
    
    By Definition \ref{d:closed-distance}, we have
    \begin{align*}
    	\gamma(M,Y)&=\inf_{
    		\begin{subarray}{l}
    			x+y\in M\\x\in X\setminus\{0\},y\in Y
    		\end{subarray}}\frac{\dist(x+y, Y)}{\dist(x+y,M\cap Y)}\\ 
    	&=\inf_{x\in\dom(M)\setminus\{0\}}\frac{\dist(x, Y)}{\dist(x,Mx)}\\ 
        &\ge\inf_{x\in\dom(M)\setminus\{0\}}
        \frac{\gamma(X, Y)\|x\|}{a(M )\|x\|}\\
        &=\frac{\gamma(X, Y)}{a(M)}.		  		
    \end{align*}   
\end{proof}

The following theorem shows that the functions $a(M)$ and $\|M\|$ of $M$ are continuous functions for bounded linear relations (cf. \cite[(1) of Lemma 0.1 ]{Ne68} for operator case).

\begin{theorem}\label{t:closed-domain}
	Let $Z=X\oplus Y$ be a Banach space with closed linear subspaces $X$, $Y$. We view two linear subspaces $M$ and $N$ as linear relations between $X$ and $Y$.  
	Denote by $P\in\Bb(Z)$ the projection  on $X$ along $Y$. Denote by 
	\begin{align*}
		\eta(X,Y):=\|P\|+\|I-P\|\ge 1.
	\end{align*}
	Then we have the following estimates
	\begin{align}
		\label{e:continuity-alpha}
		a(N)&\ge\frac{\frac{1-\delta(N0,M0)}{1+\delta(N0,M0)}a(M)-\eta(X,Y)(2+a(M))\delta(M,N)}{1+\|P\|(2+a(M))\delta(M,N)},\\
		\label{e:continuity-beta}
		\|N\|&\ge\frac{\frac{1-\delta(N0,M0)}{1+\delta(N0,M0)}\|M\|-\|I-P\|(1+\|M\|)\delta(M,N)}{1+\|P\|(1+\|M\|)\delta(M,N)}.
	\end{align} 
\end{theorem}

\begin{proof}	
	1. The case $\dom(M)=\{0\}$. 
	
	In this case, we have $a(M)=\|M\|=0$ and our estimates hold. 
	\newline 2.
	The estimate of $a(N)$ when $\dom(M)\ne\{0\}$. 
	
	For each $\varepsilon>0$, there exist an $x_1\in\dom(M)\setminus\{0\}$ and a $y_1\in Mx_1$ such that 
	\begin{align*}
		\|x_1-y_1\|&\le(1+\varepsilon)a(M)\|x_1\|,\\
		\dist(x_1,Mx_1)&\ge(1-\varepsilon)a(M)\|x_1\|.   
	\end{align*}
	Then we have 
	\begin{align*}
		\|x_1+y_1\|\le (2+(1+\varepsilon)a(M))\|x_1\|.
	\end{align*}
	By Definition \ref{d:closed-distance}, there exist an $x_2\in\dom(N)$ and a $y_2\in Nx_2$ such that
	\begin{align*}
		d\ :&=\|x_1+y_1-x_2-y_2\|\le(1+\varepsilon)\delta(M,N)\|x_1+y_1\|\\
		&\le (2+(1+\varepsilon)a(M))(1+\varepsilon)\delta(M,N)\|x_1\|.  	
	\end{align*}    
	By definition we have
	\begin{align*}   	
		\|x_1-x_2\|&\le\|P\|d,\\
		\|y_1-y_2\|&\le\|I-P\|d,\\
		\|x_2\|&\le \|P\|d+\|x_1\|.
	\end{align*}
	Then we have
	\begin{align*}
		\|x_2-y_2&|\le\|x_1-y_1\|+\|x_1-x_2\|+\|y_1-y_2\|\\
		&\le (1+\varepsilon)a(M)\|x_1\|+\eta(X,Y)d.
	\end{align*}
	By Definition \ref{d:closed-distance}, for each $y_{20}\in N0$, there exists a $y_{10}\in M0$ such that 
	\begin{align*}
		\|y_{20}-y_{10}\|\le(1+\varepsilon)\delta(N0,M0)\|y_{20}\|.  	
	\end{align*} 
	Then we have
	\begin{align*}
		\|&x_2-y_2-y_{20}\|\\
		&\ge \|x_1-y_1-y_{10}\|-\|x_1-x_2\|-\|y_1-y_2\|-\|y_{20}-y_{10}\|\\
		&\ge \dist(x_1,Mx_1)-\|x_1-x_2\|-\|y_1-y_2\|-\|y_{20}-y_{10}\|\\
		&\ge(1-\varepsilon) a(M)\|x_1\|-\eta(X,Y)d-(1+\varepsilon)\delta(N0,M0)\|y_{20}\|,
	\end{align*}
	and 
	\begin{align*}
		\|x_2-y_2-y_{20}\|&
		\ge\|y_{20}\|-\|x_2-y_2\|\\
		&\ge \|y_{20}\|-(1+\varepsilon)a(M)\|x_1\|-\eta(X,Y)d.
	\end{align*}
	The right hand sides of the above two equalities are equal when
	\begin{align*}
		\|y_{20}\|=2(1+(1+\varepsilon)\delta(N0,M0))^{-1}a(M)\|x_1\|\ge 0.
	\end{align*} 
	We obtain the estimate free of $\|y_{20}\|$ and the right hand side takes the common value of these right hand sides:
	\begin{align*}
		\|x_2-y_2-y_{20}\|
		\ge &2(1+(1+\varepsilon)\delta(N0,M0))^{-1}a(M)\|x_1\|-\\
		&
		\left( (1+\varepsilon)a(M)\|x_1\|+\eta(X,Y)d\right).
	\end{align*}
	By definition
	and let $\varepsilon\to 0+$, we obtain \eqref{e:continuity-alpha}.
	\newline 3. The estimate of $\|N\|$ when $\dom(M)\ne\{0\}$. 
	
	For each $\varepsilon>0$, there exist an $x_1\in\dom(M)\setminus\{0\}$ and a $y_1\in Mx_1$ such that 
	\begin{align*}
		\|y_1\|&\le(1+\varepsilon)\|M\|\|x_1\|,\\
		\inf\{\|y\|;y\in Mx_1\}&\ge(1-\varepsilon)\|M\|\|x_1\|.   
	\end{align*}
	Then we have 
	\begin{align*}
		\|x_1+y_1\|\le (1+(1+\varepsilon)\|M\|)\|x_1\|.
	\end{align*}
	By Definition \ref{d:closed-distance}, there exist an $x_2\in\dom(N)$ and a $y_2\in Nx_2$ such that
	\begin{align*}
		d:&=\|x_1+y_1-x_2-y_2\|\le(1+\varepsilon)\delta(M,N)\|x_1+y_1\|\\
		&\le (1+(1+\varepsilon)\|M\|)(1+\varepsilon)\delta(M,N)\|x_1\|  	
	\end{align*}    
	By definition we have
	\begin{align*}   	
		\|x_1-x_2\|&\le\|P\|d,\\
		\|y_1-y_2\|&\le\|I-P\|d,\\
		\|x_2\|&\le \|P\|d+\|x_1\|.
	\end{align*}
	Then we have
	\begin{align*}
		\|y_2\|&\le\|y_1\|+\|y_1-y_2\|\\
		&\le (1+\varepsilon)\|M\|\|x_1\|+\|I-P\|d.
	\end{align*}
	By Definition \ref{d:closed-distance}, for each $y_{20}\in N0$ there exists a $y_{10}\in M0$ such that 
	\begin{align*}
		\|y_{20}-y_{10}\|\le(1+\varepsilon)\delta(N0,M0)\|y_{20}\|.  	
	\end{align*} 
	Then we have
	\begin{align*}
		\|&y_2+y_{20}\|\\
		&\ge \|y_1+y_{10}\|-\|y_1-y_2\|-\|y_{20}-y_{10}\|\\
		&\ge \inf\{\|y\|;y\in Mx_1\}-\|y_1-y_2\|-\|y_{20}-y_{10}\|\\
		&\ge(1-\varepsilon)\|M\|\|x_1\|-\|I-P\|d-(1+\varepsilon)\delta(N0,M0)\|y_{20}\|,
	\end{align*}
	and 
	\begin{align*}
		\|y_2+y_{20}\|&
		\ge\|y_{20}\|-\|y_2\|\\
		&\ge \|y_{20}\|-(1+\varepsilon)\|M\|\|x_1\|-\|I-P\|d.
	\end{align*}
	The right hand side of the above two equalities are equal when
	\begin{align*}
		\|y_{20}\|=2(1+(1+\varepsilon)\delta(N0,M0))^{-1}\|M\|\|x_1\|\ge 0.
	\end{align*} 
	We obtain the estimate free of $\|y_{20}\|$ and the right hand side takes the common value of the right hand side:
	\begin{align*}
		\|y_2+y_{20}\|
		\ge &2(1+(1+\varepsilon)\delta(N0,M0))^{-1}\|M\|\|x_1\|-\\
		&
		\left( (1+\varepsilon)\|M\|\|x_1\|+\|I-P\|d\right).
	\end{align*}
	By definition
	and let $\varepsilon\to 0+$, we obtain \eqref{e:continuity-beta}.
\end{proof}

We need the following notion for the study of the continuity of the bounded linear relations.

\begin{definition}\label{d:b-lr}
	Let $Z=X\oplus Y$ be a Banach space with closed linear subspaces $X$, $Y$. We view two linear subspaces $M$ and $N$ as linear relations between $X$ and $Y$. For $s>\delta(\dom(M),\dom(N))$ and $t>\|M\|$, we define
	\begin{align}\label{e:b-lr}
		\begin{aligned}
		&b(s,t,M,N)\ :=\\
		&\begin{cases}
			\sup\limits_{\substack{x_1\in\dom(M)\setminus\{0\}\\
					x_2\in \dom(N)\\ \|x_2-x_1\|\le s\|x_1\|
				}}
			\inf\limits_{\substack{y_1\in Mx_1\\ \|y_1\|\le t\|x_1\|}}\frac{\dist(y_1,Nx_2)}{\|x_1\|}&\text{ if $\dom(M)\ne\{0\}$},\\
			0& \text{ if $\dom(M)=\{0\}$}.
		\end{cases}
	\end{aligned}
	\end{align}
\end{definition}

The following theorem gives a criterion for the continuity of families of bounded linear relations with continuously varying domains and indeterminacies (cf. \cite[(2) of Lemma 0.1 ]{Ne68} for operator case).

\begin{theorem}\label{t:continuous-bounded}
	Let $Z=X\oplus Y$ be a Banach space with closed linear subspaces $X$, $Y$. We view two linear subspaces $M$ and $N$ as linear relations between $X$ and $Y$. Denote by $P\in\Bb(Z)$ the projection on $X$ along $Y$. Denote by 
	\begin{align*}
		C(s,t,M,N)\;:=s\|P\|+b(s,t,M,N)\|P\|+\delta(M0,N0)(\|I-P\|+t\|P\|).
	\end{align*} 
    Then for $s>\delta(\dom(M),\dom(N))$ and $t>\|M\|$, 
	we have
	\begin{align}\label{e:continuous-bounded-1}
		\delta(M,N)\le& \max\{C(s,t,M,N),\delta(M0,N0)\},\\
		\label{e:continuous-bounded-2}
		b(s,t,M,N)\le&(t+1)(\|I-P\|+\|P\|\|N\|)\delta(M,N)+s\|N\|.
	\end{align}
\end{theorem}

\begin{proof} 1. The estimate of $\delta(M,N)$.
	
	If $M=\{0\}$, we have $\delta(M,N)=0$ and \eqref{e:continuous-bounded-1} holds. Now assume that $M\ne\{0\}$.
	Let $z_1\;:=x_1+y_1\ne 0$ be in $M$ with $x_1\in X$ and $y_1\in Y$. Then we have 
	\begin{align*}
		\|x_1\|\le\|P\|\|z_1\|,\quad \|y_1\|\le\|I-P\|\|z_1\|.
	\end{align*}
	
	If $x_1=0$, we have $y_1\in M0$ and 
	\begin{align*}
		\dist(z_1,N)=\dist(y_1,N)\le\dist(y_1,N0)\le\delta(M0,N0)\|y_1\|. 
	\end{align*}

	Now we assume that $x_1\ne 0$. Since $t>\|M\|$, the set $A\ :=\{y\in Mx_1;\;\|y\|\le t\|x_1\|\}$ is nonempty. Since $s>\delta(\dom(M),\dom(N))$, there exists an $x_2\in\dom (N)$ such that 
	\begin{align*}
		\|x_2-x_1\|\le s\|x_1\|\le s\|P\|\|z_1\|.
	\end{align*}
	For such an $x_2$ we have 
	\begin{align*}
		\inf_{y_3\in A}\dist(y_3,Nx_2)\le b(s,t,M,N)\|x_1\|.
	\end{align*}
	Thus for each $\varepsilon>0$, there exist a $y_3\in A$ and a $y_2\in Nx_2$ such that
	\begin{align*}
		\|y_3-y_2\|\le (b(s,t,M,N)+\varepsilon)\|x_1\|\le (b(s,t,M,N)+\varepsilon)\|P\|\|z_1\|.
	\end{align*}
	Set $y_{10}\ :=y_1-y_3$. Then we have $y_{10}\in M0$ and 
	\begin{align*}
		\|y_{10}\|\le\|y_1\|+\|y_3\|\le\|y_1\|+t\|x_1\|\le (\|I-P\|+t\|P\|)\|z_1\|.
	\end{align*}
	Then there exists a $y_{20}\in N0$ such that
	\begin{align*}
		\|y_{10}-y_{20}\|\le (\delta(M0,N0)+\varepsilon)\|y_{10}\|\le (\delta(M0,N0)+\varepsilon)(\|I-P\|+t\|P\|)\|z_1\|.
	\end{align*}
    Denote by $z_2\;:=x_2+y_2+y_{20}\in N$,    
    Then we have
    \begin{align*}
    	\dist(z_1,N)&\le\|z_1-z_2\|\\
    	&\le\|x_1-x_2\|+\|y_3-y_2\|+\|y_{10}-y_{20}\|\\
    	&\le\left(C(s,t,M,N)+\varepsilon((1+t)\|P\|+\|I-P\|)\right)\|z_1\|. 
    \end{align*}
    On letting $\varepsilon\to 0$, we have $\dist(z_1,N)\le C(s,t,M,N)\|z_1\|$.         
    By the definition of the gap, \eqref{e:continuous-bounded-1} follows.
    \newline 2. The estimate of $b(s,t,M,N)$ for $s>\delta(\dom(M),\dom(N))$ and $t>\|M\|$.
    
    Given an $x_1\in\dom(M)\setminus\{0\}$, there exists a $y_1\in Y$ such that there hold $z_1\;:=x_1+y_1\in M$ and 
    $\|y_1\|\le t\|x_1\|$. Then we have $\|z_1\|\le (t+1)\|x_1\|$.
    For each $\varepsilon>0$, there exists a $z_4\;:=x_4+y_4\in N$ with $x_4\in X$ and $y_4\in Y$ such that 
    \begin{align*}
    	\|z_4-z_1\|\le(\delta(M,N)+\varepsilon)\|z_1\|\le
    	(\delta(M,N)+\varepsilon)(t+1)\|x_1\|.
    \end{align*}
    Then we have 
    \begin{align*}
    	\|x_4-x_1\|&\le\|P\|\|z_1\|\le\|P\|(\delta(M,N)+\varepsilon)(t+1)\|x_1\|,\\ \|y_4-y_1\|&\le\|I-P\|\|z_1\|\le\|I-P\|(\delta(M,N)+\varepsilon)(t+1)\|x_1\|.
    \end{align*}    
        
    Let an $x_2\in\dom (N)$ be given such that 
    \begin{align*}
    	\|x_2-x_1\|\le s\|x_1\|.
    \end{align*}
    
    For such a $x_2$ we have 
    \begin{align*}
    	\|x_2-x_4\|&\le \|x_2-x_1\|+\|x_4-x_1\|\\
    	&\le \left(s+\|P\|(\delta(M,N)+\varepsilon)(t+1)\right)\|x_1\|.
    \end{align*}
    Since $\varepsilon>0$, there exist a $y_2\in Nx_2$  such that
    \begin{align*}
    	\|y_2-y_4\|&\le (\|N\|+\varepsilon)\|x_2-x_4\|.
    \end{align*}    
    
    By the above estimates we obtain that 
    \begin{align*}
    	\dist&(y_1,Nx_2)\le\|y_1-y_2\|\\
    	\le&\|y_1-y_4\|+\|y_2-y_4\|\\
    	\le&\|I-P\|(\delta(M,N)+\varepsilon)(t+1)\|x_1\|+(\|N\|+\varepsilon)\\
    	&\left(s+\|P\|(\delta(M,N)+\varepsilon)(t+1)\right)\|x_1\|. 
    \end{align*}
    On letting $\varepsilon\to 0$, we have
    \begin{align*}
    	\dist(y_1,Nx_2)
    	\le\left((t+1)(\|I-P\|+\|P\|\|N\|)\delta(M,N)+s\|N\|\right)\|x_1\|. 
    \end{align*}        
    By the definition of $b(s,t, M, N)$, \eqref{e:continuous-bounded-2} follows.
\end{proof}

We have the following estimate which generalize \cite[Theorem IV.2.29]{Ka95}.

\begin{lemma}\label{l:operator-space}
	Let $X$, $Y$ be two Banach spaces. Let $A\in\Bb(X,Y)$ be a bounded operator with bounded $C=A^{-1}|_{AM}:AM\to M$, $B\in\Bb(X,Y)$ be a bounded operator. Let $M$, $N$ be two closed linear subspace of $X$.
	Set
	\begin{align*}
		\kappa\ :=\|C\|^{-1}(1-\delta(N,M))-\|A\|
		\delta(N,M)-\|A-B\|.
	\end{align*}	
	Then the following hold. 
	\begin{itemize}
		\item[(a)] The space $AM$ is closed, and we have
		\begin{align}\label{e:operator-space}
			\delta(AM,BN)\le \|C\|(\|A-B\|+\|B\|\delta(M,N)).
		\end{align}
		\item[(b)] Assume that $\kappa>0$. Then the operator $D=B^{-1}|_{BN}:BN\to N$ is bounded with $\|D\|\le \kappa^{-1}$, $BN$ is closed and $\delta(N,M)\le (\|A\|\|C\|+1)^{-1}\le\frac{1}{2}$.
	\end{itemize}	
\end{lemma}

\begin{proof} 1.
	Since $A\in\Bb(X,Y)$ is a bounded operator with bounded    $C=A^{-1}|_{AM}:AM\to M$ and $M$ is closed, $AM$ is closed.	
	
	If $M=\{0\}$, we have $\delta(M,N)=0$ and \eqref{e:operator-space} holds.
	
	Now we assume $M\ne\{0\}$. For each $x\in M\setminus\{0\}$ we have      
	\begin{align*}
		\dist(x,N)\le \delta(M,N)\|x\|.
	\end{align*}
	Then we have 
	\begin{align*}
		\dist(Ax,BN)&\le \|Ax-Bx\|+\dist(Bx,BN)\\
		&\le \|Ax-Bx\|+\|B\|\dist(x,N)\\
		&\le \|A-B\|\|x\|+\|B\|\|x\|\delta(M,N).
	\end{align*}
	Thus we obtain
	\begin{align*}
		\delta(AM,BN)&=\sup_{x\in M\setminus\{0\}}\frac{\dist(Ax,BN)}{\|Ax\|} \\   	
		&\le\|C\|(\|A-B\|+\|B\|\delta(M,N)).
	\end{align*}
	\newline (b) Since $\kappa>0$, we have $C\ne 0$ and $M\ne \{0\}$. Then we have $\|A\|\|C\|\ge 1$ and $\delta(N,M)\le (\|A\|\|C\|+1)^{-1}\le\frac{1}{2}$. 
	
	If $N=\{0\}$, we have $D=0$ and $BN=\{0\}$.
	
	Now we assume that $N\ne\{0\}$.  For each $x_2\in N\setminus\{0\}$, we have
	\begin{align*}
		\dist(x_2,M)\le\delta(N,M)\|x_2\|.
	\end{align*}
	Then for each $\varepsilon>0$, there exists an $x_1\in M$ such that 
	\begin{align*}
		\|x_2-x_1\|\le(\delta(N,M)+\varepsilon)\|x_2\|.
	\end{align*}
	Then we have 
	\begin{align*}
		(1-\delta(N,M)-\epsilon)\|x_2\|\le\|x_1\|\le(1+\delta(N,M)+\varepsilon)\|x_2\|.
	\end{align*}
	Consequently we have 
	\begin{align*}
		\|Bx_2\|&\ge\|Ax_1\|-\|Ax_2-Ax_1\|-\|Ax_2-Bx_2\|\\
		&\ge \|C\|^{-1}\|x_1\|-\|A\|\|x_2-x_1\|-\|A-B\|\|x_2\|\\ &\ge \left(\|C\|^{-1}(1-\delta(N,M)-\varepsilon)-\|A\|(\delta(N,M)+\varepsilon)-\|A-B\|\right)\|x_2\|.
	\end{align*}
	On letting $\varepsilon\to 0+$, we have $\|Bx_2\|\ge\kappa\|x_2\|$. Then we obtain $\|D\|\le\kappa^{-1}$. By (a), the space $BN$ is closed.
\end{proof}

We want to study the case when $A$, $B$ may have nontrivial kernels. We have the following theorem (cf. \cite[(3) of Lemma 0.1 ]{Ne68} for operator case). 

\begin{theorem}\label{t:operator-space}
	Let $B$ be a topological space. Let $X$, $Y$ be a Banach spaces with direct sum decomposition into continuous families of closed linear subspaces 
	\begin{align}\label{e:XY-direct-decomposition}
		\begin{cases}		
			X=X_1(b)\oplus X_0(b),  &\\
			Y=Y_1(b)\oplus Y_0(b)	&
		\end{cases}
	\end{align}
	for $b\in B$
	respectively. Let $\{M(b)\}_{b\in B}$ be a continuous family of closed linear subspace of $X$. For each $b\in B$, we view $M(b)$ as a closed linear relation between $X_1(b)$ and $X_0(b)$. Assume that for each $b\in B$, the space $M(b)+X_0(b)$ is closed. Let $\{A(b)\in\Bb(X,Y)\}_{b\in B}$ be a continuous family of linear operators. Assume that for each $b\in B$, the operator $A(b)$ has has the form
	\begin{align}\label{e:A-form-XY}
		A(b)=\left(\begin{array}{cc}
			A_{11}(b)& 0  \\
			A_{01}(b)& A_{00}(b)
		\end{array}\right)
	\end{align}
	under direct sum decompositions \eqref{e:XY-direct-decomposition} with $A_{11}(b)^{-1}\in\Bb(Y_1(b),X_1(b))$.
	Assume that there is a continuous family $\{N_0(b)\in\Ss(Y_0(b))\}_{b\in B}$ with
	\begin{align}\label{e:N0b}		
		N_0(b)\supset A(b)M(b)\cap Y_0(b)
	\end{align}
	for each $b\in B$.
	For each $b\in B$, we set 
	\begin{align}\label{e:Nb}
		N(b):=A(b)M(b)+N_0(b).  
	\end{align}
	Then the family $\{N(b)\}_{b\in B}$ is a continuous family of closed linear subspaces of $Y$ with continuous varying closed domains and indeterminacies.
\end{theorem}

\begin{proof}
	1. Since the statement is local, by \cite[Lemma I.4.10]{Ka95} we can assume that $X_j(b)=X_j$ and $Y_j(b)=Y_j$ for $j=1,0$. 
	\newline 2. Boundedness of $M(b)$ and continuity of families 
	\begin{align}\label{e:family-Nb0}
		\{\dom (N(b))\}_{b\in B},\quad \{N(b)0\}_{b\in B}.
	\end{align}
	
	Since $\dom(M(b))=(M(b)+X_0)\cap X_1$, by \cite[Proposition A.1.1]{BoZh18} we have $M(b)+X_0=\dom(M(b))\oplus X_0$. Since for each $b\in B$, the space $M(b)+X_0(b)$ is closed and $M(b)+X_0(b)+X_1(b)=X$, by \cite[Proposition A.3.13]{BoZh18}, the families $\{M(b)+X_0(b)\}_{b\in B}$ and $\{\dom(M(b))\}_{b\in B}$ are continuous families of closed linear subspaces of $X$. By Lemma \ref{l:closed-domain}, $M(b)$ is bounded for each $b\in B$. By Lemma \ref{l:operator-space}, the families
	\begin{align*}
		\{A(b)\dom(M(b))\}_{b\in B},\quad \{A_{11}(b)\dom(M(b))\}_{b\in B}
	\end{align*}
	are continuous families of closed linear subspaces of $Y$.
	
	By \eqref{e:A-form-XY}, \eqref{e:N0b} and \eqref{e:Nb}, we have 
	\begin{align*}		
		\dom (N(b))&=A_{11}(b)\dom(M(b)),\\
		N(b)0&=N(b)\cap Y_0=N_0(b).
	\end{align*}	
    Then the families 
    \begin{align*}
    	\{\dom (N(b))\}_{b\in B},\quad \{N(b)0\}_{b\in B}
    \end{align*}
    are continuous families of closed linear subspaces of $X$ and $Y$ respectively.	
	\newline 3. The estimate for $\|N(b)\|$ with $b\in B$ and the closeness of $N(b)$.
	
	If $\dom(M(b))=\{0\}$, we have $\|N(b)\|=\|M(b)\|=0$. If $\dom(M(b))\ne\{0\}$, we have 
	\begin{align*}
		\|N(b)\|&\le\|A(b)M(b)\|=\sup\limits_{x_1\in\dom(M)\setminus\{0\}}\inf\limits_{x_0\in M(b)x_1}
		\frac{\|A_{01}(b)x_1+A_{00}(b)x_0\|}{\|A_{11}(b)x_1\|}\\
		&\le
		\|A_{01}(b)A_{11}(b)^{-1}\|+\|A_{00}(b)A_{11}(b)^{-1}\|\|M(b)\|<+\infty.
	\end{align*}
    By Lemma \ref{l:criteria-bounded-lr}, the space $N(b)$ is closed in $Y$.
    \newline 4. Continuity of the family $\{N(b)\}_{b\in B}$. 
    
    For each $b\in B$, we define two linear operators $T(b)\colon X\times Y\to Y$ and $S(b)\colon X\times Y\supset\dom(S(b))\to Y$ by
    \begin{align*}
    	T(b)(x,y)\ :=&A(b)x+y,\quad x\in X,y\in Y,\\
    	S(b)\ :=&T(b)|_{M(b)\times N_0(b)}.    	
    \end{align*}
    Then $\{T(b)\in\Bb(X\times Y,Y)\}_{b\in B}$ is a continuous family, and $\image S(b)=N(b)$.
    
    Since there hold $\Gr(T(b))+M(b)\times N_0(b)\times Y=X\times Y\times Y$ and $\Gr(S(b))=\Gr(T(b))\cap (M(b)\times N_0(b)\times Y)$, by \cite[Proposition A.3.13]{BoZh18}, $\{\Gr(S(b))\}_{b\in B}$ is a continuous family. Note that by \eqref{e:N0b} we have
    \begin{align*}
    	\ker S(b)=&\{(x,y);\;A(b)x+y=0,x\in M(b),y\in N_0(b)\}\\
    	=&\{(x_0,y_0);\;A_{00}(b)x_0+y_0=0,x_0\in X_0,y_0\in N_0(b)\}\\
    	=&\Gr(-A_{00}(b)).
    \end{align*} 
    By \cite[Proposition A.3.13]{BoZh18}, $\{N(b)\}_{b\in B}$ is a continuous family. 
\end{proof}